\documentclass[11pt]{amsart}
\usepackage[margin=28mm]{geometry}
\usepackage{amsmath,amssymb,mathtools,mathrsfs}
\usepackage[T1]{fontenc}
\usepackage{lmodern,microtype}
\usepackage{needspace,etoolbox}
\usepackage{array,booktabs}
\usepackage[hidelinks]{hyperref}
\usepackage{cite}
\numberwithin{equation}{section}
\newtheorem{theorem}{Theorem}[section]
\newtheorem{lemma}[theorem]{Lemma}

\newtheorem{corollary}[theorem]{Corollary}
\theoremstyle{remark}
\newtheorem{remark}[theorem]{Remark}
\BeforeBeginEnvironment{lemma}{\Needspace{7\baselineskip}}
\newcommand{\R}{\mathbb R}
\newcommand{\1}{\mathbf 1}

\newcommand{\supp}{\operatorname{supp}}

\newcommand{\dd}{\,d}

\title[Uniqueness and nondegeneracy for 1D fractional equations]{Uniqueness and nondegeneracy of positive solutions for a class of fractional Schrödinger equations in $\mathbb R$}
\date{}

\author{Yongsheng Jiang}
\address{(Y. Jiang)School of Statistics and Mathematics, Zhongnan University of Economics and Law, Wuhan 430073, P.R. China}
\email{jiangys@zuel.edu.cn}
\author{Xiaoyu Zeng}
\address{(X.Zeng)Center for Mathematical Sciences, Wuhan University of Technology, Wuhan 430070, P.R. China}
\email{xyzeng@whut.edu.cn}
\author{Huansong Zhou}
\address{(H.S. Zhou)Center for Mathematical Sciences, Wuhan University of Technology, Wuhan 430070, P.R. China}
\email{hszhou@whut.edu.cn}
\begin{document}
\begin{abstract}
We study positive solutions $u\in H^s(\mathbb R)$ of the fractional
Schr\"odinger equation $(-\Delta)^su+\lambda u=f(u)$ in $\mathbb R$, where
$0<s<1$, $\lambda>0$, and $f$ is a general nonlinearity such that
$f(t)/t$ is increasing. We prove that, under a single concavity condition
on an auxiliary function, positive solutions are unique up to
translations. The result applies to all positive solutions rather than
only to ground states, and it allows $f'$ to be unbounded below at zero.
It covers pure powers, sums of powers, logarithmic and saturable
nonlinearities, for suitable ranges of the parameters. We further prove that every positive solution is
nondegenerate, with the kernel of the linearized operator spanned by $u'$,
and that it has Morse index one. If moreover $f(\sqrt h)/\sqrt h$ is
concave, then for $1/2\le s<1$ a positive solution with prescribed $L^2$
norm is unique together with its Lagrange multiplier. The proofs are
purely variational, using neither a continuation argument in $s$ nor the
nondegeneracy of the linearized operator. Our method rests
on the convexity of the action along an interpolation between symmetric
decreasing functions that preserves a suitable integral quantity, and
this convexity allows us to identify every positive solution as a
constrained minimizer. For $1/2\le s<1$ we use the height interpolation
of Delgadino, Yan and Yao, while for $0<s<1/2$ we introduce a new
two-level interpolation.
\end{abstract}
\maketitle
\tableofcontents
\section{Introduction}\label{sec:intro}

We study positive solutions of
\begin{equation}\label{eq:main}
 (-\Delta)^su+\lambda u=f(u)\quad\text{in }\R,
 \qquad u\in H^s(\R),\qquad u>0\ \text{a.e.\ in }\R,
\end{equation}
where $0<s<1$, $\lambda>0$, and $(-\Delta)^s$ is the fractional
Laplacian with Fourier symbol $|\omega|^{2s}$. A solution of
\eqref{eq:main} gives the standing wave $\psi(t,x)=e^{i\lambda t}u(x)$
of the fractional nonlinear Schr\"odinger equation
$i\partial_t\psi=(-\Delta)^s\psi-f(|\psi|)\psi/|\psi|$, which arises in
the fractional quantum mechanics of Laskin \cite{Las00,Las02}. When
$f(t)=t^p$, solutions of \eqref{eq:main} also give solitary waves of
dispersion generalized Benjamin--Ono equations \cite{Wei87,KMR}, and
the case $s=1/2$, $f(t)=t^2$ is the Benjamin--Ono equation for internal
waves in deep water \cite{Ben,Ono}. Uniqueness and nondegeneracy of
positive solutions are basic ingredients in the stability theory of
these waves \cite{GSS,Wei87}, in the analysis of blow-up \cite{KMR},
and in the construction of concentrating solutions by
Lyapunov--Schmidt reduction \cite{DdPW}.

For $s=1$ the uniqueness of positive solutions is classical. In one
dimension it follows, under suitable assumptions on $f$, from the first
integral $|u'|^2=\lambda u^2-2F(u)$, where $F$ is the primitive of $f$
with $F(0)=0$. In higher dimensions it holds for pure powers and for classes of
nonlinearities, by the radial symmetry of positive solutions \cite{GNN}
and the ODE methods of Coffman \cite{Cof}, Kwong \cite{Kw}, McLeod
\cite{McL} and Serrin--Tang \cite{ST}. For $0<s<1$ these tools are not
available, since $(-\Delta)^s$ is nonlocal and \eqref{eq:main} has no
first integral, even in one dimension.

For $0<s<1$, Amick and Toland \cite{AT} proved uniqueness in the
Benjamin--Ono case $s=1/2$, $f(t)=t^2$ by complex-analytic methods. For
$f(t)=t^p$, Frank and Lenzmann \cite{FL} proved that the ground state of
\eqref{eq:main} is unique up to translations and nondegenerate for all
$0<s<1$ and all $H^s$-subcritical $p$. Fall and Valdinoci \cite{FV}
treated $s$ close to $1$ in $\R^N$, and Frank, Lenzmann and Silvestre
\cite{FLS} extended \cite{FL} to all dimensions. Here a ground state is a
positive minimizer of the associated Weinstein functional (see
\cite{FL,FLS} for the precise definitions). The proofs in \cite{FL,FLS}
first establish the nondegeneracy of the linearized operator via the
Caffarelli--Silvestre extension \cite{CS}, and then connect the ground
state to the case $s=1$ by a continuation argument in $s$. Recently, Fall
and Weth \cite{FW} proved that for $N=1$ and $p=2$ there is exactly one
even positive solution for every $1/6<s<1$. Since positive solutions are
even up to translations \cite[Proposition~3.1]{FLS}, this gives uniqueness
up to translations and extends the Amick--Toland theorem to all these $s$
without complex analysis. Apart from the case $f(t)=t^2$ treated in \cite{AT,FW}, all these results
concern ground states. For $f(t)=t^p$, Frank \cite{Fr} asked whether every
positive solution is a ground state. For
$p\ne2$ the uniqueness of all positive solutions, and not only of ground
states, remains open for $0<s<1$ (see the discussion in \cite{FW}).

The second open question concerns the nonlinearity. The arguments of
\cite{FL,FLS,FW} rely on the special structure of pure powers, and the
extension to other nonlinearities is described in \cite{FLS} as a
challenging open problem. So far, uniqueness of ground states has been obtained only in a few
cases, such as double-power nonlinearities with a small coefficient
\cite{BDG} and the logarithmic nonlinearity \cite{APZ}, both by
approximation by pure powers. The logarithmic nonlinearity, which goes
back to \cite{BBM}, is the model case of nonlinearities with $f'(t)\to-\infty$ as $t\to0^+$, for which the
potential $\lambda-f'(u)$ of the linearized operator is unbounded, and this
is the main difficulty in \cite{APZ}. For such nonlinearities the
uniqueness of all positive solutions
remains open for $0<s<1$.

The purpose of this paper is to address both questions in one
dimension. We prove that, for every $0<s<1$ and for a class of nonlinearities, all positive solutions of \eqref{eq:main} are
unique up to translations and nondegenerate, without assuming that they
are ground states. The method is variational and relies neither on the
special structure of pure powers nor on a continuation argument in $s$.

\medskip
\noindent\textit{Hypotheses.}
Throughout, $0<s<1$ and $\lambda>0$. We assume that
\begin{equation}\label{eq:f-regularity}\tag{f1}
 f\in C([0,\infty))\cap C^\infty((0,\infty)),\qquad f(0)=0,
\end{equation}
and that
\begin{equation}\label{eq:f-assumptions}\tag{f2}
 \limsup_{t\to0^+}f'(t)<\lambda,\qquad
 tf'(t)>f(t)\quad\text{for }t>0.
\end{equation}
The second condition in \eqref{eq:f-assumptions} says that $f(t)/t$ is
strictly increasing, and the first one that $f(t)<\lambda t$ for small
$t>0$. Hence, if $f(t)>\lambda t$ for some $t>0$, there is a unique
$t_\lambda>0$ with
\begin{equation}\label{eq:t-lambda}
 f(t_\lambda)=\lambda t_\lambda,\qquad f(t)<\lambda t\ \ (0<t<t_\lambda),
 \qquad f(t)>\lambda t\ \ (t>t_\lambda).
\end{equation}
Our main hypothesis is a concavity condition. For a level $\mu\ge0$ and
an exponent $r\ge1$ we set $\beta=1-1/r$ and define the auxiliary
function
\begin{equation}\label{eq:shifted-b}
 b_{r,\mu}(h)=\frac{f(\mu+h^{1/r})-\lambda(\mu+h^{1/r})}{h^\beta},
 \qquad h>0,\qquad\text{and}\qquad b_r=b_{r,0}.
\end{equation}
We require that $b_{r,\mu}$ be concave on $(0,\infty)$ for some pair
$(\mu,r)$ in one of the following two cases, which depend on the
value of $s$.
\begin{equation}\label{eq:caseA}\tag{f3}
\begin{gathered}
 \frac12\le s<1,\qquad \mu=0,\qquad 2\le r\le\frac1{1-s},\\[2pt]
 b_r''(h)\le0\quad(h>0),\qquad\text{and}\qquad
 b_r''(h)<0\quad(h>0)\ \text{ if } r=\frac1{1-s}.
\end{gathered}
\end{equation}
\begin{equation}\label{eq:caseB}\tag{f4}
\begin{gathered}
 0<s<\frac12,\qquad \mu>0,\qquad f(\mu)\le\lambda\mu,\qquad
 1\le r\le\frac1{1-s},\\[2pt]
 b_{r,\mu}''(h)\le0\quad(h>0).
\end{gathered}
\end{equation}

\medskip
We now state our main results.

\begin{theorem}\label{thm:main}
Let $0<s<1$ and $\lambda>0$, assume \eqref{eq:f-regularity} and
\eqref{eq:f-assumptions}, and suppose that there are $\mu\ge0$ and
$r\ge1$ such that either \eqref{eq:caseA} or \eqref{eq:caseB} holds.
Then \eqref{eq:main} has at most one positive solution in $H^s(\R)$,
up to translations. More precisely, for any two such solutions
$u_1,u_2$, there exists $x_0\in\R$ such that
\begin{equation*}
 u_1(x)=u_2(x-x_0)\qquad\text{for a.e. }x\in\R.
\end{equation*}
\end{theorem}

Theorem~\ref{thm:main} concerns all positive solutions. Nothing is assumed about the value of their action, their Morse index or
their symmetry, and in this respect
it differs from the results of \cite{FL,FLS}, which concern ground
states. The only other results of this kind for $0<s<1$ are those of
\cite{AT,FW}, which concern $f(t)=t^2$. It is a pure uniqueness
statement. Existence is neither assumed nor asserted and is a separate
question \cite{FL,FQT,Ar,DSZ}. The structural conditions
\eqref{eq:f-regularity} and \eqref{eq:f-assumptions} require that
$f(t)/t$ be increasing and that $f'$ be bounded above near zero, but
impose no lower bound on $f'$. Nonlinearities with $f'(t)\to-\infty$ as
$t\to0^+$, such as $f(t)=at\log t^2$, are therefore covered by the same
argument as pure powers. The main hypothesis is the concavity of
$b_{r,\mu}$. It holds for pure powers, sums of powers,
logarithmic and saturable nonlinearities (Corollaries~\ref{cor:powers}
and~\ref{cor:concave}).

For a positive solution $u$, let $L_u$ denote the self-adjoint
realization of $(-\Delta)^s+\lambda-f'(u)$ associated with its
quadratic form on
\begin{equation}\label{eq:intro-form-domain}
 D_u=\left\{\phi\in H^s(\R):
       \int_\R(\lambda-f'(u))_+\phi^2\dd x<\infty\right\}.
\end{equation}
The operator domain and its relation to the distributional equation
are specified in Section~\ref{sec:nondegeneracy}. Since $f'$ may be
unbounded below at zero, $D_u$ may be a proper subspace of $H^s(\R)$.

\begin{theorem}\label{thm:nondegeneracy}
Assume the hypotheses of Theorem~\ref{thm:main}, and let $u$ be
a positive solution. Then $u'\in D(L_u)$ and
\begin{equation}\label{eq:full-kernel}
 \ker L_u=\operatorname{span}\{u'\}.
\end{equation}
Moreover, the quadratic form of $L_u$ has Morse index one.
\end{theorem}

Theorem~\ref{thm:nondegeneracy} holds for every positive solution, for
every $0<s<1$ and for all nonlinearities covered by
Theorem~\ref{thm:main}, including the logarithmic ones, for which the
potential $\lambda-f'(u)$ is unbounded. In particular, every positive
solution has Morse index one, the property by which ground states are
defined in \cite{FLS}. For general nonlinearities both assertions of
Theorem~\ref{thm:nondegeneracy} are new.

Theorems~\ref{thm:main} and~\ref{thm:nondegeneracy} apply to the
nonlinearities most commonly studied in the literature, as the
following two corollaries show. 

\begin{corollary}\label{cor:powers}
Let $\lambda>0$ and
\begin{equation}\label{eq:mixed-f}
 f(t)=at\log t^2+\sum_{j=1}^mc_jt^{p_j},\qquad
 a\ge0,\quad c_j>0,\quad p_j>1,\quad m\ge0,
\end{equation}
with $f(0)=0$, where $a>0$ if $m=0$. Let $p_*=\max_jp_j$ if $m\ge1$ and
$p_*=1$ if $m=0$. Then \eqref{eq:f-regularity} and
\eqref{eq:f-assumptions} hold for every $\lambda>0$, and the hypotheses
of Theorem~\ref{thm:main} hold in each of the following cases.
\begin{enumerate}
\item[(i)] $1/2\le s<1$, $a>0$ and $p_*\le3$.
\item[(ii)] $1/2\le s<1$, $a=0$, and there is $r\in[2,1/(1-s)]$ such
that $r-1\le p_j\le2r-1$ for all $j$, where $p_j<3$ for some $j$ if
$s=1/2$.
\item[(iii)] $1/3\le s<1/2$ and $p_*\le(1+s)/(1-s)$.
\item[(iv)] $1/4\le s<1/3$ and $p_*\le3/2$.
\item[(v)] $1/5\le s<1/4$ and $p_*\le5/4$.
\item[(vi)] $1/6\le s<1/5$ and $m=0$.
\end{enumerate}
Consequently, in each case positive solutions of \eqref{eq:main} are
unique up to translations, nondegenerate and of Morse index one.
\end{corollary}

\begin{remark}[Pure powers]\label{rem:pure-powers}
For $f(t)=t^p$, Corollary~\ref{cor:powers} applies whenever
\begin{equation}\label{eq:pure-power-range}
\begin{gathered}
 1<p\le\frac{1+s}{1-s}\ \ \Bigl(p<3\text{ if }s=\tfrac12\Bigr)
 \ \text{ and }\ s\ge\tfrac13,\\[2pt]
 1<p\le\tfrac32\ \text{ and }\ \tfrac14\le s<\tfrac13,\qquad
 1<p\le\tfrac54\ \text{ and }\ \tfrac15\le s<\tfrac14,
\end{gathered}
\end{equation}
and then every positive solution of \eqref{eq:main} is a translate of
the ground state of \cite{FL}. This answers the question of \cite{Fr}
in the range \eqref{eq:pure-power-range}. For $s>1/2$ this range
contains the $L^2$-critical exponent $p=1+4s$ and part of the
$L^2$-supercritical range, since $(1+s)/(1-s)>1+4s$, while for $s<1/2$
it is contained in the $L^2$-subcritical range. The case $s=1/2$,
$p=2$ is the Benjamin--Ono equation, for which uniqueness was proved
in \cite{AT,FW}.
\end{remark}

For the logarithmic nonlinearity $f(t)=at\log t^2$ the substitution
$u(x)=e^{\lambda/(2a)}Q\bigl((2a)^{1/(2s)}x\bigr)$ transforms
\eqref{eq:main} into the equation $(-\Delta)^sQ=Q\log Q$ of
\cite{APZ}. By Corollary~\ref{cor:powers}, its positive solutions are
unique up to translations and nondegenerate for every $s\ge1/6$,
without the minimality assumption of \cite{APZ}.

\begin{corollary}\label{cor:concave}
Let $\lambda>0$, and let $f$ satisfy \eqref{eq:f-regularity} and
\eqref{eq:f-assumptions}. Suppose that $f(t)>\lambda t$ for some $t>0$, so that the level
$t_\lambda$ of \eqref{eq:t-lambda} is defined, and that $f$ is concave
on $[t_\lambda,\infty)$. Then, for every
$0<s<1/2$, \eqref{eq:caseB} holds with $\mu=t_\lambda$ and $r=1$, and
the conclusions of Theorems~\ref{thm:main} and~\ref{thm:nondegeneracy}
hold.
\end{corollary}
\begin{remark}[A saturable nonlinearity]\label{rem:saturable}
The nonlinearity $f(t)=ct^3/(1+t^2)$ with $c>0$ occurs in nonlinear
optics as a saturable nonlinearity. It is bounded by a linear function,
so it is neither covered by \cite{FL,FLS} nor a perturbation of a pure
power. It is concave on $[t_\lambda,\infty)$ exactly when
$\lambda\ge3c/4$, and $f(t)>\lambda t$ for some $t>0$ exactly when
$\lambda<c$. Hence the conclusions of Theorems~\ref{thm:main}
and~\ref{thm:nondegeneracy} hold for every $0<s<1/2$ if
$3c/4\le\lambda<c$, and for every $1/2\le s<1$ by case
\eqref{eq:caseA} with $r=2$. For $\lambda\ge c$ we have
$f(t)<\lambda t$ for all $t>0$, and there is no positive solution.
\end{remark}

\medskip

Since the $L^2$ norm is conserved by the Schr\"odinger flow, we also
consider \eqref{eq:main} with prescribed $\int_\R u^2\dd x$, where
$\lambda$ becomes a Lagrange multiplier \cite{CL,Je,So,LZ}. For $s=1$ and
$N\ge2$, uniqueness at prescribed $L^2$ norm has been proved recently, for
double-power nonlinearities in the $L^2$-subcritical and $L^2$-critical
cases, by following a global branch of solutions \cite{SZ}. For $s=1$ and $N=1$, positive
solutions are unique for each $\lambda$, and uniqueness at prescribed
$L^2$ norm follows from the strict monotonicity in $\lambda$ of the
$L^2$ norm of the positive solution. This monotonicity holds for
$f(t)=c_1t^p+c_2t^q$ with $c_1,c_2>0$ and $1<p<q\le5$, see \cite{KLT}
and the references therein. Our method compares two
such solutions directly.

\begin{theorem}\label{thm:normalized}
Let $1/2\le s<1$. Assume \eqref{eq:f-regularity} and
\begin{equation}\label{eq:normalized-assumptions}
 tf'(t)>f(t)\quad(t>0),\qquad
 \limsup_{t\to0^+}f'(t)<+\infty.
\end{equation}
Define $g(h)=f(\sqrt h)/\sqrt h$ for $h>0$, and suppose that
\begin{equation}\label{eq:normalized-concavity}
 g''(h)\le0\quad(h>0),\qquad
 g''(h)<0\quad(h>0)\ \text{if }s=1/2.
\end{equation}
Let $M>0$, and let $(\lambda_1,u_1)$ and $(\lambda_2,u_2)$ be two
solutions of
\begin{equation}\label{eq:normalized-problem}
 (-\Delta)^su+\lambda u=f(u)\quad\text{in }\R,\qquad
 u\in H^s(\R),\quad u>0\ \text{a.e.},\quad \int_\R u^2\dd x=M,
\end{equation}
with multipliers
\begin{equation}\label{eq:normalized-multiplier}
 \lambda_1,\lambda_2\in\R,\qquad
 \lambda_1,\lambda_2>\limsup_{t\to0^+}f'(t).
\end{equation}
Then $\lambda_1=\lambda_2$, and there is $x_0\in\R$ such that
$u_1(x)=u_2(x-x_0)$ for almost every $x\in\R$.
\end{theorem}

\begin{remark}\label{rem:normalized}
(i) For pure powers, uniqueness at prescribed $L^2$ norm follows from
uniqueness for fixed $\lambda$ by scaling. Let $0<s<1$ and $f(t)=t^p$
with $p$ in the range \eqref{eq:pure-power-range} and $p\ne1+4s$. By
Corollary~\ref{cor:powers} and \cite{FL}, for every $\lambda>0$ the
positive solutions of \eqref{eq:main} are the translates of
$Q_\lambda(x)=\lambda^{1/(p-1)}Q(\lambda^{1/(2s)}x)$, where $Q$ is the
ground state for $\lambda=1$, and
$\|Q_\lambda\|_2^2=\lambda^{\frac2{p-1}-\frac1{2s}}\|Q\|_2^2$, where the
exponent is nonzero. Hence for every $M>0$ there is exactly one
$\lambda>0$ for which \eqref{eq:normalized-problem} has a solution, and
this solution is unique up to translations. In particular, for these
nonlinearities the conclusion of Theorem~\ref{thm:normalized} holds
also for $0<s<1/2$.

(ii) For nonlinearities which are not homogeneous, no such scaling is
available. Uniqueness for each fixed $\lambda$ then does not by itself
give uniqueness at prescribed $L^2$ norm, since two solutions with the
same $L^2$ norm may have different multipliers. The approach for $s=1$
and $N=1$ recalled above requires the strict monotonicity of
$\lambda\mapsto\|u_\lambda\|_2^2$, which, as far as we know, has not
been established for such nonlinearities when $0<s<1$.
Theorem~\ref{thm:normalized} needs neither this monotonicity nor any
information on the set of possible multipliers, since it compares two
solutions directly. It applies to \eqref{eq:mixed-f} in case (i), and in
case (ii) with $p_j\le3$ for all $j$; for instance, to
$f(t)=t^p+t^q$ with $1<p<q\le3$, and to
$f(t)=at\log t^2+\sum_jc_jt^{p_j}$ with $a>0$ and $p_j\le3$, in the
latter case for all real multipliers. We do not know whether
Theorem~\ref{thm:normalized} extends to nonhomogeneous nonlinearities
when $0<s<1/2$.
\end{remark}

\medskip
Our approach is variational and is based on the action $S$ defined in
\eqref{eq:action} below. It rests on two properties of positive
solutions. Every positive solution is a strict maximizer of $S$ along the
ray $c\mapsto cu$, and it is a strict
minimizer of $S$ among symmetric decreasing comparison functions with
the same value of $\mathcal N(v)=\int_\R(v-\mu)_+^r\dd x$. Comparing two
positive solutions at each other's value of $\mathcal N$, these two
properties give uniqueness. The minimality rests on the convexity of
$S$ along interpolations which preserve $\mathcal N$. For $1/2\le s<1$
we use the height interpolation of \cite{DYY}, together with a new
positive series for the kinetic energy of the $r$th root of a density.
For $0<s<1/2$ we introduce a two-level interpolation, which treats the
parts of a function below and above $\mu$ separately. The same
convexity gives the positivity of the second variation in even
directions tangent to the constraint, and hence the nondegeneracy in
even directions, while odd directions are treated by a comparison of
quadratic forms. No continuation argument and no a priori information
on the linearized operator are needed. The structure of the argument
is summarized in Section~\ref{sec:further}, where further applications
are discussed.

\smallskip
\noindent\textit{Organization of the paper.}
In Section~\ref{sec:prelim} we fix the notation and collect the basic
properties of positive solutions, namely boundedness, decay,
smoothness, symmetry and strict monotonicity.
Section~\ref{sec:height} introduces the height interpolation and proves
the convexity of the kinetic energy along it, and Section~\ref{sec:two-level} introduces the
two-level interpolation for $0<s<1/2$. In Section~\ref{sec:coercivity} we establish the coercivity of the
second variation of $S$ in even directions tangent to the constraint. Theorem~\ref{thm:main} is proved in
Section~\ref{sec:proof-main}, and Theorems~\ref{thm:nondegeneracy}
and~\ref{thm:normalized} are proved in Section~\ref{sec:nondegeneracy}.
Section~\ref{sec:further} discusses extensions of the method to further
problems.
Appendix~\ref{app:convexity} establishes the convexity of the monomial
functions used in Sections~\ref{sec:height} and~\ref{sec:two-level},
Appendix~\ref{app:examples} gives the proofs of
Corollaries~\ref{cor:powers} and~\ref{cor:concave}, and
Appendix~\ref{app:prelim} contains the proofs of the lemmas stated in
Section~\ref{sec:prelim}.

\section{Preliminaries}\label{sec:prelim}

%\subsection*{Notation and the action}
We use the unitary Fourier transform
\begin{equation*}
 \widehat v(\omega)=(2\pi)^{-1/2}\int_\R e^{-ix\omega}v(x)\dd x,
\end{equation*}
and for $\sigma\ge0$ we let $H^\sigma(\R)$ be the space of $v\in L^2(\R)$
with $|\omega|^\sigma\widehat v\in L^2(\R)$, see \cite{DNPV} for
background. A positive solution of \eqref{eq:main}
is a function $u\in H^s(\R)$, positive almost everywhere, such that
$f(u)\in L^1_{\rm loc}(\R)$ and $(-\Delta)^su+\lambda u=f(u)$ holds in
the sense of distributions. The
kinetic energy and its singular integral representation are
\begin{equation*}
 T_s(v)=\int_\R|\omega|^{2s}|\widehat v(\omega)|^2\dd\omega
 =\frac{c_{1,s}}2\iint_{\R^2}
       \frac{(v(x)-v(y))^2}{|x-y|^{1+2s}}\dd x\dd y,\qquad
 c_{1,s}=\frac{4^s\Gamma(1/2+s)}{\sqrt\pi\,|\Gamma(-s)|},
\end{equation*}
$T_s(v,w)$ denotes the corresponding bilinear form, and
$\|v\|_{H^s}^2=T_s(v)+\|v\|_2^2$. All functions are real-valued, $|E|$
is the Lebesgue measure of a measurable set $E\subset\R$ and $\1_E$ its
indicator function, $t_+=\max\{t,0\}$, $C_c^\infty(\R)$ is the space of
smooth functions with compact support, and $C$ denotes a positive
constant whose value may change from line to line. We call a function
\emph{symmetric decreasing} if it is even and nonincreasing on
$(0,\infty)$. Finally, $\beta=1-1/r$ as in \eqref{eq:shifted-b}, so
that $1/2\le\beta\le s$ in case \eqref{eq:caseA} and $0\le\beta\le s$
in case \eqref{eq:caseB}.

Formally, \eqref{eq:main} is the Euler--Lagrange equation of the action
\begin{equation}\label{eq:action}
 S(v)=\frac12T_s(v)+\frac\lambda2\|v\|_2^2-\int_\R F(v)\dd x,
 \qquad F(t)=\int_0^tf(\tau)\dd\tau,
\end{equation}
defined for nonnegative $v\in H^s(\R)$ with $F(v)\in L^1(\R)$. Since
$S$ need not be differentiable on $H^s(\R)$ when $f'$ is unbounded
below at zero, we never use derivatives of $S$ on $H^s(\R)$. Instead we
use the weak formulation \eqref{eq:weak-all-tests} of
Lemma~\ref{lem:solution-estimates}, and we differentiate $S$ only along
explicit curves. The quantity preserved by the interpolations of
Sections~\ref{sec:height} and~\ref{sec:two-level} is
$\mathcal N(v)=\int_\R(v-\mu)_+^r\dd x$.

%\subsection*{Properties of positive solutions}
We show that every positive solution is bounded, smooth and strictly
positive, decays at infinity, and after a translation is even and
strictly decreasing on $(0,\infty)$. For locally Lipschitz
nonlinearities such properties are well known, see \cite{FQT,FLS} for
regularity, decay and symmetry, and \cite{CLL} for the method of moving
planes for the fractional Laplacian. Since $f'$ may tend to $-\infty$ at
zero, these results cannot be quoted directly, and we only use an upper
bound for $f'$ on bounded subsets of $(0,\infty)$. The proofs, which
follow the methods of \cite{CLL,FLS,RS}, are collected in
Appendix~\ref{app:prelim}.

\begin{lemma}\label{lem:growth}
Assume \eqref{eq:f-regularity}, \eqref{eq:f-assumptions} and \eqref{eq:caseA}.
There is $C>0$ such that
\begin{equation*}
 |f(t)|+t|f'(t)|\le Ct^{2r-1}\qquad(t\ge1).
\end{equation*}
In case \eqref{eq:caseB}, there is $C>0$ such that
\begin{equation*}
 f(t)_+\le C(t+t^{2r-1})\quad(t\ge0),\qquad
 1\le2r-1\le\frac{1+s}{1-s}\le\frac1{1-2s}.
\end{equation*}
\end{lemma}

Fix $m_*>0$, an upper bound for the values of the functions under
consideration. By \eqref{eq:f-regularity} and the
first condition in \eqref{eq:f-assumptions}, we may choose
\begin{equation}\label{eq:aux-definition}
 K>\max\{\lambda,\sup_{0<t\le m_*}f'(t)\},\qquad
 f_K(t)=Kt-f(t),\qquad F_K(t)=\frac K2t^2-F(t)\quad(0\le t\le m_*).
\end{equation}
The constant $K$ may be increased when $m_*$ is increased.

\begin{lemma}\label{lem:scalar}
Assume \eqref{eq:f-regularity} and \eqref{eq:f-assumptions}.
Let $m_*>0$ and choose $K,f_K,F_K$ as in \eqref{eq:aux-definition}.
Then
\begin{equation}\label{eq:scalar-bounds}
 0<tf_K'(t)<f_K(t),\qquad
 \frac12tf_K(t)\le F_K(t)\le tf_K(t)\quad(0<t\le m_*).
\end{equation}
If $c>0$ and $t,ct\in[0,m_*]$, then
\begin{equation}\label{eq:scalar-scaling}
 f_K(ct)\le\max\{1,c\}f_K(t),\qquad
 F_K(ct)\le c\max\{1,c\}\,tf_K(t).
\end{equation}
\end{lemma}

For $0<s<1/2$ we shall only use the growth condition
\begin{equation}\label{eq:small-growth}
 0<s<\frac12,\qquad\lambda>0,\qquad
 f(t)_+\le C_1(t+t^{p_0})\quad(t\ge0),
\end{equation}
where $C_1>0$ and $1\le p_0\le1/(1-2s)$; in case \eqref{eq:caseB} it
holds with $p_0=2r-1$ by Lemma~\ref{lem:growth}. Thus
Lemmas~\ref{lem:solution-estimates} and~\ref{lem:symmetry} apply, for
instance, to $f(t)=t^2$ for $1/4\le s<1/2$. The next lemma provides the
two tools which replace the embedding into $L^\infty(\R)$ when
$s<1/2$: the kernel of $((-\Delta)^s+\lambda)^{-1}$ and a barrier.

\begin{lemma}\label{lem:kernel}
Let $0<s<1/2$ and $\lambda>0$.
\begin{enumerate}
\item[(i)] There is a function $G_\lambda\in L^1(\R)$ with
$\sqrt{2\pi}\,\widehat{G_\lambda}(\omega)=(|\omega|^{2s}+\lambda)^{-1}$,
which is positive, symmetric decreasing, has $\int_\R G_\lambda\dd x=1/\lambda$, and
satisfies $G_\lambda(x)\le C|x|^{2s-1}$ for $0<|x|\le1$ and
$G_\lambda(x)\le C|x|^{-1-2s}$ for $|x|\ge1$. Consequently,
$G_\lambda\in L^q(\R)$ for $1\le q<1/(1-2s)$, and
\begin{equation}\label{eq:kernel-modulus}
 \int_\R|G_\lambda(y+h)-G_\lambda(y)|\dd y\le C|h|^{2s}
 \qquad(|h|\le1).
\end{equation}
\item[(ii)] Let $\Psi_\sigma(x)=(\sigma^2+x^2)^{-(1+2s)/2}$ for $\sigma>0$.
There is $C_\Psi>0$ such that
$(-\Delta)^s\Psi_\sigma\ge-C_\Psi\sigma^{-2s}\Psi_\sigma$ on $\R$ for
every $\sigma>0$.
\end{enumerate}
\end{lemma}

\begin{lemma}\label{lem:solution-estimates}
Assume \eqref{eq:f-regularity}, \eqref{eq:f-assumptions} and \eqref{eq:caseA}, or assume
\eqref{eq:f-regularity}, \eqref{eq:f-assumptions} and
\eqref{eq:small-growth}, and let $u$ be a positive solution of \eqref{eq:main}. Then the following hold.
\begin{enumerate}
\item[(i)] $u\in L^\infty(\R)$.
\item[(ii)] For any $m_*\ge\|u\|_\infty$ and $K,f_K,F_K$ chosen as in
\eqref{eq:aux-definition},
\begin{equation}\label{eq:a-L2}
 \|f_K(u)\|_2\le(K-\lambda)\|u\|_2.
\end{equation}
Moreover, $f(u),uf'(u)\in L^2(\R)$, $F(u),uf(u),u^2f'(u)\in L^1(\R)$, and
\begin{equation}\label{eq:weak-all-tests}
 T_s(u,\psi)+\int_\R(\lambda u-f(u))\psi\dd x=0
 \qquad\text{for every }\psi\in H^s(\R).
\end{equation}
\item[(iii)] Under \eqref{eq:caseA},
$u\in H^{2s}(\R)\cap C_0(\R)$ and $u'\in L^2(\R)$,
where $C_0(\R)$ denotes the continuous functions vanishing at infinity.
Under \eqref{eq:small-growth}, $u=G_\lambda*f(u)$, and the continuous
representative of $u$ belongs to $C_0(\R)$ and satisfies
$u(x)\le C(1+|x|)^{-1-2s}$ for $x\in\R$; in particular,
$u\in L^1(\R)$.
\item[(iv)] The continuous representative is strictly positive at every
point and belongs to $C^\infty(\R)$.
\item[(v)] $f(t)>\lambda t$ for some $t\in(0,\|u\|_\infty]$. In
particular, $\|u\|_\infty>\mu$ whenever $\mu>0$ and $f(\mu)\le\lambda\mu$.
\end{enumerate}
\end{lemma}

\begin{lemma}\label{lem:symmetry}
Assume \eqref{eq:f-regularity}, \eqref{eq:f-assumptions} and \eqref{eq:caseA}, or assume
\eqref{eq:f-regularity}, \eqref{eq:f-assumptions} and
\eqref{eq:small-growth}, and let $u$ be a positive solution of \eqref{eq:main}. After a translation,
the following hold.
\begin{enumerate}
\item[(i)] $u$ is even and strictly decreasing on $(0,\infty)$.
\item[(ii)] $u'(x)<0$ for every $x>0$.
\item[(iii)] $u''(0)<0$.
\end{enumerate}
\end{lemma}

\section{Height interpolation and convexity of the kinetic energy} \label{sec:height}

The interpolation has to make the potential energy
$\int_\R(\frac\lambda2u^2-F(u))\dd x$ convex. This
fails for the linear interpolation of $u^2$, which gives uniqueness for
sublinear problems \cite{BF,BO,DSa}, and we use instead the height
interpolation of Delgadino, Yan and Yao \cite{DYY}, see also
\cite{DV}. Along it the concavity of $b_{r,\mu}$ becomes the
convexity of the potential energy, since
$\frac\lambda2u^2-F(u)=-G_{r,\mu}\bigl((u-\mu)_+^r\bigr)$ where
$u\ge\mu$, with $G_{r,\mu}(h)=F(\mu+h^{1/r})-\frac\lambda2(\mu+h^{1/r})^2$
and $G_{r,\mu}'=b_{r,\mu}/r$.

The height function and its linear interpolation are due to
\cite{DYY}. Lemmas~\ref{lem:height-coordinates}--\ref{lem:layer}
collect the properties of heights used below. Some of them are
contained in \cite{DYY}, in a different form. We need them for heights
which are only assumed to be convex, for densities with compact
support, which occur in the two-level interpolation of
Section~\ref{sec:two-level}, and for integrals involving the $r$th root
$\rho_H^{1/r}$ rather than $\rho_H$ itself, and we therefore include
the proofs.

The new ingredient of this section is the series for
$T_s(\rho_H^{1/r})$ in Lemma~\ref{lem:kinetic-series}(iii), and the
resulting convexity of the kinetic energy along the height
interpolation in Lemma~\ref{lem:kinetic-convexity}.

Throughout this section, $0<s<1$, $r\ge1$, $\beta=1-1/r$ and $M>0$. The
\emph{mass} of a density $\rho$ is $\int_\R\rho\dd x$. For
$\rho=(u-\mu)_+^r$ it equals $\mathcal N(u)$, and for $r=2$ and $\mu=0$
it is the squared $L^2$ norm in Theorem~\ref{thm:normalized}. The
variables $\xi,\zeta\in(0,M)$ are mass coordinates, and for a convex
function $H$ on $[0,M]$, $H'$ denotes the right derivative.

\begin{lemma}\label{lem:height-coordinates}
\begin{enumerate}
\item[(i)] Let $\rho$ be a nonnegative symmetric decreasing function in
$L^1(\R)\cap L^\infty(\R)$ with $\int_\R\rho\dd x=M$. Then the function
$G(h)=\int_\R\min\{\rho,h\}\dd x$ is continuous, concave and strictly
increasing on $[0,\|\rho\|_\infty]$, with $G(0)=0$ and $G(\|\rho\|_\infty)=M$. Its
inverse $H\colon[0,M]\to[0,\|\rho\|_\infty]$ (the \emph{height function} of $\rho$) is
convex, and
\begin{equation}\label{eq:height-derivative-general}
 H'(\xi)=\frac1{|\{\rho>H(\xi)\}|}\qquad(0\le\xi<M),
\end{equation}
which is positive for $0<\xi<M$ and equals $1/|\{\rho>0\}|\in[0,\infty)$
at $\xi=0$, with $1/\infty=0$.
\item[(ii)] Let $H:[0,M]\to[0,\infty)$ be continuous and convex, with
$H(0)=0$ and $H'>0$ on $(0,M)$, and define
\begin{equation}\label{eq:reconstruction}
 R_H(\xi)=\frac1{2H'(\xi)}\quad(0<\xi<M),\qquad
 \rho_H(x)=\int_0^M\1_{\{|x|<R_H(\xi)\}}H'(\xi)\dd\xi,
 \qquad u_H=\rho_H^{1/r}.
\end{equation}
Then $\rho_H$ is symmetric decreasing, $0\le\rho_H\le H(M)$,
\begin{equation}\label{eq:superlevels}
 \{\rho_H>H(\xi)\}=(-R_H(\xi),R_H(\xi))\qquad(0<\xi<M),
\end{equation}
$\int_\R\rho_H\dd x=M$, and $H$ is the height function of $\rho_H$. Moreover, for
every $x\in\R$,
\begin{equation}\label{eq:monotone-root}
 u_H(x)=\frac1r\int_0^M\1_{\{|x|<R_H(\xi)\}}
 \frac{H'(\xi)}{H(\xi)^\beta}\dd\xi,
\end{equation}
and every measurable $\psi\ge0$ on $(0,\infty)$ satisfies
\begin{equation}\label{eq:monotone-cake}
 \int_0^\infty\psi(h)\,|\{\rho_H>h\}|\dd h=\int_0^M\psi(H(\xi))\dd\xi
 \quad\text{in }[0,\infty],
\end{equation}
and so does every measurable $\psi$ for which one side is absolutely
convergent.
\item[(iii)] In the situation of (i), $\rho_H=\rho$ at every point where
$\rho$ is continuous, hence almost everywhere. In particular, a symmetric
decreasing density in $L^1\cap L^\infty$ is determined almost everywhere by its height function.
\end{enumerate}
\end{lemma}
\begin{proof}
(i) By the layer cake formula,
\begin{equation*}
 G(h)=\int_\R\int_0^h\1_{\{t<\rho(y)\}}\dd t\dd y
 =\int_0^h|\{\rho>t\}|\dd t\qquad(0\le h\le\|\rho\|_\infty).
\end{equation*}
The function $t\mapsto|\{\rho>t\}|$ is nonincreasing and right
continuous, finite for $t>0$ since $|\{\rho>t\}|\le M/t$ by Chebyshev's
inequality, and positive for $t<\|\rho\|_\infty$. Hence $G$ is finite, concave
and strictly increasing, $G(0)=0$, $G(\|\rho\|_\infty)=\int_\R\rho\dd x=M$, and
the right derivative of $G$ is $|\{\rho>h\}|$ for $0<h<\|\rho\|_\infty$ and
$|\{\rho>0\}|\in(0,\infty]$ at $h=0$. The inverse of a strictly
increasing concave function is convex, and for the right derivatives,
with $\xi=G(h)$ and $\eta=G(k)$, where $k\to h^+$ corresponds to
$\eta\to\xi^+$ because $G$ is continuous and strictly increasing,
\begin{equation*}
 H'(\xi)=\lim_{\eta\to\xi^+}\frac{H(\eta)-H(\xi)}{\eta-\xi}
 =\lim_{k\to h^+}\frac{k-h}{G(k)-G(h)}=\frac1{|\{\rho>h\}|}
 =\frac1{|\{\rho>H(\xi)\}|},
\end{equation*}
which is \eqref{eq:height-derivative-general}.

(ii) Since $H$ is convex, $H'$ is nondecreasing and right continuous, so
$R_H$ is nonincreasing and right continuous with values in $(0,\infty)$
on $(0,M)$; and $H$ is absolutely continuous with $H(\xi)=\int_0^\xi H'$,
indeed $H'$ is finite, nonnegative (as $H\ge0=H(0)$) and nondecreasing
on $[0,M)$, so $H(\xi)-H(\varepsilon)=\int_\varepsilon^\xi H'$ for
$0<\varepsilon<\xi<M$, and the continuity of $H$ at $0$ and $M$ gives the
formula on $[0,M]$ with $\int_0^MH'=H(M)<\infty$.
For $x\in\R$ let
\begin{equation*}
 \xi(x)=\sup\{\zeta\in(0,M):R_H(\zeta)>|x|\},\qquad\sup\emptyset=0.
\end{equation*}
By monotonicity of $R_H$, the set $\{\zeta\in(0,M):R_H(\zeta)>|x|\}$ is an
initial segment of $(0,M)$, and by right continuity it does not contain
its endpoint; thus it equals $(0,\xi(x))$, and
\begin{equation}\label{eq:rho-H-formula}
 \rho_H(x)=\int_0^{\xi(x)}H'(\zeta)\dd\zeta=H(\xi(x)).
\end{equation}
As $\xi(x)$ is nonincreasing in $|x|$ and $H$ is strictly increasing,
$\rho_H$ is symmetric decreasing with $0\le\rho_H\le H(M)$, and for $0<\xi_0<M$,
\begin{equation*}
 \rho_H(x)>H(\xi_0)\iff\xi(x)>\xi_0\iff\xi_0\in(0,\xi(x))\iff R_H(\xi_0)>|x|,
\end{equation*}
which is \eqref{eq:superlevels}; in particular
$|\{\rho_H>H(\xi)\}|=2R_H(\xi)=1/H'(\xi)$ for $0<\xi<M$, and
$\{\rho_H>h\}=\emptyset$ for $h\ge H(M)$. For nonnegative measurable
$\psi$, the change of variables $h=H(\xi)$, which is valid because $H$
is absolutely continuous and increasing, gives
\begin{equation*}
 \int_0^\infty\psi(h)\,|\{\rho_H>h\}|\dd h
 =\int_0^M\psi(H(\xi))\,|\{\rho_H>H(\xi)\}|\,H'(\xi)\dd\xi
 =\int_0^M\psi(H(\xi))\dd\xi,
\end{equation*}
which is \eqref{eq:monotone-cake}; the signed version follows by
splitting $\psi$ into its positive and negative parts. With
$\psi=\1_{(0,H(\xi_0))}$ and the layer cake formula,
\begin{equation*}
 \int_\R\min\{\rho_H,H(\xi_0)\}\dd x
 =\int_0^{H(\xi_0)}|\{\rho_H>h\}|\dd h
 =\int_0^M\1_{\{H(\xi)<H(\xi_0)\}}\dd\xi=\xi_0\qquad(0<\xi_0\le M).
\end{equation*}
Moreover $\|\rho_H\|_\infty=H(M)$, since $\rho_H\le H(M)$ and
$|\{\rho_H>H(\xi)\}|=1/H'(\xi)>0$ for $0<\xi<M$. Hence $H$ is the height function
of $\rho_H$, and $\xi_0=M$ gives $\int_\R\rho_H\dd x=M$. Finally, the same change of variables with
$g(h)=\frac1rh^{-\beta}$ gives, for every $x$,
\begin{equation*}
 u_H(x)=H(\xi(x))^{1/r}=\int_0^{H(\xi(x))}\frac1rh^{-\beta}\dd h
 =\int_0^{\xi(x)}\frac{H'(\zeta)}{rH(\zeta)^\beta}\dd\zeta
 =\frac1r\int_0^M\1_{\{|x|<R_H(\zeta)\}}\frac{H'(\zeta)}{H(\zeta)^\beta}\dd\zeta,
\end{equation*}
which is \eqref{eq:monotone-root}.

(iii) Let $\rho$ be continuous at $x$ and $0\le h<\|\rho\|_\infty$. If $h<\rho(x)$,
then $\rho>h$ on a neighborhood of $x$, so $|\{\rho>h\}|>2|x|$ by
monotonicity; if $h\ge\rho(x)$, then $\{\rho>h\}\subset[-|x|,|x|]$. By
\eqref{eq:height-derivative-general}, $R_H(\xi)=|\{\rho>H(\xi)\}|/2$,
hence
\begin{equation*}
 \{\zeta\in(0,M):R_H(\zeta)>|x|\}=\{\zeta\in(0,M):H(\zeta)<\rho(x)\},
\end{equation*}
so that $\xi(x)=G(\rho(x))$ and \eqref{eq:rho-H-formula} gives
$\rho_H(x)=\rho(x)$. A monotone function has at most countably many
discontinuities.
\end{proof}

For the densities $\rho=u^r$ of positive smooth solutions the height function is
smooth, and we shall need its derivatives explicitly. We call
$H\in C([0,M])\cap C^2((0,M))$ \emph{admissible} if
\begin{equation}\label{eq:admissible}
 H(0)=0,\qquad H'>0,\qquad H''>0,\qquad
 \lim_{\xi\to0^+}H'(\xi)=0,\text{ and } \lim_{\xi\rightarrow M^-}H'(\xi)=\infty.
\end{equation}

\begin{lemma}\label{lem:reconstruction}
(i) Let $0<R_*\le\infty$, and let $\rho\in C(\R)\cap L^1(\R)$ be even,
positive and smooth on $(-R_*,R_*)$, with $\rho=0$ outside $(-R_*,R_*)$,
$\rho'(x)<0$ for $0<x<R_*$ and $\int_\R\rho\dd x=M$. Then its height function $H$
satisfies $H\in C^\infty((0,M))$, $H(0)=0$, $H(M)=\rho(0)$,
$\lim_{\xi\to M^-}H'(\xi)=\infty$ and $\lim_{\xi\to0^+}H'(\xi)=1/(2R_*)$,
with $1/\infty=0$. In particular $H$ is admissible if $R_*=\infty$. The
coordinate $\xi(x)=G(\rho(x))$ satisfies, for $0<x<R_*$,
\begin{equation}\label{eq:mass-coordinate}
 \xi(x)=2x\rho(x)+2\int_x^{R_*}\rho(y)\dd y,\qquad
 H(\xi(x))=\rho(x),\qquad \xi'(x)=2x\rho'(x)<0,
\end{equation}
\begin{equation}\label{eq:height-derivatives}
 H'(\xi(x))=\frac1{2x},\qquad
 H''(\xi(x))=-\frac1{4x^3\rho'(x)}>0.
\end{equation}
(ii) Conversely, let $H$ be admissible. Then $R_H:(0,M)\to(0,\infty)$ is a
strictly decreasing bijection, and $\rho_H$ is positive, continuous and
strictly decreasing on $(0,\infty)$, with $\rho_H(0)=H(M)$,
$\rho_H(R_H(\xi))=H(\xi)$ and $\rho_H'(R_H(\xi))=-2H'(\xi)^3/H''(\xi)<0$ for
$0<\xi<M$.
\end{lemma}
\begin{proof}
(i) We first show that $x\rho(x)\to0$ as $x\to R_*^-$. If $R_*=\infty$,
monotonicity and integrability give
\begin{equation*}
 0\le\frac x2\rho(x)\le\int_{x/2}^x\rho(y)\dd y\le\int_{x/2}^\infty\rho(y)\dd y
 \longrightarrow0\quad\text{as }x\to\infty,
\end{equation*}
and if $R_*<\infty$, then $\rho(x)\to\rho(R_*)=0$ by continuity. For
$0<x<R_*$ the strict monotonicity gives $\{\rho>\rho(x)\}=(-x,x)$, and
splitting the integral defining $G(\rho(x))$ at $\pm x$,
\begin{equation*}
 \xi(x)=\int_{-x}^x\rho(x)\dd y+2\int_x^{R_*}\rho(y)\dd y
 =2x\rho(x)+2\int_x^{R_*}\rho(y)\dd y.
\end{equation*}
Hence, $\xi'(x)=2\rho(x)+2x\rho'(x)-2\rho(x)=2x\rho'(x)<0$, and
$x\rho(x)\to0$ gives $\xi(x)\to M$ as $x\to0^+$ and $\xi(x)\to0$ as
$x\to R_*^-$. Thus $\xi$ is a smooth decreasing bijection of $(0,R_*)$
onto $(0,M)$ with nonvanishing derivative, $H=\rho\circ\xi^{-1}$ is
smooth on $(0,M)$, and the chain rule gives
\begin{equation*}
 H'(\xi(x))=\frac{\rho'(x)}{\xi'(x)}=\frac1{2x},\qquad
 H''(\xi(x))=\frac1{\xi'(x)}\frac{d}{dx}\Bigl(\frac1{2x}\Bigr)
 =-\frac1{4x^3\rho'(x)}>0.
\end{equation*}
Letting $x\to R_*^-$ and $x\to0^+$ gives $H'(\xi)\to1/(2R_*)$ as
$\xi\to0^+$ and $H'(\xi)\to\infty$ as $\xi\to M^-$, and
$H(M)=\rho(0)$, $H(0)=0$.

(ii) By \eqref{eq:admissible}, $R_H'=-H''/(2H'^2)<0$, $R_H(\xi)\to\infty$
as $\xi\to0^+$ and $R_H(\xi)\to0$ as $\xi\to M^-$, so $R_H$ is a
decreasing bijection with $C^1$ inverse. In the notation of the proof
of Lemma~\ref{lem:height-coordinates}(ii), $\xi(x)=R_H^{-1}(x)$ for
$x>0$, so by \eqref{eq:rho-H-formula}
$\rho_H(x)=H(R_H^{-1}(x))$, $\rho_H(0)=H(M)$ and $\rho_H(R_H(\xi))=H(\xi)$;
hence $\rho_H$ is positive, continuous and strictly decreasing on
$(0,\infty)$, with $\rho_H(x)\to H(M)$ as $x\to0^+$ and $\rho_H(x)\to0$ as
$x\to\infty$, and
\begin{equation*}
 \rho_H'(R_H(\xi))=\frac{H'(\xi)}{R_H'(\xi)}=-\frac{2H'(\xi)^3}{H''(\xi)}<0.
 \qedhere
\end{equation*}
\end{proof}

\begin{lemma}\label{lem:height-interpolation}
Let $H_0,H_1$ satisfy the hypotheses of
Lemma~\ref{lem:height-coordinates}(ii), and for $0\le\theta\le1$ let
\begin{equation}\label{eq:interpolation}
 H_\theta=(1-\theta)H_0+\theta H_1,\qquad
 u_\theta=\rho_{H_\theta}^{1/r}.
\end{equation}
Then $H_\theta$ also satisfies the hypotheses of
Lemma~\ref{lem:height-coordinates}(ii), $u_\theta$ is symmetric decreasing,
$\int_\R u_\theta^r\dd x=M$, and
\begin{equation*}
 \|u_\theta\|_\infty^r=H_\theta(M)
   =(1-\theta)\|u_0\|_\infty^r+\theta\|u_1\|_\infty^r.
\end{equation*}
If $H_0$ and $H_1$ are admissible, then so is $H_\theta$, and $u_\theta$ is
positive, continuous and strictly decreasing on $(0,\infty)$.
\end{lemma}
\begin{proof}
The affine combination is continuous and convex on $[0,M]$, with
$H_\theta(0)=0$ and
\begin{equation*}
 H_\theta'=(1-\theta)H_0'+\theta H_1'>0\quad\text{on }(0,M).
\end{equation*}
In the admissible case also $H_\theta''=(1-\theta)H_0''+\theta H_1''>0$,
and the endpoint limits of $H_0',H_1'$ in \eqref{eq:admissible} pass to
$H_\theta'$. The remaining assertions follow from
Lemma~\ref{lem:height-coordinates}(ii) and
Lemma~\ref{lem:reconstruction}(ii): $\int_\R u_\theta^r\dd x=\int_\R\rho_{H_\theta}\dd x=M$,
and $\operatorname{ess\,sup}\rho_{H_\theta}=H_\theta(M)$ by
\eqref{eq:superlevels}, since $\rho_{H_\theta}\le H_\theta(M)$ and
$|\{\rho_{H_\theta}>H_\theta(\xi)\}|=1/H_\theta'(\xi)>0$ for every $\xi<M$; and
$H_i(M)=\|u_i\|_\infty^r$ for $i=0,1$.
\end{proof}

The next lemma is the layer cake formula in mass coordinates, see
\cite{DYY}. Nonnegative integrals are understood in
$[0,\infty]$; they are not assumed to be finite, and will be shown to be finite in Lemmas~\ref{lem:kinetic-convexity} and~\ref{lem:path-finiteness}
from the finite endpoint energies.

\begin{lemma}\label{lem:layer}
Let $H$ satisfy the hypotheses of Lemma~\ref{lem:height-coordinates}(ii)
and let $J$ be absolutely continuous on $[0,H(M)]$ with $J(0)=0$. If
$J'\ge0$ almost everywhere, or if $\int_0^M|J'(H(\xi))|\dd\xi<\infty$, then
\begin{equation}\label{eq:layer}
 \int_\R J(\rho_H(x))\dd x=\int_0^M J'(H(\xi))\dd\xi,
\end{equation}
in $[0,\infty]$ in the first case and with $J(\rho_H)\in L^1(\R)$ in
the second. In particular, for every $q>0$,
\begin{equation}\label{eq:powers}
 \int_\R u_H^q\dd x=\frac qr\int_0^M H(\xi)^{q/r-1}\dd\xi\quad\text{in }[0,\infty].
\end{equation}
\end{lemma}
\begin{proof}
Since $J(\rho_H(x))=\int_0^{\rho_H(x)}J'(h)\dd h$, it follows from 
\eqref{eq:monotone-cake}  and  $J'\ge0$ that
\begin{equation*}
 \int_\R J(\rho_H)\dd x
 =\int_\R\int_0^\infty\1_{\{h<\rho_H(x)\}}J'(h)\dd h\dd x
 =\int_0^\infty J'(h)\,|\{\rho_H>h\}|\dd h
 =\int_0^MJ'(H(\xi))\dd\xi.
\end{equation*}
In the second case the same identity with $|J'|$ in place of $J'$ shows
that the double integral converges absolutely, and Fubini's theorem
gives \eqref{eq:layer}. For \eqref{eq:powers} take $J(h)=h^{q/r}$, which
is absolutely continuous on $[0,H(M)]$ with $J'(h)=\frac qrh^{q/r-1}>0$.
\end{proof}

For a measurable real function $v$ we define $T_s(v)\in[0,\infty]$ by the
double integral
\begin{equation}\label{eq:kinetic-double}
 T_s(v)=\frac{c_{1,s}}2\iint_{\R^2}
       \frac{(v(x)-v(y))^2}{|x-y|^{1+2s}}\dd x\dd y,
\end{equation}
which agrees with the kinetic energy of Section~\ref{sec:prelim} when
$v\in H^s(\R)$ and is finite together with $\|v\|_2$ exactly when
$v\in H^s(\R)$ \cite{DNPV}. The following lemma is the
main computation of the paper. Part (i) computes the interaction of
two intervals centered at the origin, part (ii) expresses the kinetic energy of an integral of indicator
functions of such intervals through these interactions, and
part (iii) combines them for $u_H=\rho_H^{1/r}$. The convexity of $T_s(\rho_{H_\theta})$ along the height interpolation
is proved in \cite{DV}. Here the kinetic energy is evaluated at
$\rho_{H_\theta}^{1/r}$ instead; the series in (iii) has convex terms,
which gives its convexity, and
single terms of it give the lower bounds for the second variation in
Section~\ref{sec:coercivity}.
\Needspace{13\baselineskip}
\begin{lemma}\label{lem:kinetic-series}
Let $0<s<1$ and define
\begin{equation*}
 d_{0,s}=1,\qquad d_{k,s}=\frac{\prod_{j=0}^{2k-1}(2s+j)}{(2k+1)!}
 \quad(k\ge1),\qquad C_{s,r}=\frac{2^{2s+1}c_{1,s}}{sr^2}.
\end{equation*}
\begin{enumerate}
\item[(i)] For $A\ge B\ge0$ with $A>0$ let
$\mathcal B_s(A,B)=c_{1,s}\int_{-B}^B\int_{\R\setminus(-A,A)}|x-y|^{-1-2s}\dd y\dd x$
and $\mathcal B_s(B,A)=\mathcal B_s(A,B)$. Then
\begin{equation}\label{eq:small-interaction}
 \mathcal B_s(A,B)=\frac{c_{1,s}}s\int_{A-B}^{A+B}t^{-2s}\dd t
 =\frac{2c_{1,s}}s\sum_{k=0}^\infty d_{k,s}A^{-2s-2k}B^{2k+1}
 \quad\text{in }[0,\infty],
\end{equation}
which is finite if $A>B$, or if $s<1/2$. For $s<1/2$ one has
$\sum_{k\ge0}d_{k,s}=2^{-2s}/(1-2s)$ and
$\mathcal B_s(A,B)=\frac{c_{1,s}}{s(1-2s)}\bigl((A+B)^{1-2s}-(A-B)^{1-2s}\bigr)$;
for $s=1/2$ and $A>B$, $\mathcal B_s(A,B)=2c_{1,1/2}\log\frac{A+B}{A-B}$.
\item[(ii)] Let $(\Theta,\nu)$ be a $\sigma$-finite measure space,
$R:\Theta\to(0,\infty)$ measurable, and suppose that
$v(x)=\int_\Theta\1_{\{|x|<R(\vartheta)\}}\dd\nu(\vartheta)$ is finite for
almost every $x$. Then
\begin{equation}\label{eq:small-layer}
 T_s(v)=\int_\Theta\int_\Theta \mathcal B_s\bigl(R(\vartheta),R(\vartheta')\bigr)
 \dd\nu(\vartheta)\dd\nu(\vartheta')\in[0,\infty].
\end{equation}
\item[(iii)] Let $H$ satisfy the hypotheses of
Lemma~\ref{lem:height-coordinates}(ii) and $u_H=\rho_H^{1/r}$. Then
\begin{equation}\label{eq:kinetic-series}
 T_s(u_H)=C_{s,r}\sum_{k=0}^\infty d_{k,s}
 \int_{0<\xi<\zeta<M}
 \frac{H'(\xi)^{1+2s+2k}}
      {H'(\zeta)^{2k}H(\xi)^\beta H(\zeta)^\beta}\dd\xi\dd\zeta
 \quad\text{in }[0,\infty].
\end{equation}
\end{enumerate}
\end{lemma}
\begin{proof}
(i) For $-B<x<B$, 
\begin{equation*}
 \int_A^\infty\frac{\dd y}{|x-y|^{1+2s}}=\int_{A-x}^\infty t^{-1-2s}\dd t
 =\frac{(A-x)^{-2s}}{2s},\qquad
 \int_{-\infty}^{-A}\frac{\dd y}{|x-y|^{1+2s}}=\frac{(A+x)^{-2s}}{2s},
\end{equation*}
and the substitutions $t=A-x$ and $t=A+x$ give
\begin{equation*}
 \mathcal B_s(A,B)=\frac{c_{1,s}}{2s}\int_{-B}^B\bigl((A-x)^{-2s}+(A+x)^{-2s}\bigr)\dd x
 =\frac{c_{1,s}}s\int_{A-B}^{A+B}t^{-2s}\dd t,
\end{equation*}
which is finite unless $A=B$ and $s\ge1/2$; evaluating the integral
gives the closed forms for $s<1/2$ and for $s=1/2$. For $0\le B<A$, the
substitution $t=A(1+z)$ gives
\begin{equation*}
 \mathcal B_s(A,B)=\frac{c_{1,s}}sA^{1-2s}\int_{-B/A}^{B/A}(1+z)^{-2s}\dd z.
\end{equation*}
The binomial series
$(1+z)^{-2s}=1+\sum_{n\ge1}\frac{(-1)^n\prod_{j=0}^{n-1}(2s+j)}{n!}z^n$
converges uniformly on $[-q,q]$ for $0\le q<1$; integrating term by
term, the odd powers integrate to zero and
$\int_{-q}^qz^{2k}\dd z=2q^{2k+1}/(2k+1)$, so that
\begin{equation}\label{eq:integrated-binomial-series}
 \frac12\int_{-q}^q(1+z)^{-2s}\dd z
 =\sum_{k=0}^\infty\frac{\prod_{j=0}^{2k-1}(2s+j)}{(2k)!\,(2k+1)}\,q^{2k+1}
 =\sum_{k=0}^\infty d_{k,s}q^{2k+1}
 \qquad(0\le q<1),
\end{equation}
since $(2k)!(2k+1)=(2k+1)!$. With $q=B/A$ this is the series in
\eqref{eq:small-interaction}. As $B\to A^-$, the left side of
\eqref{eq:small-interaction} increases to $\frac{c_{1,s}}s\int_0^{2A}t^{-2s}\dd t$,
and the right side increases to the series at $B=A$ by monotone
convergence (each term $d_{k,s}A^{-2s-2k}B^{2k+1}$ is nonnegative and
nondecreasing in $B$), so the identity extends to $A=B$ in $[0,\infty]$; for
$s<1/2$ and $A=B=1$ it gives
$\sum_kd_{k,s}=\frac12\int_{-1}^1(1+z)^{-2s}\dd z=2^{-2s}/(1-2s)$.

(ii) Write $\chi_\vartheta=\1_{\{|\cdot|<R(\vartheta)\}}$. If $|x|\le|y|$,
then $\chi_\vartheta(x)-\chi_\vartheta(y)\ge0$ for every $\vartheta$, and
if $|x|\ge|y|$ these differences are all nonpositive. Hence, whenever
$v(x),v(y)<\infty$,
\begin{equation*}
 (v(x)-v(y))^2=\int_\Theta\int_\Theta
 \bigl(\chi_\vartheta(x)-\chi_\vartheta(y)\bigr)
 \bigl(\chi_{\vartheta'}(x)-\chi_{\vartheta'}(y)\bigr)
 \dd\nu(\vartheta)\dd\nu(\vartheta')
\end{equation*}
with a nonnegative integrand, and Tonelli's theorem allows us to
insert this into \eqref{eq:kinetic-double} and integrate first in
$(x,y)$. For radii $A=\max\{R(\vartheta),R(\vartheta')\}$ and
$B=\min\{R(\vartheta),R(\vartheta')\}$, the two differences are
simultaneously nonzero exactly when one of the points lies in $(-B,B)$
and the other outside $(-A,A)$, and their product is then one:
\begin{equation*}
 \bigl(\chi_\vartheta(x)-\chi_\vartheta(y)\bigr)
 \bigl(\chi_{\vartheta'}(x)-\chi_{\vartheta'}(y)\bigr)
 =\1_{\{|x|<B,\,|y|\ge A\}}+\1_{\{|y|<B,\,|x|\ge A\}}.
\end{equation*}
Interchanging $x$ and $y$ in the second term,
\begin{equation*}
 \frac{c_{1,s}}2\iint_{\R^2}
 \frac{(\chi_\vartheta(x)-\chi_\vartheta(y))(\chi_{\vartheta'}(x)-\chi_{\vartheta'}(y))}
 {|x-y|^{1+2s}}\dd x\dd y
 =c_{1,s}\int_{-B}^B\int_{\R\setminus(-A,A)}\frac{\dd y\dd x}{|x-y|^{1+2s}}
 =\mathcal B_s(A,B),
\end{equation*}
also when $A=B$. This is \eqref{eq:small-layer}.

(iii) By \eqref{eq:monotone-root}, $u_H$ is of the form in (ii) with
$\Theta=(0,M)$, $R=R_H$ and $\dd\nu=w\dd\xi$, where
$w=r^{-1}H^{-\beta}H'$, and $u_H\le H(M)^{1/r}$. Since $\mathcal B_s$ is symmetric and 
$R_H(\xi)\ge R_H(\zeta)$ for $\xi<\zeta$, \eqref{eq:small-layer} gives
\begin{equation}\label{eq:kinetic-triangular}
 T_s(u_H)=2\int_{0<\xi<\zeta<M}
 \mathcal B_s\bigl(R_H(\xi),R_H(\zeta)\bigr)w(\xi)w(\zeta)\dd\xi\dd\zeta.
\end{equation}
Inserting \eqref{eq:small-interaction} with $A=1/(2H'(\xi))$ and
$B=1/(2H'(\zeta))$, for which
\begin{equation*}
 A^{-2s-2k}B^{2k+1}=2^{2s-1}\,\frac{H'(\xi)^{2s+2k}}{H'(\zeta)^{2k+1}},\qquad
 w(\xi)w(\zeta)=\frac{H'(\xi)H'(\zeta)}{r^2H(\xi)^\beta H(\zeta)^\beta},
\end{equation*}
we obtain
\begin{equation*}
 \mathcal B_s\bigl(R_H(\xi),R_H(\zeta)\bigr)w(\xi)w(\zeta)
 =\frac{2^{2s}c_{1,s}}{sr^2}\sum_{k=0}^\infty d_{k,s}\,
 \frac{H'(\xi)^{1+2s+2k}}
 {H'(\zeta)^{2k}H(\xi)^\beta H(\zeta)^\beta}.
\end{equation*}
Substituting this into \eqref{eq:kinetic-triangular} and applying
monotone convergence to the nonnegative partial sums gives
\eqref{eq:kinetic-series} with $C_{s,r}=2\cdot2^{2s}c_{1,s}/(sr^2)$.
\end{proof}

\begin{lemma}\label{lem:kinetic-convexity}
Let $0<s<1$ and $1\le r\le1/(1-s)$, so that $\beta\le s$. Let $H_0,H_1$
satisfy the hypotheses of Lemma~\ref{lem:height-coordinates}(ii), and
let $u_\theta$ be given by \eqref{eq:interpolation}. Then
$\theta\mapsto T_s(u_\theta)\in[0,\infty]$ is convex on $[0,1]$; in
particular,
\begin{equation}\label{eq:kinetic-convexity}
 T_s(u_\theta)\le(1-\theta)T_s(u_0)+\theta T_s(u_1),
 \qquad 0\le\theta\le1.
\end{equation}
\end{lemma}
\begin{proof}
Fix $0<\xi<\zeta<M$ and $k\ge0$. The arguments $H_\theta'(\xi)$,
$H_\theta'(\zeta)$, $H_\theta(\xi)$, $H_\theta(\zeta)$ are positive and
affine in $\theta$, and the integrand of the $k$th term in
\eqref{eq:kinetic-series} is the monomial $X^P/(Y_1^{q_1}Y_2^{q_2}Y_3^{q_3})$
of Lemma~\ref{lem:monomial} with $ P=1+2s+2k$, $q_1=2k$,  $ q_2=q_3=\beta$,
\begin{gather*}
 X=H_\theta'(\xi),\qquad Y_1=H_\theta'(\zeta),\qquad Y_2=H_\theta(\xi),\text{ and }
 Y_3=H_\theta(\zeta)
\end{gather*}
for which $P-1-\sum_iq_i=2(s-\beta)\ge0$. Hence each integrand is a
convex function of $\theta$, being the composition of a convex function
with an affine map, and so are its integral over $0<\xi<\zeta<M$ and the
sum over $k$, with values in $[0,\infty]$. In particular,
\eqref{eq:kinetic-convexity} holds.
\end{proof}

The next lemma collects what is needed in case \eqref{eq:caseA}: the
action $S$ is finite and convex along the height interpolation of two
positive solutions, and the interpolation is continuous in $H^s$ at the
endpoints.

\begin{lemma}\label{lem:path-finiteness}
Assume \eqref{eq:f-regularity}, \eqref{eq:f-assumptions} and
\eqref{eq:caseA}. Let $H_0,H_1$ be admissible on $[0,M]$ with
$u_0,u_1\in H^s(\R)$ and $F(u_0),F(u_1)\in L^1(\R)$, and let $u_\theta$ be
given by \eqref{eq:interpolation}. Then,
\begin{enumerate}
\item[(i)] For $0\le\theta\le1$, $u_\theta\in H^s(\R)$,
$F(u_\theta)\in L^1(\R)$, $b_r(H_\theta)\in L^1(0,M)$, and
\begin{equation}\label{eq:height-action}
 S(u_\theta)=\frac12T_s(u_\theta)
                    -\frac1r\int_0^M b_r(H_\theta(\xi))\dd\xi.
\end{equation}
\item[(ii)] $\theta\mapsto S(u_\theta)$ is convex on $[0,1]$: for
$\theta_0,\theta_1,t\in[0,1]$,
\begin{equation}\label{eq:action-convexity}
 S(u_{(1-t)\theta_0+t\theta_1})
       \le(1-t)S(u_{\theta_0})+tS(u_{\theta_1}).
\end{equation}
\item[(iii)] $\|u_\theta-u_0\|_{H^s}+\|u_\theta-u_0\|_\infty\to0$ as
$\theta\to0^+$, and $\|u_\theta-u_1\|_{H^s}+\|u_\theta-u_1\|_\infty\to0$ as
$\theta\to1^-$.
\end{enumerate}
\end{lemma}
\begin{proof}
\emph{Step 1. $u_\theta\in H^s(\R)$.}
Since $r\ge2$, the function $h\mapsto h^{2/r-1}$ is convex on $(0,\infty)$
(constant if $r=2$), and \eqref{eq:powers} with $q=2$ gives
\begin{equation*}
 \|u_\theta\|_2^2=\frac2r\int_0^MH_\theta^{2/r-1}\dd\xi
 \le\frac2r\int_0^M\bigl((1-\theta)H_0^{2/r-1}+\theta H_1^{2/r-1}\bigr)\dd\xi
 =(1-\theta)\|u_0\|_2^2+\theta\|u_1\|_2^2<\infty.
\end{equation*}
Together with \eqref{eq:kinetic-convexity}, this proves $u_\theta\in H^s(\R)$.

\emph{Step 2. $F(u_\theta)\in L^1(\R)$.}
Let $m_*=\max\{\|u_0\|_\infty,\|u_1\|_\infty\}$, so that $0<u_\theta\le m_*$ by
Lemma~\ref{lem:height-interpolation}. Choose $K,f_K,F_K$ as
in \eqref{eq:aux-definition} on $[0,m_*]$, and define $\Phi_K(h)=F_K(h^{1/r})$ for
$0\le h\le m_*^r$. Since $F_K(t)=\int_0^tf_K$, the substitution $t=\tau^{1/r}$ gives
\begin{equation*}
 \Phi_K(h)=\int_0^{h^{1/r}}f_K(t)\dd t=\int_0^h\frac{f_K(\tau^{1/r})}{r\tau^\beta}\dd\tau,
\end{equation*}
so $\Phi_K$ is absolutely continuous with $\Phi_K(0)=0$, and with $t=h^{1/r}$,
\begin{equation*}
 \Phi_K'(h)=\frac{f_K(t)}{rt^{r-1}}>0,\qquad
 \Phi_K''(h)=\frac{tf_K'(t)-(r-1)f_K(t)}{r^2t^{2r-1}}<0,
\end{equation*}
where the last sign follows from $tf_K'(t)<f_K(t)$ in \eqref{eq:scalar-bounds}
and $r\ge2$. For $0<c\le1$, the monotonicity of $f_K$ gives
\begin{equation*}
 \Phi_K'(ch)=\frac{f_K(c^{1/r}h^{1/r})}{rc^\beta h^\beta}\le c^{-\beta}\Phi_K'(h),
\end{equation*}
and since $H_\theta\ge(1-\theta)H_0$ and $\Phi_K'$ is decreasing,
\begin{equation}\label{eq:Phi-height-bound}
 \Phi_K'(H_\theta(\xi))\le\Phi_K'((1-\theta)H_0(\xi))\le(1-\theta)^{-\beta}\Phi_K'(H_0(\xi)),
 \qquad 0\le\theta<1.
\end{equation}
Lemma~\ref{lem:layer} with $J=\Phi_K$ (applied first to $H_0$
and then to $H_\theta$) gives
\begin{equation}\label{eq:A-interpolation-bound}
 \begin{split}
 \int_\R F_K(u_\theta)\dd x=\int_0^M\Phi_K'(H_\theta(\xi))\dd\xi
 &\le(1-\theta)^{-\beta}\int_0^M\Phi_K'(H_0(\xi))\dd\xi
 =(1-\theta)^{-\beta}\int_\R F_K(u_0)\dd x\\
 &\le(1-\theta)^{-\beta}\Bigl(\frac K2\|u_0\|_2^2+\|F(u_0)\|_1\Bigr),
 \end{split}
\end{equation}
which is finite; at $\theta=1$ the integral is finite by the
assumption on $u_1$. Since, by \eqref{eq:scalar-bounds},
\begin{equation*}
 |F(t)|\le\frac K2t^2+F_K(t),\qquad t|f(t)|\le Kt^2+tf_K(t)\le Kt^2+2F_K(t)
 \qquad(0<t\le m_*),
\end{equation*}
both $F(u_\theta)$ and $u_\theta f(u_\theta)$ are integrable.

\emph{Step 3. A representation of $S$ in terms of height functions.}
From $f_K(t)=Kt-f(t)$ and the formula for $\Phi_K'$,
\begin{equation*}
 b_r(h)=\frac{f(h^{1/r})-\lambda h^{1/r}}{h^\beta}
 =\frac{(K-\lambda)h^{1/r}-f_K(h^{1/r})}{h^\beta}
 =(K-\lambda)h^{2/r-1}-r\Phi_K'(h),
\end{equation*}
so that, by \eqref{eq:powers} and \eqref{eq:A-interpolation-bound},
\begin{equation}\label{eq:b-integrability}
 \begin{split}
 \int_0^M|b_r(H_\theta(\xi))|\dd\xi
 &\le(K-\lambda)\int_0^MH_\theta^{2/r-1}\dd\xi+r\int_0^M\Phi_K'(H_\theta)\dd\xi\\
 &=\frac{r(K-\lambda)}2\|u_\theta\|_2^2+r\int_\R F_K(u_\theta)\dd x<\infty.
 \end{split}
\end{equation}
The function $J(h)=\frac\lambda2h^{2/r}-F(h^{1/r})=\frac{\lambda-K}2h^{2/r}+\Phi_K(h)$
is absolutely continuous on $[0,m_*^r]$ with $J(0)=0$ and
\begin{equation*}
 J'(h)=\frac{\lambda-K}rh^{2/r-1}+\Phi_K'(h)=-\frac1rb_r(h),
\end{equation*}
so \eqref{eq:b-integrability} is the hypothesis of the second case of
Lemma~\ref{lem:layer}, and \eqref{eq:layer} gives
\begin{equation*}
 \frac\lambda2\|u_\theta\|_2^2-\int_\R F(u_\theta)\dd x
 =\int_\R J(\rho_{H_\theta})\dd x=\int_0^MJ'(H_\theta(\xi))\dd\xi
 =-\frac1r\int_0^Mb_r(H_\theta(\xi))\dd\xi.
\end{equation*}
Adding $\frac12T_s(u_\theta)$ gives \eqref{eq:height-action}.

\emph{Step 4. Convexity.}
For $\theta_t=(1-t)\theta_0+t\theta_1$ the affine definition gives
$H_{\theta_t}=(1-t)H_{\theta_0}+tH_{\theta_1}$, with admissible height functions by
Lemma~\ref{lem:height-interpolation}. Lemma~\ref{lem:kinetic-convexity}
applied to $H_{\theta_0},H_{\theta_1}$ gives
$T_s(u_{\theta_t})\le(1-t)T_s(u_{\theta_0})+tT_s(u_{\theta_1})$, and the
concavity of $b_r$ gives, for $0<\xi<M$,
\begin{equation*}
 -b_r(H_{\theta_t}(\xi))\le-(1-t)b_r(H_{\theta_0}(\xi))-tb_r(H_{\theta_1}(\xi)).
\end{equation*}
All three local terms are integrable by Step~3; integrating, multiplying
by $1/r$ and adding half the kinetic inequality gives
\eqref{eq:action-convexity} by \eqref{eq:height-action}.

\emph{Step 5. Pointwise and uniform convergence at the endpoint.}
By symmetry it suffices to consider $\theta\to0^+$. We have
\begin{gather*}
 \|H_\theta-H_0\|_{C([0,M])}=\theta\|H_1-H_0\|_{C([0,M])}\to0,\\
 R_{H_\theta}(\xi)=\frac1{2((1-\theta)H_0'(\xi)+\theta H_1'(\xi))}\to R_{H_0}(\xi),
 \qquad 0<\xi<M.
\end{gather*}
Fix $x>0$ and write $\xi_\theta(x)=R_{H_\theta}^{-1}(x)$. For
$0<\varepsilon<\min\{\xi_0(x),M-\xi_0(x)\}$, strict decrease of $R_{H_0}$
gives $R_{H_0}(\xi_0(x)-\varepsilon)>x>R_{H_0}(\xi_0(x)+\varepsilon)$, the
same inequalities hold for $R_{H_\theta}$ when $\theta$ is small, and
inverting them gives $\xi_0(x)-\varepsilon<\xi_\theta(x)<\xi_0(x)+\varepsilon$.
Thus $\xi_\theta(x)\to\xi_0(x)$, and
\begin{equation*}
 u_\theta(x)^r=H_\theta(\xi_\theta(x))\longrightarrow H_0(\xi_0(x))=u_0(x)^r
 \quad\text{as }\theta\to0^+\text{, for every }x>0,
\end{equation*}
while $u_\theta(0)^r=H_\theta(M)\to H_0(M)=u_0(0)^r$; by evenness,
$u_\theta\to u_0$ pointwise on $\R$. Since all these functions are
nonincreasing in $|x|$ and $u_0$ is continuous, the convergence is
uniform on compact sets: given $L,\varepsilon>0$, choose a partition
$0=x_0<\dots<x_n=L$ with $u_0(x_j)-u_0(x_{j+1})\le\varepsilon$; then
monotonicity gives
$|u_\theta(x)-u_0(x)|\le\max_j|u_\theta(x_j)-u_0(x_j)|+\varepsilon$ on
$[0,L]$, hence $\limsup_{\theta\to0^+}\|u_\theta-u_0\|_{L^\infty([-L,L])}\le\varepsilon$.
On $\{|x|\ge L\}$, both $u_\theta$ and $u_0$ are small uniformly in
$\theta$: since $u_\theta$ is symmetric decreasing with
$\int_\R u_\theta^r\dd y=M$, we have
$2|x|u_\theta(x)^r\le\int_{-|x|}^{|x|}u_\theta^r\dd y\le M$, and
likewise for $u_0$, so that
$\sup_{|x|\ge L}|u_\theta(x)-u_0(x)|\le2(M/(2L))^{1/r}$ for every
$\theta\in[0,1]$. Letting
first $\theta\to0^+$ and then $L\to\infty$ proves $\|u_\theta-u_0\|_\infty\to0$.

\emph{Step 6. Convergence in $L^2$ and in $H^s$.}
For $0\le\theta\le1/2$ we have $H_\theta\ge H_0/2$, hence
$0<H_\theta^{2/r-1}\le2^{1-2/r}H_0^{2/r-1}$ with
$\int_0^MH_0^{2/r-1}\dd\xi=\frac r2\|u_0\|_2^2<\infty$, and dominated
convergence in \eqref{eq:powers} gives
\begin{equation*}
 \|u_\theta\|_2^2=\frac2r\int_0^MH_\theta(\xi)^{2/r-1}\dd\xi
 \longrightarrow\frac2r\int_0^MH_0(\xi)^{2/r-1}\dd\xi=\|u_0\|_2^2\quad\text{as }\theta\to0^+.
\end{equation*}
Since $u_\theta,u_0\ge0$, we have $(u_\theta-u_0)^2\le|u_\theta^2-u_0^2|$, and
pointwise convergence with $\min\{u_\theta^2,u_0^2\}\le u_0^2\in L^1$ gives
\begin{equation*}
 \|u_\theta-u_0\|_2^2\le\int_\R|u_\theta^2-u_0^2|\dd x
 =\|u_\theta\|_2^2+\|u_0\|_2^2-2\int_\R\min\{u_\theta^2,u_0^2\}\dd x
 \longrightarrow0\quad\text{as }\theta\to0^+.
\end{equation*}
Finally, \eqref{eq:kinetic-convexity} and the $L^2$ bound make
$\{u_\theta\}_{0\le\theta\le1}$ bounded in $H^s$; every weak limit along a subsequence
equals $u_0$ by the strong $L^2$ convergence, so $u_\theta\rightharpoonup u_0$
in $H^s$, and weak lower semicontinuity together with
\eqref{eq:kinetic-convexity} gives
\begin{equation*}
 T_s(u_0)\le\liminf_{\theta\to0^+}T_s(u_\theta)
 \le\limsup_{\theta\to0^+}T_s(u_\theta)
 \le\lim_{\theta\to0^+}\bigl((1-\theta)T_s(u_0)+\theta T_s(u_1)\bigr)=T_s(u_0).
\end{equation*}
Convergence of the norms and weak convergence give
$\|u_\theta-u_0\|_{H^s}^2=\|u_\theta\|_{H^s}^2+\|u_0\|_{H^s}^2-2\langle u_\theta,u_0\rangle_{H^s}\to0$.
\end{proof}

\section{The two-level interpolation}\label{sec:two-level}

Throughout this section, $0<s<1/2$, $\lambda>0$, $\mu>0$ and
$r\in[1,1/(1-s)]$, so that $0\le\beta\le s$. For such $s$ the exponent
$r$ must be smaller than $2$, and then $b_r$ cannot be concave near
$h=0$ when $f(t)/t$ has a finite limit at zero. To overcome this difficulty, we introduce a new interpolation. Above
the level $\mu$ it interpolates the height functions of
Section~\ref{sec:height}, while below $\mu$ it interpolates linearly
the powers $R(t)^{1-2s}$ of the half-lengths
$R(t)=\frac12|\{u>t\}|$ of the superlevel sets $\{u>t\}$, $0<t<\mu$.

For $M>0$, let
\begin{equation}\label{eq:class-A}
 \mathcal A_M=\Bigl\{v\in H^s(\R)\cap L^1(\R)\cap L^\infty(\R):\
 v\ge0\text{ symmetric decreasing},\ \mathcal N(v)=M\Bigr\},
\end{equation}
where $\mathcal N(v)=\int_\R(v-\mu)_+^r\dd x$ as in
Section~\ref{sec:prelim}. 
Let $v\in\mathcal A_M$. Then $|\{v>\mu\}|>0$ and $|\{v>t\}|\le\|v\|_1/t<\infty$
for $t>0$. We define
\begin{equation}\label{eq:two-level-data}
 R_v(t)=\frac12|\{v>t\}|\quad(0<t<\mu),\qquad \rho_v=(v-\mu)_+^r,
\end{equation}
and we let $H_v$ be the height function of $\rho_v$ (Lemma~\ref{lem:height-coordinates}(i)).
We set $R_{H_v}=1/(2H_v')$ on $[0,M)$ and $w_v=r^{-1}H_v^{-\beta}H_v'$ on
$(0,M)$. The function $R_v$ is nonincreasing with values in $(0,\infty)$,
$H_v(M)=(\|v\|_\infty-\mu)^r$, $H_v'(0)=1/|\{v>\mu\}|>0$ by
\eqref{eq:height-derivative-general}, and since
$\{v>t\}\supset\{v>\mu\}=\{\rho_v>0\}\supset\{\rho_v>H_v(\xi)\}$,
\begin{equation}\label{eq:two-level-nest}
 R_v(t)\ge\frac12|\{v>\mu\}|=R_{H_v}(0)\ge R_{H_v}(\xi)
 \qquad(0<t<\mu,\ 0\le\xi<M).
\end{equation}
By the layer cake formula applied to $\min\{v,\mu\}$,
Lemma~\ref{lem:height-coordinates}(iii) applied to $\rho_v$, and
\eqref{eq:monotone-root} applied to $\rho_v^{1/r}=(v-\mu)_+$, we have, for
almost every $x\in\R$,
\begin{equation}\label{eq:two-level-layers}
 v(x)=\min\{v(x),\mu\}+(v(x)-\mu)_+
 =\int_0^\mu\1_{\{|x|<R_v(t)\}}\dd t
 +\int_0^M\1_{\{|x|<R_{H_v}(\xi)\}}w_v(\xi)\dd\xi.
\end{equation}

Let $u_0,u_1\in\mathcal A_M$, and write $R_i=R_{u_i}$ and $H_i=H_{u_i}$
for $i=0,1$. For $0\le\theta\le1$ we define
\begin{equation}\label{eq:two-level-def}
 \begin{split}
 &U_\theta(t)=(1-\theta)R_0(t)^{1-2s}+\theta R_1(t)^{1-2s},\qquad
 R_\theta(t)=U_\theta(t)^{1/(1-2s)}\qquad(0<t<\mu),\\
 &H_\theta=(1-\theta)H_0+\theta H_1,\qquad
 R_{H_\theta}=\frac1{2H_\theta'}\ \text{ on }[0,M),\qquad
 w_\theta=\frac1r\,H_\theta^{-\beta}H_\theta'\ \text{ on }(0,M),
 \end{split}
\end{equation}
and
\begin{equation}\label{eq:two-level}
 u_\theta(x)=\int_0^\mu\1_{\{|x|<R_\theta(t)\}}\dd t
 +\int_0^M\1_{\{|x|<R_{H_\theta}(\xi)\}}w_\theta(\xi)\dd\xi,
 \qquad x\in\R.
\end{equation}
We call $\{u_\theta\}_{0\le\theta\le1}$ the two-level interpolation between
$u_0$ and $u_1$. The height function $H_\theta$ satisfies the hypotheses of
Lemma~\ref{lem:height-coordinates}(ii), with
$H_\theta'=(1-\theta)H_0'+\theta H_1'>0$ on $[0,M)$, and
\eqref{eq:monotone-root} identifies the second integral in
\eqref{eq:two-level} with $\rho_{H_\theta}^{1/r}$.

\begin{lemma}\label{lem:two-level}
Let $u_0,u_1\in\mathcal A_M$ and $0\le\theta\le1$.
\begin{enumerate}
\item[(i)] $R_\theta(t)\ge R_{H_\theta}(0)\ge R_{H_\theta}(\xi)$ for
$0<t<\mu$ and $0\le\xi<M$.
\item[(ii)] $u_\theta$ is symmetric decreasing, and
\begin{equation*}
 \min\{u_\theta,\mu\}=\int_0^\mu\1_{\{|\cdot|<R_\theta(t)\}}\dd t,\qquad
 (u_\theta-\mu)_+=\rho_{H_\theta}^{1/r}.
\end{equation*}
\item[(iii)] $|\{u_\theta>t\}|=2R_\theta(t)$ for all but countably many
$t\in(0,\mu)$, $\{u_\theta>t\}=\{\rho_{H_\theta}>(t-\mu)^r\}$ for
$t\ge\mu$, $\mathcal N(u_\theta)=M$, $u_\theta\in L^1(\R)\cap L^\infty(\R)$,
and $\|u_\theta\|_\infty\le\max\{\|u_0\|_\infty,\|u_1\|_\infty\}$.
\item[(iv)] $u_\theta=u_i$ almost everywhere for $\theta=i\in\{0,1\}$.
\end{enumerate}
\end{lemma}
\begin{proof}
(i) For $q\in\R\setminus\{0\}$ and $x_0,x_1>0$, let
$\mathcal P_q(x_0,x_1)=((1-\theta)x_0^q+\theta x_1^q)^{1/q}$. The identity $H_\theta'(0)=(1-\theta)H_0'(0)+\theta H_1'(0)$
gives
$R_{H_\theta}(0)^{-1}=(1-\theta)R_{H_0}(0)^{-1}+\theta R_{H_1}(0)^{-1}$,
that is, $R_{H_\theta}(0)=\mathcal P_{-1}(R_{H_0}(0),R_{H_1}(0))$. Power means are nondecreasing in each argument and in $q$. It follows from \eqref{eq:two-level-nest}  and   $0<1-2s<1$ that 
\begin{equation*}
 R_\theta(t)=\mathcal P_{1-2s}(R_0(t),R_1(t))
 \ge \mathcal P_{1-2s}(R_{H_0}(0),R_{H_1}(0))
 \ge \mathcal P_{-1}(R_{H_0}(0),R_{H_1}(0))=R_{H_\theta}(0),
\end{equation*}
and $R_{H_\theta}$ is nonincreasing because $H_\theta'$ is nondecreasing.

(ii) By (i), if $|x|<R_{H_\theta}(\xi)$ for some $\xi$, then
$|x|<R_\theta(t)$ for every $t\in(0,\mu)$, so the first integral in
\eqref{eq:two-level} equals $\mu$ on the support of the second one, and
the second one equals $\rho_{H_\theta}(x)^{1/r}$ by
\eqref{eq:monotone-root}. Hence $u_\theta=\mu+\rho_{H_\theta}^{1/r}$ where
$\rho_{H_\theta}>0$, and $u_\theta=\int_0^\mu\1_{\{|x|<R_\theta(t)\}}\dd t\le\mu$
elsewhere, which gives both identities. Both integrands in
\eqref{eq:two-level} are nonincreasing in $|x|$ and even, so $u_\theta$
is symmetric decreasing.

(iii) By (ii), for $0<t<\mu$ the set $\{u_\theta>t\}$ consists of the
points $x$ with $|\{\tau\in(0,\mu):R_\theta(\tau)>|x|\}|>t$. Since
$R_\theta$ is nonincreasing, $\{\tau:R_\theta(\tau)>|x|\}$ is an initial
segment of $(0,\mu)$ of length
$\tau_*(x)=\sup\{\tau\in(0,\mu):R_\theta(\tau)>|x|\}$ (with $\sup\emptyset=0$),
so that $\min\{u_\theta(x),\mu\}=\tau_*(x)$. If $\tau_*(x)>t$, then
$R_\theta(\tau)>|x|$ for some $\tau>t$, and $R_\theta(t)\ge R_\theta(\tau)>|x|$.
Conversely, if $|x|<R_\theta(t)$ and $R_\theta$ is continuous at $t$,
then $R_\theta(\tau)>|x|$ for all $\tau$ in a neighborhood of $t$, so
$\tau_*(x)>t$. Hence $\{u_\theta>t\}=\{|x|<R_\theta(t)\}$ is an interval of
length $2R_\theta(t)$ at every continuity point $t$ of $R_\theta$, and the
monotone function $R_\theta$ has at most countably many
discontinuities. For $t\ge\mu$, (ii) gives
$\{u_\theta>t\}=\{\rho_{H_\theta}^{1/r}>t-\mu\}=\{\rho_{H_\theta}>(t-\mu)^r\}$,
and $\mathcal N(u_\theta)=\int_\R\rho_{H_\theta}\dd x=M$ by
Lemma~\ref{lem:height-coordinates}(ii). Next,
\begin{equation*}
 (u_\theta-\mu)_+\le H_\theta(M)^{1/r}
 =\bigl((1-\theta)H_0(M)+\theta H_1(M)\bigr)^{1/r}
 \le\max_iH_i(M)^{1/r}=\max_i\bigl(\|u_i\|_\infty-\mu\bigr),
\end{equation*}
which gives the bound for $\|u_\theta\|_\infty$; and since
$R_\theta\le\max\{R_0,R_1\}\le R_0+R_1$ and the second integral in
\eqref{eq:two-level} is supported in $[-R_{H_\theta}(0),R_{H_\theta}(0)]$
and bounded by $H_\theta(M)^{1/r}$,
\begin{equation*}
 \|u_\theta\|_1\le\int_0^\mu2R_\theta(t)\dd t+2R_{H_\theta}(0)H_\theta(M)^{1/r}
 \le\|u_0\|_1+\|u_1\|_1+2R_{H_\theta}(0)H_\theta(M)^{1/r}<\infty,
\end{equation*}
where $\int_0^\mu2R_i(t)\dd t=\|\min\{u_i,\mu\}\|_1\le\|u_i\|_1$.

(iv) For $\theta=i$, \eqref{eq:two-level-def} gives $R_\theta=R_i$,
$H_\theta=H_i$, $R_{H_\theta}=R_{H_i}$ and $w_\theta=w_{u_i}$, so the right
side of \eqref{eq:two-level} is the right side of
\eqref{eq:two-level-layers} for $v=u_i$, which equals $u_i$ almost
everywhere.
\end{proof}

\begin{lemma}\label{lem:two-level-convexity}
Let $u_0,u_1\in\mathcal A_M$, and let $u_\theta$ be given by
\eqref{eq:two-level}.
\begin{enumerate}
\item[(i)] The function $\theta\mapsto T_s(u_\theta)$ is convex on $[0,1]$,
and $T_s(u_\theta)\le(1-\theta)T_s(u_0)+\theta T_s(u_1)<\infty$. In
particular, $u_\theta\in\mathcal A_M$ for $0\le\theta\le1$.
\item[(ii)] If $u_0\ne u_1$, then $\theta\mapsto T_s(u_\theta)$ is not
affine on $[0,1]$.
\item[(iii)] Assume in addition \eqref{eq:f-regularity},
\eqref{eq:f-assumptions} and \eqref{eq:caseB}. Then, for
$0\le\theta\le1$,
\begin{equation}\label{eq:two-level-action}
 S(u_\theta)=\frac12T_s(u_\theta)
 +2\int_0^\mu(\lambda t-f(t))R_\theta(t)\dd t
 -\frac1r\int_0^Mb_{r,\mu}(H_\theta(\xi))\dd\xi,
\end{equation}
where all terms are finite, and $\theta\mapsto S(u_\theta)$ is convex on
$[0,1]$. It is not affine if $u_0\ne u_1$.
\end{enumerate}
\end{lemma}
\begin{proof}
(i) By \eqref{eq:two-level}, $u_\theta$ is of the form of
Lemma~\ref{lem:kinetic-series}(ii). Indeed, we take as $\Theta$ the disjoint
union of the two intervals $(0,\mu)$ and $(0,M)$, formally
$\Theta=\bigl((0,\mu)\times\{1\}\bigr)\cup\bigl((0,M)\times\{2\}\bigr)$, with
the measure $\nu$ given by $\dd t$ on the first interval and by
$w_\theta(\xi)\dd\xi$ on the second, and with $R(t,1)=R_\theta(t)$ and
$R(\xi,2)=R_{H_\theta}(\xi)$. The function $u_\theta$ is finite
everywhere by Lemma~\ref{lem:two-level}(iii). The double integral in
\eqref{eq:small-layer} splits into four parts according to the
intervals in which $\vartheta$ and $\vartheta'$ lie. Since $\mathcal B_s$ is
symmetric, the two mixed parts are equal, and since the diagonals $\{t=t'\}$, $\{\xi=\zeta\}$ have measure zero, \eqref{eq:small-layer}
gives $T_s(u_\theta)=2(I_1(\theta)+I_2(\theta)+I_3(\theta))$ with
\begin{align*}
 I_1(\theta)&=\iint_{0<t<t'<\mu}\mathcal B_s\bigl(R_\theta(t),R_\theta(t')\bigr)
 \dd t\dd t',\\
 I_2(\theta)&=\int_0^\mu\int_0^M\mathcal B_s\bigl(R_\theta(t),R_{H_\theta}(\xi)\bigr)
 w_\theta(\xi)\dd\xi\dd t,\\
 I_3(\theta)&=\iint_{0<\xi<\zeta<M}\mathcal B_s\bigl(R_{H_\theta}(\xi),R_{H_\theta}(\zeta)\bigr)
 w_\theta(\xi)w_\theta(\zeta)\dd\xi\dd\zeta.
\end{align*}
In $I_1$ we have $R_\theta(t)\ge R_\theta(t')$; in $I_2$,
$R_\theta(t)\ge R_{H_\theta}(\xi)$ by Lemma~\ref{lem:two-level}(i); in
$I_3$, $R_{H_\theta}(\xi)\ge R_{H_\theta}(\zeta)$. We expand $\mathcal B_s$ by
\eqref{eq:small-interaction} with the larger radius in the first
argument. In $I_1$, with $R_\theta=U_\theta^{1/(1-2s)}$,
\begin{equation}\label{eq:two-level-T1}
 \mathcal B_s\bigl(R_\theta(t),R_\theta(t')\bigr)
 =\frac{2c_{1,s}}s\sum_{k=0}^\infty d_{k,s}\,R_\theta(t)^{-2s-2k}R_\theta(t')^{2k+1}
 =\frac{2c_{1,s}}s\sum_{k=0}^\infty d_{k,s}\,
 U_\theta(t)^{-\frac{2s+2k}{1-2s}}\,U_\theta(t')^{\frac{2k+1}{1-2s}}.
\end{equation}
In $I_2$, with $R_{H_\theta}(\xi)^{2k+1}=(2H_\theta'(\xi))^{-2k-1}$ and
$w_\theta=r^{-1}H_\theta^{-\beta}H_\theta'$,
\begin{equation}\label{eq:two-level-T2}
 \begin{split}
 \mathcal B_s\bigl(R_\theta(t),R_{H_\theta}(\xi)\bigr)w_\theta(\xi)
 &=\frac{2c_{1,s}}s\sum_{k=0}^\infty d_{k,s}\,
 \frac{U_\theta(t)^{-\frac{2s+2k}{1-2s}}}{(2H_\theta'(\xi))^{2k+1}}\cdot
 \frac{H_\theta'(\xi)}{rH_\theta(\xi)^\beta}\\
 &=\frac{2c_{1,s}}{sr}\sum_{k=0}^\infty\frac{d_{k,s}}{2^{2k+1}}\,
 U_\theta(t)^{-\frac{2s+2k}{1-2s}}\,H_\theta'(\xi)^{-2k}\,
 H_\theta(\xi)^{-\beta},
 \end{split}
\end{equation}
while $2I_3(\theta)=T_s(\rho_{H_\theta}^{1/r})$ is the series
\eqref{eq:kinetic-series} for the height function $H_\theta$. For
$t,t'\in(0,\mu)$ and $\xi\in(0,M)$, the quantities $U_\theta(t)$,
$U_\theta(t')$, $H_\theta'(\xi)$ and $H_\theta(\xi)$ are positive and
affine in $\theta$. Each term of \eqref{eq:two-level-T1} is convex in
$\theta$ by Lemma~\ref{lem:monomial}, applied with $X=U_\theta(t')$,
$P=(2k+1)/(1-2s)$ and the single exponent $q_1=(2s+2k)/(1-2s)$, for which
$
 P-1-q_1=0$;
each term of \eqref{eq:two-level-T2} has only nonpositive exponents and
is convex by Lemma~\ref{lem:nonpositive-monomial}; and $I_3$ is convex by
Lemma~\ref{lem:kinetic-convexity}, since $r\le1/(1-s)$. Sums and
integrals of nonnegative convex functions are convex with values in
$[0,\infty]$, so $I_1,I_2,I_3$ and $T_s(u_\theta)$ are convex on $[0,1]$,
and convexity gives
$T_s(u_\theta)\le(1-\theta)T_s(u_0)+\theta T_s(u_1)<\infty$. Since $u_\theta\in L^2$ by Lemma~\ref{lem:two-level}(iii) and
$T_s(u_\theta)<\infty$, we have $u_\theta\in H^s(\R)$, see the remark
after \eqref{eq:kinetic-double}.
With Lemma~\ref{lem:two-level}(ii),(iii) we get $u_\theta\in\mathcal A_M$.

(ii) Suppose that $u_0\ne u_1$. Then $R_0\ne R_1$ on a subset of
$(0,\mu)$ of positive measure, or $H_0'\ne H_1'$ on a subset of $(0,M)$
of positive measure: otherwise $H_0=H_1$, since $H_i(0)=0$ and the $H_i$
are absolutely continuous, and \eqref{eq:two-level-layers} gives
$u_0=u_1$ almost everywhere. Hence the set
\begin{equation*}
 E=\bigl\{(t,\xi)\in(0,\mu)\times(0,M):
 (U_1(t)-U_0(t),\,H_1'(\xi)-H_0'(\xi))\ne(0,0)\bigr\}
\end{equation*}
has positive measure. Since $d_{1,s}>0$, the term with $k=1$ in
\eqref{eq:two-level-T2} is a positive multiple of
\begin{equation*}
 m_\theta(t,\xi)=U_\theta(t)^{-\frac{2s+2}{1-2s}}H_\theta'(\xi)^{-2}H_\theta(\xi)^{-\beta},
\end{equation*}
and for $(t,\xi)\in E$ the function $\theta\mapsto m_\theta(t,\xi)$ is
strictly convex on $[0,1]$ by Lemma~\ref{lem:nonpositive-monomial},
applied with $Y_0=U_\theta(t)$, $Y_1=H_\theta'(\xi)$, $Y_2=H_\theta(\xi)$,
the exponents $\frac{2s+2}{1-2s},2,\beta$,
$(b_0,b_1)=(U_1(t)-U_0(t),H_1'(\xi)-H_0'(\xi))\ne(0,0)$ and
$b_2=H_1(\xi)-H_0(\xi)$. Let $0\le\theta_0<\theta_1\le1$, $0<\tau<1$ and
$\theta_\tau=(1-\tau)\theta_0+\tau\theta_1$. Then
$m_{\theta_\tau}<(1-\tau)m_{\theta_0}+\tau m_{\theta_1}$ on $E$. The
integrals of $m_{\theta_0}$, $m_{\theta_1}$ and $m_{\theta_\tau}$ over $E$
are finite, being bounded by a constant multiple of
$I_2(\theta)\le\frac12T_s(u_\theta)<\infty$ by (i), so integrating over
$E$ gives a strict convexity inequality for the contribution of $E$ to
$I_2$. All other contributions to $T_s(u_\theta)$ are finite and convex
in $\theta$ by the proof of (i). Hence
\begin{equation*}
 T_s(u_{\theta_\tau})<(1-\tau)T_s(u_{\theta_0})+\tau T_s(u_{\theta_1}),
\end{equation*}
so $\theta\mapsto T_s(u_\theta)$ is strictly convex on $[0,1]$, and in
particular not affine.

(iii) Let $m_*=\max\{\|u_0\|_\infty,\|u_1\|_\infty\}$. By
Lemma~\ref{lem:two-level}(iii), $0\le u_\theta\le m_*$ and $u_\theta\in L^1$,
and $\lambda t-f(t)$ is bounded on $[0,m_*]$ by \eqref{eq:f-regularity}.
Since $\frac\lambda2t^2-F(t)=\int_0^t(\lambda\tau-f(\tau))\dd\tau$, the
layer cake formula and Lemma~\ref{lem:two-level}(iii) give, with
absolutely convergent integrals,
\begin{equation*}
 \begin{split}
 \int_\R\Bigl(\frac\lambda2u_\theta^2-F(u_\theta)\Bigr)\dd x
 &=\int_0^{m_*}(\lambda t-f(t))\,|\{u_\theta>t\}|\dd t\\
 &=2\int_0^\mu(\lambda t-f(t))R_\theta(t)\dd t
 +\int_\mu^{m_*}(\lambda t-f(t))\,|\{\rho_{H_\theta}>(t-\mu)^r\}|\dd t.
 \end{split}
\end{equation*}
In the second integral we substitute $t=\mu+h^{1/r}$, for which
\begin{equation*}
 (\lambda t-f(t))\dd t
 =\bigl(\lambda(\mu+h^{1/r})-f(\mu+h^{1/r})\bigr)\frac{h^{1/r-1}}r\dd h
 =-\frac1rb_{r,\mu}(h)\dd h.
\end{equation*}
The second integral converges absolutely, since $\lambda t-f(t)$ is
bounded on $[\mu,m_*]$ and
$|\{\rho_{H_\theta}>(t-\mu)^r\}|\le|\{\rho_{H_\theta}>0\}|=1/H_\theta'(0)<\infty$,
and this is preserved by the substitution. Hence
\eqref{eq:monotone-cake} in its signed version, applied to
$\psi(h)=-r^{-1}b_{r,\mu}(h)\1_{\{h<(m_*-\mu)^r\}}$, gives
\begin{equation*}
 \int_\mu^{m_*}(\lambda t-f(t))\,|\{\rho_{H_\theta}>(t-\mu)^r\}|\dd t
 =-\frac1r\int_0^\infty b_{r,\mu}(h)\,|\{\rho_{H_\theta}>h\}|\dd h
 =-\frac1r\int_0^Mb_{r,\mu}(H_\theta(\xi))\dd\xi.
\end{equation*}
Adding $T_s(u_\theta)/2$, which is finite by (i), proves
\eqref{eq:two-level-action}. In the middle term of \eqref{eq:two-level-action},
$\lambda t-f(t)\ge0$ for $0<t<\mu$, because $f(\mu)\le\lambda\mu$ and
$f(t)/t$ is increasing, and $R_\theta(t)=U_\theta(t)^{1/(1-2s)}$ is
convex in $\theta$, as the composition of the convex increasing function
$U\mapsto U^{1/(1-2s)}$, $1/(1-2s)>1$, with the affine function
$\theta\mapsto U_\theta(t)$; hence the middle term is convex in $\theta$.
By \eqref{eq:caseB}, $-b_{r,\mu}$ is convex,
and $H_\theta(\xi)>0$ is affine in $\theta$ for $\xi>0$; hence the last
term in \eqref{eq:two-level-action} is convex in $\theta$, and so is
$S(u_\theta)$. If $u_0\ne u_1$, it is not affine, being the sum of a
convex function and the function $\frac12T_s(u_\theta)$, which is convex
and not affine by (ii).
\end{proof}

\begin{lemma}\label{lem:two-level-endpoint}
Let $u_0,u_1\in\mathcal A_M$, assume that $u_0$ is continuous, and let
$u_\theta$ be given by \eqref{eq:two-level}. Then
$\|u_\theta-u_0\|_{L^2(\R)}^2=o(\theta)$ as $\theta\to0^+$.
\end{lemma}
\begin{proof}
Write $u_\theta^{\mathrm{lo}}=\min\{u_\theta,\mu\}$ and
$u_\theta^{\mathrm{up}}=(u_\theta-\mu)_+=\rho_\theta^{1/r}$, $\rho_\theta=\rho_{H_\theta}$.
By continuity of $u_0$, $\{u_0>t\}=(-R_0(t),R_0(t))$ for every
$t\in(0,\mu)$, and Lemma~\ref{lem:height-coordinates}(iii) applies at
every point, so the right side of \eqref{eq:two-level} with $\theta=0$ equals $u_0$, that is,
$u_0^{\mathrm{lo}}=\min\{u_0,\mu\}$ and $\rho_0=(u_0-\mu)_+^r$
everywhere. Since $u_\theta-u_0=(u_\theta^{\mathrm{lo}}-u_0^{\mathrm{lo}})+(u_\theta^{\mathrm{up}}-u_0^{\mathrm{up}})$,
\begin{equation*}
 \|u_\theta-u_0\|_2^2\le2\|u_\theta^{\mathrm{lo}}-u_0^{\mathrm{lo}}\|_2^2
 +2\|u_\theta^{\mathrm{up}}-u_0^{\mathrm{up}}\|_2^2,
\end{equation*}
and, by Hölder's inequality and, for the upper part, the elementary
inequality $|a^{1/r}-b^{1/r}|\le|a-b|^{1/r}$ for $a,b\ge0$ and $r\ge1$,
\begin{gather*}
 \|u_\theta^{\mathrm{lo}}-u_0^{\mathrm{lo}}\|_2^2\le\|u_\theta^{\mathrm{lo}}-u_0^{\mathrm{lo}}\|_\infty
 \|u_\theta^{\mathrm{lo}}-u_0^{\mathrm{lo}}\|_1,\\
 \|u_\theta^{\mathrm{up}}-u_0^{\mathrm{up}}\|_2^2\le\int_\R|\rho_\theta-\rho_0|^{2/r}\dd x
 \le\|\rho_\theta-\rho_0\|_\infty^{2/r-1}\|\rho_\theta-\rho_0\|_1,
\end{gather*}
where $2/r-1>0$ because $r\le1/(1-s)<2$. It therefore suffices to prove
\begin{equation}\label{eq:two-level-four}
 \|u_\theta^{\mathrm{lo}}-u_0^{\mathrm{lo}}\|_1=O(\theta),\quad
 \|u_\theta^{\mathrm{lo}}-u_0^{\mathrm{lo}}\|_\infty\to0,\quad
 \|\rho_\theta-\rho_0\|_1=O(\theta),\quad
 \|\rho_\theta-\rho_0\|_\infty\to0\quad\text{as }\theta\to0^+.
\end{equation}

\emph{The lower part.} For $x_0,x_1>0$ the function
$g(\theta)=((1-\theta)x_0^{1-2s}+\theta x_1^{1-2s})^{1/(1-2s)}$ satisfies
\begin{equation*}
 g'(\theta)=\frac{g(\theta)^{2s}\,(x_1^{1-2s}-x_0^{1-2s})}{1-2s},\qquad
 |g'(\theta)|\le\frac{\max\{x_0,x_1\}}{1-2s},
\end{equation*}
since $g(\theta)\le\max\{x_0,x_1\}$ and $|x_1^{1-2s}-x_0^{1-2s}|\le\max\{x_0,x_1\}^{1-2s}$.
Hence $|R_\theta(t)-R_0(t)|\le\theta(1-2s)^{-1}(R_0(t)+R_1(t))$. By
the first formula in Lemma~\ref{lem:two-level}(ii), applied to $u_\theta$
and to $u_0$, and since the two intervals $\{|x|<R_\theta(t)\}$ and $\{|x|<R_0(t)\}$ are
both centered at the origin, so that their indicators differ on a set
of measure $2|R_\theta(t)-R_0(t)|$, 
\begin{equation*}
 \begin{split}
 \|u_\theta^{\mathrm{lo}}-u_0^{\mathrm{lo}}\|_1
 &\le\int_0^\mu\int_\R\bigl|\1_{\{|x|<R_\theta(t)\}}-\1_{\{|x|<R_0(t)\}}\bigr|\dd x\dd t\\
 &=2\int_0^\mu|R_\theta(t)-R_0(t)|\dd t
 \le\frac{\theta\,(\|u_0\|_1+\|u_1\|_1)}{1-2s},
 \end{split}
\end{equation*}
where $\int_0^\mu2R_i(t)\dd t\le\|u_i\|_1$ as in the proof of
Lemma~\ref{lem:two-level}(iii).
For the uniform convergence, $R_\theta\ge((1-\theta)R_0^{1-2s})^{1/(1-2s)}=(1-\theta)^{1/(1-2s)}R_0$
gives $u_\theta^{\mathrm{lo}}(x)\ge u_0^{\mathrm{lo}}(c_\theta x)$ with
$c_\theta=(1-\theta)^{-1/(1-2s)}$. Since $u_0^{\mathrm{lo}}$ is continuous,
nonincreasing in $|x|$ and tends to zero at infinity, for $\varepsilon>0$
there is $L$ with $u_0^{\mathrm{lo}}\le\varepsilon$ on $\{|x|\ge L\}$, while for
$|x|\le L$ the points $x$ and $c_\theta x$ are at distance at most
$(c_\theta-1)L\to0$; by uniform continuity of $u_0^{\mathrm{lo}}$ on $\R$,
\begin{equation*}
 \sup_x\,(u_0^{\mathrm{lo}}-u_\theta^{\mathrm{lo}})_+
 \le\sup_x\bigl(u_0^{\mathrm{lo}}(x)-u_0^{\mathrm{lo}}(c_\theta x)\bigr)
 \le\varepsilon+\sup_{|x|\le L}\bigl|u_0^{\mathrm{lo}}(x)-u_0^{\mathrm{lo}}(c_\theta x)\bigr|
 \longrightarrow\varepsilon\quad\text{as }\theta\to0^+,
\end{equation*}
and $\varepsilon$ is arbitrary.
Conversely, fix $\delta\in(0,\mu)$. For $t\in[\delta,\mu)$ we have
$R_1(t)\le R_1(\delta)$ and $R_0(t)\le R_0(\delta)$, hence
\begin{gather*}
 R_\theta(t)\le\bigl(R_0(t)^{1-2s}+\theta R_1(\delta)^{1-2s}\bigr)^{1/(1-2s)}
 \le R_0(t)+e_\delta(\theta),\\
 e_\delta(\theta)=\sup_{0\le z\le R_0(\delta)}
 \Bigl[\bigl(z^{1-2s}+\theta R_1(\delta)^{1-2s}\bigr)^{1/(1-2s)}-z\Bigr],
\end{gather*}
and $e_\delta(\theta)\to0$ as $\theta\to0^+$ by uniform continuity on the
compact interval. Thus, at a point $x$, the layers $t\in[\delta,\mu)$
with $|x|<R_\theta(t)$ satisfy $|x|-e_\delta(\theta)<R_0(t)$, so
$u_\theta^{\mathrm{lo}}(x)\le\delta+u_0^{\mathrm{lo}}((|x|-e_\delta(\theta))_+)$, and the
uniform continuity of $u_0^{\mathrm{lo}}$ gives
$\limsup_{\theta\to0^+}\sup_x(u_\theta^{\mathrm{lo}}-u_0^{\mathrm{lo}})_+\le\delta$. Since
$\delta$ is arbitrary, $\|u_\theta^{\mathrm{lo}}-u_0^{\mathrm{lo}}\|_\infty\to0$  as $\theta\to0^+$.

\emph{The upper part.} Let $h_*=\min\{H_0'(0),H_1'(0)\}>0$, so that
$H_\theta'\ge h_*$ on $[0,M)$ because $H_i'$ is nondecreasing. From \eqref{eq:reconstruction},
\begin{equation*}
 \begin{split}
 |\rho_\theta(x)-\rho_0(x)|
 &\le\int_0^M\bigl|\1_{\{|x|<R_{H_\theta}(\xi)\}}H_\theta'(\xi)
 -\1_{\{|x|<R_{H_0}(\xi)\}}H_0'(\xi)\bigr|\dd\xi\\
 &\le\int_0^M\Bigl(|H_\theta'-H_0'|\1_{\{|x|<R_{H_\theta}\}}
 +H_0'\,\bigl|\1_{\{|x|<R_{H_\theta}\}}-\1_{\{|x|<R_{H_0}\}}\bigr|\Bigr)\dd\xi,
 \end{split}
\end{equation*}
and integrating in $x$, with $|H_\theta'-H_0'|=\theta|H_1'-H_0'|$ and
$|2R_{H_\theta}-2R_{H_0}|=|H_\theta'-H_0'|/(H_\theta'H_0')$,
\begin{equation*}
 \begin{split}
 \|\rho_\theta-\rho_0\|_1
 &\le\int_0^M\Bigl(|H_\theta'-H_0'|\,2R_{H_\theta}
 +H_0'\,\bigl|2R_{H_\theta}-2R_{H_0}\bigr|\Bigr)\dd\xi\\
 &=2\theta\int_0^M\frac{|H_1'-H_0'|}{H_\theta'}\dd\xi
 \le\frac{2\theta}{h_*}\bigl(H_0(M)+H_1(M)\bigr).
 \end{split}
\end{equation*}
For the uniform estimate, $H_\theta'\ge(1-\theta)H_0'$ gives
$R_{H_\theta}\le R_{H_0}/(1-\theta)$, hence
$\1_{\{|x|<R_{H_\theta}(\xi)\}}\le\1_{\{(1-\theta)|x|<R_{H_0}(\xi)\}}$ and
\begin{equation*}
 \rho_\theta(x)\le\int_0^M\1_{\{(1-\theta)|x|<R_{H_0}(\xi)\}}
 \bigl((1-\theta)H_0'(\xi)+\theta H_1'(\xi)\bigr)\dd\xi
 \le(1-\theta)\rho_0((1-\theta)x)+\theta H_1(M).
\end{equation*}
Conversely, let $\Lambda>1$ and $E_\Lambda=\{\xi\in(0,M):H_1'(\xi)>\Lambda H_0'(\xi)\}$.
Then $\int_{E_\Lambda}H_0'\dd\xi\le\Lambda^{-1}\int_0^MH_1'\dd\xi=\Lambda^{-1}H_1(M)$,
and on $(0,M)\setminus E_\Lambda$ we have
$H_\theta'\le(1+\theta(\Lambda-1))H_0'$, that is,
$R_{H_\theta}\ge R_{H_0}/(1+\theta(\Lambda-1))$. Keeping only the
contribution of $(1-\theta)H_0'$ on $(0,M)\setminus E_\Lambda$,
with $c_\Lambda=1+\theta(\Lambda-1)$,
\begin{equation*}
 \begin{split}
 \rho_\theta(x)&\ge(1-\theta)\int_{(0,M)\setminus E_\Lambda}
 \1_{\{c_\Lambda|x|<R_{H_0}(\xi)\}}H_0'(\xi)\dd\xi
 \ge(1-\theta)\Bigl(\rho_0(c_\Lambda x)-\int_{E_\Lambda}H_0'(\xi)\dd\xi\Bigr)\\
 &\ge(1-\theta)\rho_0(c_\Lambda x)-\Lambda^{-1}H_1(M)
 \ge\rho_0(c_\Lambda x)-\theta H_0(M)-\Lambda^{-1}H_1(M),
 \end{split}
\end{equation*}
where we used $\rho_0(c_\Lambda x)=\int_0^M\1_{\{c_\Lambda|x|<R_{H_0}(\xi)\}}H_0'\dd\xi$
from \eqref{eq:reconstruction}, the bound on $\int_{E_\Lambda}H_0'$, and
$0\le\rho_0\le H_0(M)$.
We combine the two bounds. Let $L_0=R_{H_0}(0)$. By
\eqref{eq:superlevels} and the right continuity of $H_0'$ at $0$,
$\{\rho_0>0\}=\bigcup_{0<\xi<M}(-R_{H_0}(\xi),R_{H_0}(\xi))=(-L_0,L_0)$,
so the continuous function $\rho_0$ has compact support $[-L_0,L_0]$
and is uniformly continuous on $\R$, that is, the quantity
$\omega_{\rho_0}(\delta)=\sup\{|\rho_0(x)-\rho_0(y)|:|x-y|\le\delta\}$ tends to
zero as $\delta\to0^+$. Likewise $\rho_\theta$ vanishes for
$|x|\ge R_{H_\theta}(0)=1/(2H_\theta'(0))\le L_0/(1-\theta)$. The upper
bound gives $\rho_\theta(x)-\rho_0(x)\le\rho_0((1-\theta)x)-\rho_0(x)+\theta H_1(M)$;
here both $\rho_\theta(x)$ and $\rho_0(x)$ vanish if $|x|\ge L_0/(1-\theta)$,
and otherwise $|(1-\theta)x-x|=\theta|x|\le\theta L_0/(1-\theta)$. The lower
bound gives
$\rho_\theta(x)-\rho_0(x)\ge\rho_0((1+\theta(\Lambda-1))x)-\rho_0(x)-\theta H_0(M)-\Lambda^{-1}H_1(M)$;
here $\rho_0(x)=0\le\rho_\theta(x)$ if $|x|\ge L_0$, and otherwise the
two arguments of $\rho_0$ are at distance $\theta(\Lambda-1)|x|\le\theta(\Lambda-1)L_0$.
Consequently
\begin{equation*}
 \|\rho_\theta-\rho_0\|_\infty
 \le\omega_{\rho_0}\Bigl(\frac{\theta L_0}{1-\theta}\Bigr)+\omega_{\rho_0}\bigl(\theta(\Lambda-1)L_0\bigr)
 +\theta\bigl(H_0(M)+H_1(M)\bigr)+\Lambda^{-1}H_1(M),
\end{equation*}
and letting first $\theta\to0^+$ and then $\Lambda\to\infty$ gives
$\|\rho_\theta-\rho_0\|_\infty\to0$. This proves \eqref{eq:two-level-four}.
\end{proof}

\section{Constrained coercivity of the second variation}\label{sec:coercivity}

Convexity along an interpolation shows that $u$ is a strict minimizer
only together with information at the endpoint $\theta=0$. In case
\eqref{eq:caseB} this information is given by
Lemma~\ref{lem:two-level-endpoint}. In case \eqref{eq:caseA}, where the
height interpolation need not be differentiable in $H^s$ at $\theta=0$
and $\lambda-f'(u)$ may be unbounded, it is given by the coercivity of
the second variation $\mathcal Q_u$ in even directions tangent to
$\{\mathcal N=\mathcal N(u)\}$, which is stable under local
perturbations of the potential. In both cases, the strict positivity of
$\mathcal Q_u$ in these directions 
also gives the nondegeneracy in Section~\ref{sec:nondegeneracy}.

\smallskip
For a real measurable potential $W$ bounded below, define
\begin{equation*}
 D(W)=\left\{\phi\in H^s(\R):\int_\R W_+\phi^2\dd x<\infty\right\},
 \qquad
 \|\phi\|_{D(W)}^2=\|\phi\|_{H^s}^2+\int_\R W_+\phi^2\dd x,
\end{equation*}
where $W_\pm=\max\{\pm W,0\}$. Since $W_-\in L^\infty$, for $\phi\in D(W)$
we have $\int_\R|W|\phi^2\dd x\le\int_\R W_+\phi^2\dd x+\|W_-\|_\infty\|\phi\|_2^2<\infty$,
so the bilinear form
$\mathcal Q_W[\phi,\psi]=T_s(\phi,\psi)+\int_\R W\phi\psi\dd x$ is finite on
$D(W)$; we write $\mathcal Q_W[\phi]=\mathcal Q_W[\phi,\phi]$, and Cauchy--Schwarz in each of the
three terms gives
\begin{equation}\label{eq:form-continuity-general}
 |\mathcal Q_W[\phi,\psi]|\le(1+\|W_-\|_\infty)\|\phi\|_{D(W)}\|\psi\|_{D(W)}
 \qquad(\phi,\psi\in D(W)).
\end{equation}
If $W$ is locally bounded, then $C_c^\infty(\R)$ is dense in $D(W)$, and
even (odd) functions in $C_c^\infty(\R)$ are dense in the even (odd) part
of $D(W)$. Indeed, let $\phi\in D(W)$. The cutoffs $\chi_R$ of
Appendix~\ref{app:prelim} satisfy $\chi_R\phi\to\phi$ in $H^s$ by
\eqref{eq:cutoff-approximation}, and
\begin{equation*}
 \int_\R W_+|\chi_R\phi-\phi|^2\dd x\le\int_{|x|>R}W_+\phi^2\dd x\longrightarrow0
 \quad\text{as }R\to\infty;
\end{equation*}
for fixed $R$, mollification with an even mollifier approximates
$\chi_R\phi$ in $H^s$ by a sequence $\{v_j\}_{j\ge1}$ of smooth functions supported in
$[-2R-1,2R+1]$, and
$\int_\R W_+|v_j-\chi_R\phi|^2\dd x\le\|W_+\|_{L^\infty([-2R-1,2R+1])}\|v_j-\chi_R\phi\|_2^2\to0$.
Both operations preserve parity, and a diagonal choice proves the claim.

\smallskip
In this section $u$ is a symmetric decreasing positive solution under the
hypotheses of Theorem~\ref{thm:main}, in either case, and
\begin{equation*}
 V=\lambda-f'(u),\qquad D_u=D(V),\qquad \mathcal Q_u=\mathcal Q_V,\qquad
 M=\mathcal N(u),\qquad \rho=(u-\mu)_+^r.
\end{equation*}
The form $\mathcal Q_u$ is the second variation of $S$ at $u$ in the
following sense. Along the smooth curves $\sigma\mapsto u_\sigma$
constructed in Lemmas~\ref{lem:realization}
and~\ref{lem:small-realization} below, one has
$\frac{d^2}{d\sigma^2}S(u_\sigma)|_{\sigma=0}
=\mathcal Q_u[\partial_\sigma u_\sigma|_{\sigma=0}]$ (Step~1 of the
proofs of Lemmas~\ref{lem:local-lower-bound} and~\ref{lem:small-local}).
Since $u$ is symmetric decreasing, smooth and positive
(Lemmas~\ref{lem:solution-estimates} and~\ref{lem:symmetry}), $V$ is even
and locally bounded. By \eqref{eq:f-assumptions} there are
$t_0,\nu_0>0$ with $f'(t)\le\lambda-\nu_0$ for $0<t\le t_0$, and since
$u\in C_0(\R)$ there is $R_0$ with $u\le t_0$ on $\{|x|\ge R_0\}$; with
$C_V=\sup_{0<t\le\|u\|_\infty}f'(t)-\lambda$ this gives
\begin{equation}\label{eq:linearized-potential-bounds}
 V\ge-C_V\quad\text{on }\R,\qquad
 V\ge\nu_0\quad\text{on }\R\setminus[-R_0,R_0].
\end{equation}
In case \eqref{eq:caseB} let $R_*>0$ be the point with $u(R_*)=\mu$,
which exists by Lemma~\ref{lem:solution-estimates}(v) and
Lemma~\ref{lem:symmetry}, and let $Z=\{-R_*,0,R_*\}$. The constraint
functional is
\begin{equation}\label{eq:constraint-functional}
 \ell(\phi)=\int_\R q_u\phi\dd x,\qquad
 q_u=\begin{cases}
 u^{r-1}&\text{in case \eqref{eq:caseA}},\\
 (u-\mu)_+^{r-1}&\text{in case \eqref{eq:caseB}, }r>1,\\
 \1_{(-R_*,R_*)}&\text{in case \eqref{eq:caseB}, }r=1,
 \end{cases}
\end{equation}
so that $r\ell(\phi)=\frac{d}{d\varepsilon}\mathcal N(u+\varepsilon\phi)|_{\varepsilon=0}$
is the derivative of $\mathcal N$ at $u$ in the direction $\phi$. The weight $q_u$
is even, nonnegative, bounded and in $L^2(\R)$: in case
\eqref{eq:caseA}, $\|u^{r-1}\|_2^2\le\|u\|_\infty^{2r-4}\|u\|_2^2$ since
$r\ge2$, and in case \eqref{eq:caseB} it has compact support. Moreover
$q_u=\rho^\beta$ on $\{u>\mu\}$, and $\ell(u)>0$. The density $\rho$ is continuous, even, positive and smooth on
$(-R_*,R_*)$, vanishes outside $(-R_*,R_*)$, and satisfies $\rho'<0$ on
$(0,R_*)$, where $R_*=\infty$ in case \eqref{eq:caseA}. Hence
Lemma~\ref{lem:reconstruction}(i) applies. The coordinate
$\xi(x)=G(\rho(x))$ is a smooth decreasing bijection of $(0,R_*)$ onto
$(0,M)$, the height function $H$ of $\rho$ is smooth on $(0,M)$ and satisfies
\eqref{eq:mass-coordinate} and \eqref{eq:height-derivatives} for
$0<x<R_*$, $H$ is admissible in case \eqref{eq:caseA}, and
$H'(0)=1/(2R_*)>0$ in case \eqref{eq:caseB}.
For $\phi\in L^2(\R)$ we define
\begin{equation}\label{eq:e-phi}
 E_\phi(x)=\frac rx\int_0^xq_u(y)\phi(y)\dd y\quad(0<x<R_*),\qquad
 \eta_\phi(\xi(x))=E_\phi(x).
\end{equation}
The function $E_\phi$ is locally absolutely continuous, and
differentiating $xE_\phi(x)=r\int_0^xq_u\phi\dd y$ gives
\begin{equation}\label{eq:e-derivative}
 xE_\phi'(x)=rq_u(x)\phi(x)-E_\phi(x)\quad\text{for a.e. }0<x<R_*.
\end{equation}
In case \eqref{eq:caseB}, since $q_u\in L^\infty$ and $\phi\in L^2$,
$E_\phi$ extends continuously to $x=R_*$. If $\phi$ is even, then
$\int_0^xq_u\phi\dd y\to\frac12\ell(\phi)$ as $x\to R_*$, so in case
\eqref{eq:caseA} $E_\phi$ vanishes for $x\ge L$ if $\ell(\phi)=0$ and
$\supp\phi\subset(-L,L)$, and in case \eqref{eq:caseB}
$E_\phi(R_*)=\frac r{2R_*}\ell(\phi)$. Finally, let
\begin{equation*}
 \mathcal T=\begin{cases}
 \{\phi\in C_c^\infty(\R):\ \phi\text{ even},\ \ell(\phi)=0\}
 &\text{in case \eqref{eq:caseA}},\\[2pt]
 \{\phi\in C_c^\infty(\R\setminus Z):\ \phi\text{ even},\ \ell(\phi)=0\}
 &\text{in case \eqref{eq:caseB}}.
 \end{cases}
\end{equation*}

\begin{theorem}\label{thm:coercivity}
Assume the hypotheses of Theorem~\ref{thm:main}, and let $u$ be a
symmetric decreasing positive solution. Set $V_0=V=\lambda-f'(u)$, and let $\{V_\theta\}_{0<\theta<\theta_*}$
be a family of even real measurable potentials such that, for constants
$C_V,R_0,\nu_0>0$ independent of $\theta$,
\begin{equation}\label{eq:secant-structure}
 V_\theta\ge-C_V\quad\text{a.e. in }\R,
 \qquad V_\theta(x)\ge\nu_0\quad\text{for a.e. }x\text{ with }|x|\ge R_0,
\end{equation}
and, for every $R>0$,
\begin{equation}\label{eq:potential-convergence}
 \lim_{\theta\to0^+}\|V_\theta-V_0\|_{L^\infty([-R,R])}=0.
\end{equation}
Then there are $c,C>0$ and $\bar\theta\in(0,\theta_*)$ such that, for every
$0\le \theta<\bar\theta$ and every even $\phi\in D(V_\theta)$,
\begin{equation}\label{eq:penalized-coercivity}
 \mathcal Q_{V_\theta}[\phi]\ge c\|\phi\|_{H^s}^2-C|\ell(\phi)|^2.
\end{equation}
In particular, $\mathcal Q_u[\phi]\ge c\|\phi\|_{H^s}^2$ for every even
$\phi\in D_u$ with $\ell(\phi)=0$. Moreover, with $\delta=\sqrt{c/(2C)}$,
for every $0\le\theta<\bar\theta$ and every even $\phi\in D(V_\theta)$,
\begin{equation}\label{eq:almost-tangent-coercivity}
 \mathcal Q_{V_\theta}[\phi]\ge\frac c2\|\phi\|_{H^s}^2
 \qquad\text{if }|\ell(\phi)|\le\delta\|\phi\|_{H^s}.
\end{equation}
\end{theorem}

The constant family $\{V_\theta\}_{0<\theta<\theta_*}$, $V_\theta=V$, satisfies the hypotheses of Theorem~\ref{thm:coercivity};
the perturbed potentials are needed in case \eqref{eq:caseA} of
Lemma~\ref{lem:own-mass}, where $V_\theta$ will be the potential of the
exact remainder of Lemma~\ref{lem:remainder} below.

%\subsection*{Smooth tangent directions in case \texorpdfstring{\eqref{eq:caseA}}{(f3)}}

\begin{lemma}\label{lem:realization}
In case \eqref{eq:caseA}, let $\phi\in\mathcal T$ and $\eta=\eta_\phi$.
Then $\eta$ is smooth on $(0,M)$, continuous on $[0,M]$ and zero near
$\xi=0$, and for $|\sigma|$ small the height function $H_\sigma=H+\sigma\eta$ is
admissible, and its reconstruction $u_\sigma=\rho_{H_\sigma}^{1/r}$ satisfies
\begin{equation*}
 u_0=u,\qquad \int_\R u_\sigma^r\dd x=M,
 \qquad \left.\partial_\sigma u_\sigma\right|_{\sigma=0}=\phi.
\end{equation*}
The curve $\sigma\mapsto u_\sigma$ is $C^2$ in $H^s(\R)$, and $u_\sigma-u$ is smooth
with a common compact support.
\end{lemma}
\begin{proof}
\emph{Step 1. The height direction.}
Write $E=E_\phi$. The substitution $y=ax$ in \eqref{eq:e-phi} gives
\begin{equation}\label{eq:tangent-primitive}
 E(x)=r\int_0^1u(ax)^{r-1}\phi(ax)\dd a,\qquad
 xE(x)=r\int_0^xu^{r-1}\phi\dd y,\qquad E(x)+xE'(x)=ru(x)^{r-1}\phi(x),
\end{equation}
so $E$ is even and smooth. Choose $L>0$ with $\supp\phi\subset(-L,L)$;
since $\phi$ is even and $\ell(\phi)=0$, $\int_0^\infty u^{r-1}\phi\dd y=\frac12\ell(\phi)=0$,
hence $E(x)=0$ for $|x|\ge L$. As $\xi$ decreases from $M$ to $0$ when
$x$ increases from $0$ to $\infty$, $\eta$ is smooth on $(0,M)$,
continuous on $[0,M]$, and zero near $\xi=0$.

\emph{Step 2. Admissibility.}
By Lemma~\ref{lem:symmetry}, $\rho'<0$ on $(0,\infty)$ and
$\rho''(0)=ru(0)^{r-1}u''(0)<0$. Define $\kappa=E'/\rho'$ away from zero.
Since $E'$ and $\rho'$ are odd with $E'(0)=\rho'(0)=0$,
\begin{equation*}
 \kappa(x)=\frac{E'(x)}{\rho'(x)}
 =\frac{\int_0^1E''(ax)\dd a}{\int_0^1\rho''(ax)\dd a}\qquad(x\ne0),
\end{equation*}
and the denominator equals $\rho''(0)<0$ at $x=0$, so this formula
extends $\kappa$ to an even smooth function with compact support and
$\kappa(0)=E''(0)/\rho''(0)$. For
$|\sigma|(\|\kappa\|_\infty+\|x\kappa'\|_\infty)<1/2$ define
\begin{equation*}
 y_\sigma(x)=\frac{x}{1+\sigma\kappa(x)},\qquad
 y_\sigma'(x)=\frac{1+\sigma\kappa(x)-\sigma x\kappa'(x)}{(1+\sigma\kappa(x))^2}>0;
\end{equation*}
this odd map is a smooth diffeomorphism of $\R$, equal to the identity
outside $[-L,L]$. Differentiating $\eta(\xi(x))=E(x)$ and using
\eqref{eq:mass-coordinate} and \eqref{eq:height-derivatives},
\begin{equation*}
 \eta'(\xi(x))=\frac{E'(x)}{\xi'(x)}=\frac{E'(x)}{2x\rho'(x)}=\frac{\kappa(x)}{2x},
 \qquad
 H_\sigma'(\xi(x))=\frac{1+\sigma\kappa(x)}{2x}=\frac1{2y_\sigma(x)},
\end{equation*}
and differentiating the last identity in $x$,
\begin{equation*}
 H_\sigma''(\xi(x))\,\xi'(x)=-\frac{y_\sigma'(x)}{2y_\sigma(x)^2}<0,
\end{equation*}
so $H_\sigma''>0$ because $\xi'<0$. Also $H_\sigma'>0$, $H_\sigma(0)=0$, and
$H_\sigma'(\xi(x))=1/(2y_\sigma(x))$ tends to $0$ as $x\to\infty$ and to $\infty$
as $x\to0^+$, which are the endpoint limits in \eqref{eq:admissible}.
Thus $H_\sigma$ is admissible, and since $R_{H_\sigma}(\xi(x))=1/(2H_\sigma'(\xi(x)))=y_\sigma(x)$,
Lemma~\ref{lem:reconstruction}(ii) gives
\begin{equation}\label{eq:varied-reconstruction}
 \rho_\sigma(y_\sigma(x))=H_\sigma(\xi(x))=\rho(x)+\sigma E(x),
 \qquad \int_\R\rho_\sigma\dd x=M,\qquad \rho_\sigma=\rho_{H_\sigma}.
\end{equation}

\emph{Step 3. Smoothness and the tangent.}
Since $y_\sigma^{-1}$ depends smoothly on $(\sigma,x)$,
\eqref{eq:varied-reconstruction} shows that
$\rho_\sigma(x)=\rho(y_\sigma^{-1}(x))+\sigma E(y_\sigma^{-1}(x))$ is smooth in $(\sigma,x)$ and
equals $\rho$ outside $[-L,L]$; after reducing $|\sigma|$,
$\rho_\sigma\ge\frac12\min_{[-L,L]}\rho>0$ on $[-L,L]$, so $u_\sigma=\rho_\sigma^{1/r}$ is
smooth in $(\sigma,x)$, and $u_\sigma-u$ together with all its $\sigma$-derivatives has support
in $[-L,L]$. For $v\in C_c^1([-L,L])$, Fourier comparison gives
$\|v\|_{H^s}\le\sqrt2\,\|v\|_{H^1}\le C_L\|v\|_{C^1}$, so the continuity of the
$\sigma$-derivatives in $C^1$ makes the curve $C^2$ in $H^s$. Finally,
since $\partial_\sigma y_\sigma|_{\sigma=0}=-x\kappa(x)$, it follows from \eqref{eq:tangent-primitive} and differentiating
\eqref{eq:varied-reconstruction} at $\sigma=0$ that 

\begin{equation*}
 \left.\partial_\sigma\rho_\sigma(x)\right|_{\sigma=0}
 =E(x)-\rho'(x)\,\partial_\sigma y_\sigma(x)|_{\sigma=0}=E(x)+x\kappa(x)\rho'(x)=E(x)+xE'(x)
 =ru(x)^{r-1}\phi(x),
\end{equation*}
and the derivative of the positive $r$th root gives
$\partial_\sigma u_\sigma|_{\sigma=0}=\frac1ru^{1-r}\partial_\sigma\rho_\sigma|_{\sigma=0}=\phi$.
\end{proof}

\begin{lemma}\label{lem:local-lower-bound}
In case \eqref{eq:caseA}, let $\phi\in\mathcal T$ and $\eta=\eta_\phi$.
If $s>\beta$, then for compact intervals $\mathcal I,\mathcal J\subset(0,M)$ of positive
length with $\max\mathcal I<\min\mathcal J$ there is $c_{\mathcal I,\mathcal J}>0$, depending only on
$u,\mathcal I,\mathcal J$, such that
\begin{equation}\label{eq:kinetic-local-bound}
 \mathcal Q_u[\phi]\ge c_{\mathcal I,\mathcal J}
             \left(\int_{\mathcal I}\eta^2\dd\xi+\int_{\mathcal J}\eta^2\dd\xi\right).
\end{equation}
If $s=\beta$, then for every compact interval $\mathcal I\subset(0,M)$,
\begin{equation}\label{eq:potential-local-bound}
 \mathcal Q_u[\phi]\ge\frac1r\min_{\xi\in \mathcal I}[-b_r''(H(\xi))]
                                  \int_{\mathcal I}\eta^2\dd\xi.
\end{equation}
\end{lemma}
\begin{proof}
\emph{Step 1. The second-variation identity.}
Let $u_\sigma$ be the curve of Lemma~\ref{lem:realization} and
$\psi=\partial_\sigma^2u_\sigma|_{\sigma=0}\in C_c^\infty(\R)$; both $\phi$ and $\psi$
belong to $D_u$. Since $u_\sigma-u$ is smooth with common compact support,
on which $u_\sigma$ has a positive lower bound, the local terms of
$S(u_\sigma)-S(u)$ may be differentiated twice under the integral, and the
kinetic term is $C^2$ because $u_\sigma$ is $C^2$ in $H^s$. Writing
$\dot u_\sigma=\partial_\sigma u_\sigma$ and $\ddot u_\sigma=\partial_\sigma^2u_\sigma$, the derivatives
of $S(u_\sigma)=\frac12T_s(u_\sigma)+\int_\R(\frac\lambda2u_\sigma^2-F(u_\sigma))\dd x$ are
\begin{align*}
 \frac{d}{d\sigma}S(u_\sigma)&=T_s(u_\sigma,\dot u_\sigma)+\int_\R\bigl(\lambda u_\sigma-f(u_\sigma)\bigr)\dot u_\sigma\dd x,\\
 \frac{d^2}{d\sigma^2}S(u_\sigma)&=T_s(\dot u_\sigma)+T_s(u_\sigma,\ddot u_\sigma)
 +\int_\R\bigl(\lambda-f'(u_\sigma)\bigr)\dot u_\sigma^2\dd x
 +\int_\R\bigl(\lambda u_\sigma-f(u_\sigma)\bigr)\ddot u_\sigma\dd x.
\end{align*}
At $\sigma=0$ we have $\dot u_0=\phi$ and $\ddot u_0=\psi$, and the two terms
containing $\psi$ add up to $T_s(u,\psi)+\int_\R(\lambda u-f(u))\psi\dd x=0$
by \eqref{eq:weak-all-tests}. Hence
\begin{equation*}
 \left.\frac{d^2}{d\sigma^2}S(u_\sigma)\right|_{\sigma=0}
 =T_s(\phi)+\int_\R(\lambda-f'(u))\phi^2\dd x
 =\mathcal Q_u[\phi].
\end{equation*}
On the other hand, \eqref{eq:height-action} applies along the affine
height curve, since $H_\sigma$ is admissible, $u_\sigma\in H^s$ and
$F(u_\sigma)\in L^1$. Moreover, on the support of $\eta$ all values $H+\sigma\eta$
lie in a compact subset of $(0,\infty)$, so
$|b_r''(H+\sigma\eta)\eta^2|\le C\eta^2\in L^1(0,M)$, and differentiating
$-\frac1r\int_0^Mb_r(H+\sigma\eta)\dd\xi$ twice under the integral gives
\begin{equation}\label{eq:height-second-variation}
 \mathcal Q_u[\phi]
 =\frac12\left.\frac{d^2}{d\sigma^2}T_s(u_\sigma)\right|_{\sigma=0}
       +\frac1r\int_0^M[-b_r''(H(\xi))]\eta(\xi)^2\dd\xi,
\end{equation}
where both terms are nonnegative by Lemma~\ref{lem:kinetic-convexity}
and \eqref{eq:caseA}.

\emph{Step 2. The kinetic bound when $s>\beta$.}
Let
\begin{equation*}
 T_{\mathcal I,\mathcal J}(\sigma)=C_{s,r}\int_{\mathcal I}\int_{\mathcal J}
 \frac{H_\sigma'(\xi)^{1+2s}}{H_\sigma(\xi)^\beta H_\sigma(\zeta)^\beta}\dd\zeta\dd\xi
\end{equation*}
be the contribution of the $k=0$ term of \eqref{eq:kinetic-series} on
$\mathcal I\times\mathcal J$, which is contained in $\{\xi<\zeta\}$
because $\max\mathcal I<\min\mathcal J$. All other terms of the series,
and the $k=0$ term outside $\mathcal I\times\mathcal J$, have integrands
which are convex in $\sigma$ by Lemma~\ref{lem:monomial}, since their
arguments $H_\sigma'(\xi)$, $H_\sigma'(\zeta)$, $H_\sigma(\xi)$,
$H_\sigma(\zeta)$ are positive and affine in $\sigma$, as in the proof of
Lemma~\ref{lem:kinetic-convexity}. And their sum is finite
since it is bounded by $T_s(u_\sigma)<\infty$; hence $T_s(u_\sigma)-T_{\mathcal I,\mathcal J}(\sigma)$ is
convex in $\sigma$. Both $T_s(u_\sigma)$ and $T_{\mathcal I,\mathcal J}(\sigma)$ are $C^2$, the latter by
differentiation under the integral on the compact rectangle, so for
small $\sigma>0$,
\begin{equation*}
 \frac{T_s(u_\sigma)+T_s(u_{-\sigma})-2T_s(u_0)}{\sigma^2}
 \ge\frac{T_{\mathcal I,\mathcal J}(\sigma)+T_{\mathcal I,\mathcal J}(-\sigma)-2T_{\mathcal I,\mathcal J}(0)}{\sigma^2},
\end{equation*}
and letting $\sigma\to0^+$,
\begin{equation}\label{eq:localized-kinetic-derivative}
 \left.\frac{d^2}{d\sigma^2}T_s(u_\sigma)\right|_{\sigma=0}\ge T_{\mathcal I,\mathcal J}''(0).
\end{equation}
Along the affine curve the arguments $H_\sigma'(\xi)$, $H_\sigma(\xi)$,
$H_\sigma(\zeta)$ have the $\sigma$-derivatives $\eta'(\xi)$, $\eta(\xi)$,
$\eta(\zeta)$, so \eqref{eq:strict-monomial} gives, pointwise on
$\mathcal I\times\mathcal J$,
\begin{equation*}
 \frac{d^2}{d\sigma^2}\Bigl[\frac{H_\sigma'(\xi)^{1+2s}}{H_\sigma(\xi)^\beta H_\sigma(\zeta)^\beta}\Bigr]_{\sigma=0}
 \ge\frac{\beta(s-\beta)}s\,
 \frac{H'(\xi)^{1+2s}}{H(\xi)^\beta H(\zeta)^\beta}
 \left(\frac{\eta(\xi)^2}{H(\xi)^2}+\frac{\eta(\zeta)^2}{H(\zeta)^2}\right).
\end{equation*}
Integrating, and combining with \eqref{eq:localized-kinetic-derivative}
and \eqref{eq:height-second-variation},
\begin{equation*}
 \mathcal Q_u[\phi]\ge\frac12T_{\mathcal I,\mathcal J}''(0)\ge\frac{C_{s,r}\beta(s-\beta)}{2s}
 \int_{\mathcal I}\int_{\mathcal J}
 \frac{H'(\xi)^{1+2s}}{H(\xi)^\beta H(\zeta)^\beta}
 \left(\frac{\eta(\xi)^2}{H(\xi)^2}
             +\frac{\eta(\zeta)^2}{H(\zeta)^2}\right)\dd\zeta\dd\xi.
\end{equation*}
The factor $H'(\xi)^{1+2s}/(H(\xi)^\beta H(\zeta)^\beta)$ has a positive
minimum on $\mathcal I\times\mathcal J$, and $H\le\max_{\mathcal I}H$ on $\mathcal I$, $H\le\max_{\mathcal J}H$ on $\mathcal J$;
integrating in the other variable gives \eqref{eq:kinetic-local-bound}
with
\begin{equation*}
 c_{\mathcal I,\mathcal J}=\frac{C_{s,r}\beta(s-\beta)}{2s}
 \min_{(\xi,\zeta)\in\mathcal I\times\mathcal J}\frac{H'(\xi)^{1+2s}}{H(\xi)^\beta H(\zeta)^\beta}
 \cdot\min\Bigl\{\frac{|\mathcal J|}{(\max_{\mathcal I}H)^2},\frac{|\mathcal I|}{(\max_{\mathcal J}H)^2}\Bigr\}>0,
\end{equation*}
which does not depend on $\phi$.

\emph{Step 3. The potential bound when $s=\beta$.}
Since $s=\beta$ means $r=1/(1-s)$, the strictness in \eqref{eq:caseA}
gives $b_r''<0$ on $(0,\infty)$. The kinetic term in
\eqref{eq:height-second-variation} is nonnegative, and so is the
integrand $-b_r''(H)\eta^2$ of the potential term. Dropping the former
and restricting the latter to $\mathcal I$ gives
\eqref{eq:potential-local-bound}, where
$\min_{\xi\in\mathcal I}[-b_r''(H(\xi))]>0$ because $-b_r''\circ H$ is
continuous and positive on the compact interval $\mathcal I$.
\end{proof}

%\subsection*{Smooth tangent directions in case \texorpdfstring{\eqref{eq:caseB}}{(f4)}}

\begin{lemma}\label{lem:small-realization}
In case \eqref{eq:caseB}, let $\phi\in\mathcal T$, and define $\Upsilon$
on $(0,\mu)$ and $\eta$ on $(0,M)$ by
\begin{equation}\label{eq:lower-direction}
 \Upsilon(u(x))=\frac{(1-2s)\,\phi(x)}{x^{2s}|u'(x)|}\quad(x>R_*),\qquad
 \eta(\xi(x))=E_\phi(x)\quad(0<x<R_*),
\end{equation}
that is, $\eta=\eta_\phi$ as in \eqref{eq:e-phi}. Then
$\Upsilon\in C_c^\infty((0,\mu))$ and $\eta\in C_c^\infty((0,M))$, and
there are $\sigma_0>0$, a compact set $K_\phi\subset\R\setminus Z$ and
a family $\{u_\sigma\}_{|\sigma|<\sigma_0}\subset\mathcal A_M$ with the
following properties.
\begin{enumerate}
\item[(i)] $u_0=u$, and $u_\sigma=u$ outside $K_\phi$.
\item[(ii)] The map $(\sigma,x)\mapsto u_\sigma(x)$ is smooth on
$(-\sigma_0,\sigma_0)\times\R$, and
$\partial_\sigma u_\sigma|_{\sigma=0}=\phi$.
\item[(iii)] The functions $R_{u_\sigma}^{1-2s}$ and $H_{u_\sigma}$ of
\eqref{eq:two-level-data} are affine in $\sigma$:
\begin{equation}\label{eq:affine-data}
 R_{u_\sigma}(t)^{1-2s}=R_u(t)^{1-2s}+\sigma\Upsilon(t)\quad(0<t<\mu),\qquad
 H_{u_\sigma}=H+\sigma\eta.
\end{equation}
\end{enumerate}
\end{lemma}

\begin{proof}
Since $\phi\in C_c^\infty(\R\setminus Z)$ is even, there are
$0<a^{\mathrm{up}}<b^{\mathrm{up}}<R_*<a^{\mathrm{lo}}<b^{\mathrm{lo}}$ with
\begin{equation*}
 \phi=0\quad\text{on }(0,a^{\mathrm{up}}]\cup[b^{\mathrm{up}},R_*)\cup(R_*,a^{\mathrm{lo}}]\cup[b^{\mathrm{lo}},\infty).
\end{equation*}

\emph{Step 1. The lower part.}  Since $u$ is smooth with $u'<0$ on $(0,\infty)$ and $u(R_*)=\mu$, we
have $\{u>t\}=(-R_u(t),R_u(t))$ for $0<t<\mu$, so that the function
$R_u$ of \eqref{eq:two-level-data} is the inverse of $u$ on
$(R_*,\infty)$. It is smooth with $R_u'(t)=1/u'(R_u(t))<0$, and by
\eqref{eq:lower-direction},

\begin{equation*}
 \Upsilon(t)=\frac{(1-2s)\,\phi(R_u(t))}{R_u(t)^{2s}|u'(R_u(t))|}
 \qquad(0<t<\mu).
\end{equation*}
Hence $\Upsilon$ is smooth with support in $[u(b^{\mathrm{lo}}),u(a^{\mathrm{lo}})]\subset(0,\mu)$,
so $\Upsilon\in C_c^\infty((0,\mu))$. Let
\begin{equation*}
 U_\sigma=R_u^{1-2s}+\sigma\Upsilon,\qquad R_\sigma=U_\sigma^{1/(1-2s)}.
\end{equation*}
On the support of $\Upsilon$ we have $R_u\ge a^{\mathrm{lo}}>R_*$, and
$(R_u^{1-2s})'=(1-2s)R_u^{-2s}R_u'$ is negative and bounded away from
zero. Hence there is $\sigma_0>0$ such that, for $|\sigma|<\sigma_0$,
$U_\sigma'<0$ and $U_\sigma>R_*^{1-2s}$ on $(0,\mu)$. Then $R_\sigma$ is
smooth in $(\sigma,t)$, strictly decreasing in $t$, equal to $R_u$
outside the support of $\Upsilon$, and maps $(0,\mu)$ onto
$(R_*,\infty)$. After reducing $\sigma_0$, there are
$R_*<\tilde a^{\mathrm{lo}}<\tilde b^{\mathrm{lo}}$ with $R_\sigma(t)\in[\tilde a^{\mathrm{lo}},\tilde b^{\mathrm{lo}}]$
for all $|\sigma|<\sigma_0$ and all $t$ in the support of $\Upsilon$.

For $|x|>R_*$ let $u_\sigma^{\mathrm{lo}}(x)\in(0,\mu)$ be the unique
solution of $R_\sigma(u_\sigma^{\mathrm{lo}}(x))=|x|$. Since
$\partial_tR_\sigma<0$, the implicit function theorem shows that
$u_\sigma^{\mathrm{lo}}$ is smooth in $(\sigma,x)$. If
$|x|\notin[\tilde a^{\mathrm{lo}},\tilde b^{\mathrm{lo}}]$, then $t=u_\sigma^{\mathrm{lo}}(x)$ lies
outside the support of $\Upsilon$, so $R_u(t)=|x|$ and
$u_\sigma^{\mathrm{lo}}(x)=u(x)$; the same holds for every $|x|>R_*$ if
$\sigma=0$. Finally, for $x>R_*$ and $t=u(x)$, so that $R_u(t)=x$,
differentiating $R_\sigma(u_\sigma^{\mathrm{lo}}(x))=x$ at $\sigma=0$ gives
\begin{equation*}
 \partial_\sigma R_\sigma(t)\big|_{\sigma=0}
 +R_u'(t)\,\partial_\sigma u_\sigma^{\mathrm{lo}}(x)\big|_{\sigma=0}=0,
 \qquad
 \partial_\sigma R_\sigma(t)\big|_{\sigma=0}=\frac{R_u(t)^{2s}\Upsilon(t)}{1-2s}
 =\frac{\phi(x)}{|u'(x)|},
\end{equation*}
by \eqref{eq:lower-direction}. Since $R_u'(t)=1/u'(x)$, we obtain
$\partial_\sigma u_\sigma^{\mathrm{lo}}(x)|_{\sigma=0}=-u'(x)\phi(x)/|u'(x)|=\phi(x)$,
and the same holds for $x<-R_*$ by evenness. This computation is what
determines the formula \eqref{eq:lower-direction} for $\Upsilon$.

\emph{Step 2. The upper part.}
Write $E=E_\phi$. Since $\phi=0$ on $(0,a^{\mathrm{up}}]$, $E=0$ on $(0,a^{\mathrm{up}}]$. For
$b^{\mathrm{up}}\le x<R_*$, since $\phi=0$ on $[b^{\mathrm{up}},R_*)$ and $q_u=0$ on
$(R_*,\infty)$, we have $\int_0^xq_u\phi\dd y=\frac12\ell(\phi)=0$, so
$E=0$ on $[b^{\mathrm{up}},R_*)$. Hence $E$ is smooth with support in $[a^{\mathrm{up}},b^{\mathrm{up}}]$,
and since $\eta(\xi(x))=E(x)$ and $\xi$ is a smooth decreasing bijection
of $(0,R_*)$ onto $(0,M)$, $\eta\in C_c^\infty((0,M))$. Moreover,
\eqref{eq:e-derivative} reads
\begin{equation}\label{eq:upper-tangent}
 E+xE'=r\rho^\beta\phi\qquad\text{on }(0,R_*).
\end{equation}
Let $\kappa=E'/\rho'$ on $(0,R_*)$; it is smooth with support in
$[a^{\mathrm{up}},b^{\mathrm{up}}]$, because $\rho'<0$ on $(0,R_*)$. After reducing $\sigma_0$
so that $\sigma_0(\|\kappa\|_\infty+\|x\kappa'\|_\infty)<1/2$, the map
$y_\sigma(x)=x/(1+\sigma\kappa(x))$ is, for $|\sigma|<\sigma_0$, a smooth
increasing diffeomorphism of $(0,R_*)$ onto itself, equal to the
identity outside $[a^{\mathrm{up}},b^{\mathrm{up}}]$, with $y_0=\mathrm{id}$ and
$\partial_\sigma y_\sigma|_{\sigma=0}=-x\kappa$.

Let $H_\sigma=H+\sigma\eta$. Differentiating $\eta(\xi(x))=E(x)$ and
using \eqref{eq:mass-coordinate} and \eqref{eq:height-derivatives}, as in
Step~2 of the proof of Lemma~\ref{lem:realization}, we obtain
\begin{equation*}
 H_\sigma'(\xi(x))=\frac{1+\sigma\kappa(x)}{2x}=\frac1{2y_\sigma(x)},\qquad
 H_\sigma''(\xi(x))\,\xi'(x)=-\frac{y_\sigma'(x)}{2y_\sigma(x)^2}<0,
\end{equation*}
so $H_\sigma''>0$ on $(0,M)$, and $H_\sigma=H$ near $\xi=0$ and near
$\xi=M$. Thus $H_\sigma$ is continuous and convex on $[0,M]$ with
$H_\sigma(0)=0$ and $H_\sigma'>0$ on $[0,M)$, and
Lemma~\ref{lem:height-coordinates}(ii) gives the density
$\rho_\sigma=\rho_{H_\sigma}$, with $\int_\R\rho_\sigma\dd x=M$ and height
$H_\sigma$. Since $R_{H_\sigma}(\xi(x))=y_\sigma(x)$,
\eqref{eq:rho-H-formula} gives
\begin{equation*}
 \rho_\sigma(y_\sigma(x))=H_\sigma(\xi(x))=\rho(x)+\sigma E(x)
 \quad(0<x<R_*),\qquad \rho_\sigma=0\ \text{ on }\{|x|\ge R_*\}.
\end{equation*}
In particular $\rho_0=\rho$. Hence $\rho_\sigma$ is smooth in
$(\sigma,x)$ on $(-\sigma_0,\sigma_0)\times(-R_*,R_*)$ and equals $\rho$
outside $\pm[a^{\mathrm{up}},b^{\mathrm{up}}]$. After reducing $\sigma_0$,
$\rho_\sigma\ge\frac12\min_{[a^{\mathrm{up}},b^{\mathrm{up}}]}\rho>0$ on $[a^{\mathrm{up}},b^{\mathrm{up}}]$, so
$u_\sigma^{\mathrm{up}}=\rho_\sigma^{1/r}$ is smooth in $(\sigma,x)$ and
equals $(u-\mu)_+$ outside $\pm[a^{\mathrm{up}},b^{\mathrm{up}}]$. Differentiating at
$\sigma=0$ and using \eqref{eq:upper-tangent},
\begin{equation*}
 \partial_\sigma\rho_\sigma|_{\sigma=0}=E+x\kappa\rho'=E+xE'=r\rho^\beta\phi,
 \qquad
 \partial_\sigma u_\sigma^{\mathrm{up}}|_{\sigma=0}
 =\frac1r\rho^{1/r-1}\cdot r\rho^\beta\phi=\phi
 \quad\text{on }(-R_*,R_*).
\end{equation*}

\emph{Step 3. Gluing.}
Define $u_\sigma=u_\sigma^{\mathrm{lo}}$ on $\{|x|>R_*\}$, $u_\sigma(\pm R_*)=\mu$,
$u_\sigma=\mu+u_\sigma^{\mathrm{up}}$ on $\{|x|<R_*\}$, and
\begin{equation*}
 K_\phi=\bigl\{x\in\R:\ a^{\mathrm{up}}\le|x|\le b^{\mathrm{up}}\bigr\}
 \cup\bigl\{x\in\R:\ \tilde a^{\mathrm{lo}}\le|x|\le\tilde b^{\mathrm{lo}}\bigr\}.
\end{equation*}

(i) By Steps~1 and~2, $u_\sigma=u$ outside $K_\phi$, and $u_0=u$.

(ii) Since $u_\sigma=u$ near $Z$, the two pieces fit together smoothly,
so $(\sigma,x)\mapsto u_\sigma(x)$ is smooth on
$(-\sigma_0,\sigma_0)\times\R$, and Steps~1 and~2 give
$\partial_\sigma u_\sigma|_{\sigma=0}=\phi$.

(iii) Since $u_\sigma-u$ is smooth with compact support,
$u_\sigma\in H^s\cap L^1\cap L^\infty$. Both pieces are nonincreasing in
$|x|$, and $u_\sigma\ge\mu$ on $[-R_*,R_*]$ while $u_\sigma<\mu$
elsewhere, so $u_\sigma$ is symmetric decreasing. Moreover
$(u_\sigma-\mu)_+^r=\rho_\sigma$, so $\mathcal N(u_\sigma)=M$ and
$H_{u_\sigma}=H_\sigma=H+\sigma\eta$. For $0<t<\mu$, since
$R_\sigma(t)>R_*$, we have $\{u_\sigma>t\}=(-R_\sigma(t),R_\sigma(t))$, so
$R_{u_\sigma}=R_\sigma$, that is,
$R_{u_\sigma}^{1-2s}=R_u^{1-2s}+\sigma\Upsilon$. Hence
$u_\sigma\in\mathcal A_M$ and \eqref{eq:affine-data} holds.
\end{proof}

\begin{lemma}\label{lem:small-local}
In case \eqref{eq:caseB}, for compact intervals $X_-\subset(R_*,\infty)$
and $X_+\subset(0,R_*)$ of positive length there is $c_{X_-,X_+}>0$ such
that
\begin{equation}\label{eq:small-local}
 \mathcal Q_u[\phi]\ge c_{X_-,X_+}\left(\int_{X_-}\phi^2\dd x
 +\int_{X_+}(E_\phi')^2\dd x\right)
 \qquad\text{for every }\phi\in\mathcal T.
\end{equation}
\end{lemma}
\begin{proof}
\emph{Step 1. The second-variation identity.}
Let $u_\sigma$ be the curve of Lemma~\ref{lem:small-realization} and
$j(\sigma)=S(u_\sigma)$. Since $u_\sigma-u$ is smooth in $(\sigma,x)$ with
support in $K_\phi$, the map $\sigma\mapsto u_\sigma-u$ is $C^2$ with
values in $C^1$ functions supported in $K_\phi$, hence in $H^s(\R)$, and
after reducing $\sigma_0$ we have $u_\sigma\ge\frac12\min_{K_\phi}u>0$ on
$K_\phi$, so $F(u_\sigma)-F(u)$ is smooth in $(\sigma,x)$ with support in
$K_\phi$. Thus $j\in C^2$, and, as in Step~1 of the proof of
Lemma~\ref{lem:local-lower-bound}, \eqref{eq:weak-all-tests} with
$\psi=\partial_\sigma^2u_\sigma|_{\sigma=0}$ gives
\begin{equation*}
 j''(0)=\mathcal Q_u[\phi]+T_s(u,\psi)+\int_\R(\lambda u-f(u))\psi\dd x=\mathcal Q_u[\phi].
\end{equation*}

\emph{Step 2. A strictly convex mixed term.}
By Lemma~\ref{lem:small-realization}(iii), the functions
$U_\sigma=R_{u_\sigma}^{1-2s}=R_u^{1-2s}+\sigma\Upsilon$ and
$H_\sigma=H_{u_\sigma}=H+\sigma\eta$ are affine in $\sigma$, so for
$\sigma_1,\sigma_2\in(-\sigma_0,\sigma_0)$ the two-level
interpolation between $u_{\sigma_1}$ and $u_{\sigma_2}$ coincides almost
everywhere with $\theta\mapsto u_{(1-\theta)\sigma_1+\theta\sigma_2}$: with
$\sigma=(1-\theta)\sigma_1+\theta\sigma_2$ one has
$(1-\theta)U_{\sigma_1}+\theta U_{\sigma_2}=U_\sigma$ and
$(1-\theta)H_{\sigma_1}+\theta H_{\sigma_2}=H_\sigma$, so the right side
of \eqref{eq:two-level} for the pair $u_{\sigma_1},u_{\sigma_2}$ at
$\theta$ equals the right side of \eqref{eq:two-level-layers} for
$u_\sigma$. Moreover
$R_{u_\sigma}(t)>R_*=R_{H_\sigma}(0)\ge R_{H_\sigma}(\xi)$, so the radii are
ordered as in the proof of Lemma~\ref{lem:two-level-convexity}. Hence
that proof, with $U_\sigma$ and $H_\sigma$ in place of $U_\theta$ and
$H_\theta$, writes
\begin{equation*}
 j(\sigma)=I_1(\sigma)+I_2(\sigma)+I_3(\sigma)
 +2\int_0^\mu(\lambda t-f(t))R_{u_\sigma}(t)\dd t
 -\frac1r\int_0^Mb_{r,\mu}(H_\sigma(\xi))\dd\xi,
\end{equation*}
where every term of the series \eqref{eq:two-level-T1},
\eqref{eq:two-level-T2}, \eqref{eq:kinetic-series} and the two
potential terms are convex in $\sigma$. For compact intervals
$\mathcal J\subset(0,\mu)$ and $\mathcal I\subset(0,M)$ of positive length let
\begin{equation*}
 m(\sigma)=\int_{\mathcal J}\int_{\mathcal I}
 U_\sigma(t)^{-\frac{2s+2}{1-2s}}H_\sigma'(\xi)^{-2}H_\sigma(\xi)^{-\beta}
 \dd\xi\dd t,
\end{equation*}
so that, by \eqref{eq:two-level-T2} with $k=1$,
$c_2m(\sigma)$ with $c_2=\frac{2c_{1,s}}{sr}\cdot\frac{d_{1,s}}{2^3}=\frac{c_{1,s}d_{1,s}}{4sr}$
is the contribution of the term $k=1$ on $\mathcal J\times\mathcal I$ to
$I_2(\sigma)$, and $j-c_2m$ is a sum of convex functions, hence convex.
On $\mathcal J\times\mathcal I$ the functions $U_0=R_u^{1-2s}$, $H'$ and $H$
are bounded and bounded away from zero, so $m$ is $C^2$ and may be
differentiated under the integral sign; since $j$ and $m$ are $C^2$ and
$j-c_2m$ is convex, $j''(0)\ge c_2m''(0)$. By
\eqref{eq:nonpositive-monomial}, applied to the affine functions
$U_\sigma(t)$, $H_\sigma'(\xi)$, $H_\sigma(\xi)$ of $\sigma$ with
derivatives $\Upsilon(t)$, $\eta'(\xi)$, $\eta(\xi)$ and exponents
$(2s+2)/(1-2s)$, $2$, $\beta$,
\begin{equation*}
 m''(0)\ge\int_{\mathcal J}\int_{\mathcal I}k_0(t,\xi)
 \left(\frac{2s+2}{1-2s}\,\frac{\Upsilon(t)^2}{U_0(t)^2}
 +\frac{2\eta'(\xi)^2}{H'(\xi)^2}\right)\dd\xi\dd t,
 \qquad k_0=U_0^{-\frac{2s+2}{1-2s}}(H')^{-2}H^{-\beta}.
\end{equation*}

\emph{Step 3. The bound in terms of $\phi$.}
Let $X_-=R_u(\mathcal J)\subset(R_*,\infty)$ and
$X_+=\xi^{-1}(\mathcal I)\subset(0,R_*)$. With $x=R_u(t)$, that is,
$t=u(x)$, we have $\dd t=|u'(x)|\dd x$ and, by the definition of $\Upsilon$,
\begin{equation*}
 \frac{\Upsilon(t)}{U_0(t)}=\frac{(1-2s)\phi(x)}{x^{2s}|u'(x)|\,x^{1-2s}}
 =\frac{(1-2s)\phi(x)}{x|u'(x)|};
\end{equation*}
with $\xi=\xi(x)$ we have $\dd\xi=|\xi'(x)|\dd x=2x|\rho'(x)|\dd x$ and,
differentiating $\eta(\xi(x))=E_\phi(x)$ and using
\eqref{eq:mass-coordinate} and \eqref{eq:height-derivatives},
\begin{equation*}
 \frac{\eta'(\xi(x))}{H'(\xi(x))}=\frac{E_\phi'(x)/\xi'(x)}{1/(2x)}
 =\frac{E_\phi'(x)}{\rho'(x)}.
\end{equation*}
Bounding $k_0$ from below by its minimum on $\mathcal J\times\mathcal I$ in
Step~2, and integrating the first summand in $\xi$ and the second in $t$,
\begin{equation*}
 m''(0)\ge\min_{\mathcal J\times\mathcal I}k_0\left(\frac{2s+2}{1-2s}|\mathcal I|
 \int_{\mathcal J}\frac{\Upsilon(t)^2}{U_0(t)^2}\dd t
 +2|\mathcal J|\int_{\mathcal I}\frac{\eta'(\xi)^2}{H'(\xi)^2}\dd\xi\right).
\end{equation*}
The two changes of variables above give
\begin{gather*}
 \int_{\mathcal J}\frac{\Upsilon(t)^2}{U_0(t)^2}\dd t
 =\int_{X_-}\frac{(1-2s)^2\phi(x)^2}{x^2|u'(x)|^2}\,|u'(x)|\dd x
 =(1-2s)^2\int_{X_-}\frac{\phi(x)^2}{x^2|u'(x)|}\dd x,\\
 \int_{\mathcal I}\frac{\eta'(\xi)^2}{H'(\xi)^2}\dd\xi
 =\int_{X_+}\frac{E_\phi'(x)^2}{\rho'(x)^2}\,2x|\rho'(x)|\dd x
 =\int_{X_+}\frac{2x\,E_\phi'(x)^2}{|\rho'(x)|}\dd x,
\end{gather*}
and with Step~1 and $\frac{2s+2}{1-2s}(1-2s)^2=2(1+s)(1-2s)$ we obtain
\begin{equation*}
 \mathcal Q_u[\phi]=j''(0)\ge c_2\min_{\mathcal J\times\mathcal I}k_0
 \left(2(1+s)(1-2s)|\mathcal I|\int_{X_-}\frac{\phi(x)^2}{x^2|u'(x)|}\dd x
 +2|\mathcal J|\int_{X_+}\frac{2x\,E_\phi'(x)^2}{|\rho'(x)|}\dd x\right).
\end{equation*}
The weights
$1/(x^2|u'|)$ on $X_-$ and $2x/|\rho'|$ on $X_+$ are continuous and
positive, hence bounded below by positive constants. As $\mathcal J$ and
$\mathcal I$ are arbitrary, so are $X_-$ and $X_+$, and
\eqref{eq:small-local} follows.
\end{proof}

%\subsection*{From smooth tangent directions to coercivity}

\begin{lemma}\label{lem:tangent-positive}
Let $\phi\in D_u$ be even with $\ell(\phi)=0$. Then the bounds
\eqref{eq:kinetic-local-bound} (if $s>\beta$) and
\eqref{eq:potential-local-bound} (if $s=\beta$) in case \eqref{eq:caseA},
and \eqref{eq:small-local} in case \eqref{eq:caseB}, hold for $\phi$, with
$\eta_\phi$ and $E_\phi$ given by \eqref{eq:e-phi}. In particular,
\begin{equation}\label{eq:strict-tangent}
 \mathcal Q_u[\phi]>0\qquad\text{if }\phi\ne0.
\end{equation}
\end{lemma}
\begin{proof}
\emph{Step 1. Approximation by $\mathcal T$.}
By the density statement following \eqref{eq:form-continuity-general},
there is a sequence $\{\psi_j\}_{j\ge1}$ of even functions in $C_c^\infty(\R)$ with $\psi_j\to\phi$ in $D_u$.
In case \eqref{eq:caseB} we also remove neighborhoods of $Z$, using
$s<1/2$. Let $\chi_0\in C_c^\infty((-1,1))$ be even with $0\le\chi_0\le1$
and $\chi_0=1$ on $[-1/2,1/2]$, and for $0<\delta<R_*/2$ define
\begin{equation*}
 \zeta_\delta(x)=\chi_0(x/\delta)+\chi_0((x-R_*)/\delta)+\chi_0((x+R_*)/\delta),
\end{equation*}
which is even, takes values in $[0,1]$, and equals one near $Z$. The
three summands have disjoint supports, $T_s$ is translation invariant
with $T_s(\chi_0(\cdot/\delta))=\delta^{1-2s}T_s(\chi_0)$ by the
substitution $(x,y)\mapsto(\delta x,\delta y)$ in
\eqref{eq:kinetic-double}, and $T_s(a+b+c)\le3(T_s(a)+T_s(b)+T_s(c))$ by
the Cauchy--Schwarz inequality; hence
\begin{equation*}
 T_s(\zeta_\delta)\le9\,\delta^{1-2s}T_s(\chi_0),\qquad
 \|\zeta_\delta\|_2^2=3\delta\|\chi_0\|_2^2.
\end{equation*}
For $\psi\in C_c^\infty(\R)$, the identity
$(\zeta_\delta\psi)(x)-(\zeta_\delta\psi)(y)=\psi(x)(\zeta_\delta(x)-\zeta_\delta(y))+\zeta_\delta(y)(\psi(x)-\psi(y))$
and $(a+b)^2\le2a^2+2b^2$ give
\begin{equation*}
 T_s(\zeta_\delta\psi)\le2\|\psi\|_\infty^2T_s(\zeta_\delta)
 +c_{1,s}\Gamma_\psi\|\zeta_\delta\|_2^2,\qquad
 \Gamma_\psi=\sup_{y\in\R}\int_\R\frac{(\psi(x)-\psi(y))^2}{|x-y|^{1+2s}}\dd x,
\end{equation*}
where $\Gamma_\psi\le\|\psi'\|_\infty^2\int_{|z|\le1}|z|^{1-2s}\dd z
+4\|\psi\|_\infty^2\int_{|z|>1}|z|^{-1-2s}\dd z<\infty$, and
$\int_\R V_+(\zeta_\delta\psi)^2\dd x\le\|V_+\|_{L^\infty(\supp\psi)}\|\psi\|_\infty^2\|\zeta_\delta\|_2^2$.
Since $1-2s>0$, all three bounds tend to zero with $\delta$. Hence
$\zeta_\delta\psi\to0$ in $D_u$ as $\delta\to0^+$, while
$(1-\zeta_\delta)\psi$ is even and belongs to $C_c^\infty(\R\setminus Z)$.
A diagonal choice gives, in both cases, a sequence $\{\widetilde\phi_j\}_{j\ge1}$ of even functions in
$C_c^\infty(\R)$ (in $C_c^\infty(\R\setminus Z)$ in case
\eqref{eq:caseB}) with $\widetilde\phi_j\to\phi$ in $D_u$. Fix an even
$\psi_*\in C_c^\infty(\R)$ with $\ell(\psi_*)=1$, supported in
$(-R_*,R_*)\setminus\{0\}$ in case \eqref{eq:caseB}, which is possible
since $q_u>0$ on an open set, and define the sequence $\{\phi_j\}_{j\ge1}$ by
\begin{equation*}
 \phi_j=\widetilde\phi_j-\ell(\widetilde\phi_j)\psi_*\in\mathcal T.
\end{equation*}
Since $q_u\in L^2$ and $\ell(\phi)=0$,
$|\ell(\widetilde\phi_j)|=|\ell(\widetilde\phi_j-\phi)|\le\|q_u\|_2\|\widetilde\phi_j-\phi\|_2\to0$,
and therefore
\begin{equation*}
 \|\phi_j-\phi\|_{D_u}\le\|\widetilde\phi_j-\phi\|_{D_u}
 +|\ell(\widetilde\phi_j)|\,\|\psi_*\|_{D_u}
 \le\bigl(1+\|q_u\|_2\|\psi_*\|_{D_u}\bigr)\|\widetilde\phi_j-\phi\|_{D_u}\longrightarrow0\quad\text{as }j\to\infty.
\end{equation*}

\emph{Step 2. Passage to the limit.}
By bilinearity and \eqref{eq:form-continuity-general},
\begin{equation*}
 |\mathcal Q_u[\phi_j]-\mathcal Q_u[\phi]|=|\mathcal Q_u[\phi_j-\phi,\phi_j+\phi]|
 \le(1+\|V_-\|_\infty)\|\phi_j-\phi\|_{D_u}\bigl(\|\phi_j\|_{D_u}+\|\phi\|_{D_u}\bigr),
\end{equation*}
which tends to zero as $j\to\infty$.
Since $q_u$ is bounded, Cauchy--Schwarz gives
\begin{equation*}
 |E_{\phi_j}(x)-E_\phi(x)|=\frac rx\Bigl|\int_0^xq_u(\phi_j-\phi)\dd y\Bigr|
 \le\frac{r\|q_u\|_\infty}{\sqrt x}\|\phi_j-\phi\|_2,
\end{equation*}
so $E_{\phi_j}\to E_\phi$ uniformly on compact subsets of $(0,R_*)$ (of
$(0,R_*]$ in case \eqref{eq:caseB}). In case \eqref{eq:caseA}, a compact
interval $\mathcal I\subset(0,M)$ corresponds to a compact interval $[a,b]\subset(0,\infty)$
of $x$, and
\begin{equation*}
 \Bigl|\int_{\mathcal I}\eta_{\phi_j}^2\dd\xi-\int_{\mathcal I}\eta_\phi^2\dd\xi\Bigr|
 \le|\mathcal I|\,\|\eta_{\phi_j}-\eta_\phi\|_{L^\infty(\mathcal I)}
 \bigl(\|\eta_{\phi_j}\|_{L^\infty(\mathcal I)}+\|\eta_\phi\|_{L^\infty(\mathcal I)}\bigr)\longrightarrow0\quad\text{as }j\to\infty.
\end{equation*}
In case \eqref{eq:caseB}, $\int_{X_-}\phi_j^2\to\int_{X_-}\phi^2$, and
\eqref{eq:e-derivative} gives, on $X_+\subset(0,R_*)$,
$E_{\phi_j}'-E_\phi'=x^{-1}\bigl(rq_u(\phi_j-\phi)-(E_{\phi_j}-E_\phi)\bigr)\to0$
in $L^2(X_+)$. Applying the respective bound to $\phi_j$ and letting
$j\to\infty$ proves it for $\phi$; the constants depend only on $u$ and
the fixed intervals.

\emph{Step 3. Strict positivity.}
Let $\phi\ne0$. In case \eqref{eq:caseA}, if $\eta_\phi\equiv0$, then
$\int_0^xu^{r-1}\phi\dd y=0$ for all $x>0$, so $u^{r-1}\phi=0$ almost
everywhere on $(0,\infty)$, and $u>0$ with evenness gives $\phi=0$. Hence
$\eta_\phi$, which is continuous, is nonzero on some compact
interval $\mathcal I\subset(0,M)$ of positive length with $\int_{\mathcal I}\eta_\phi^2\dd\xi>0$; if
$s>\beta$, choose a second compact interval $\mathcal J$ to the right of $\mathcal I$ and
apply \eqref{eq:kinetic-local-bound}, and if $s=\beta$ apply
\eqref{eq:potential-local-bound}. In case \eqref{eq:caseB}, if $\phi$
does not vanish almost everywhere on $(R_*,\infty)$, choose $X_-$ with
$\int_{X_-}\phi^2\dd x>0$. Otherwise $\phi$ does not vanish almost
everywhere on $(0,R_*)$ by evenness; if $E_\phi'=0$ almost everywhere on
$(0,R_*)$, then $E_\phi$ is constant on $(0,R_*]$, equal to
$E_\phi(R_*)=\frac r{2R_*}\ell(\phi)=0$, so that
$\int_0^xq_u\phi\dd y=0$ for all $x$, $q_u\phi=0$ and $\phi=0$ almost
everywhere on $(0,R_*)$, a contradiction; hence there is $X_+$ with
$\int_{X_+}(E_\phi')^2\dd x>0$. In both cases \eqref{eq:small-local}
gives $\mathcal Q_u[\phi]>0$.
\end{proof}

\begin{lemma}\label{lem:stable-form}
Let $0<s<1$, let $\varpi\in L^2(\R)$ be even and nonzero, and let $V$ be an
even locally bounded real potential. Let $V_0=V$ and suppose that the family
$\{V_\theta\}_{0\le\theta<\theta_*}$ of even real potentials satisfies \eqref{eq:secant-structure}
and \eqref{eq:potential-convergence}. Assume that $\mathcal Q_V[\phi]>0$ for every
nonzero even $\phi\in D(V)$ with $\int_\R \varpi\phi\dd x=0$. Then there are
$c,C,\bar\theta>0$ such that
\begin{equation*}
 \mathcal Q_{V_\theta}[\phi]\ge c\|\phi\|_{H^s}^2
                  -C\left|\int_\R \varpi\phi\dd x\right|^2
\end{equation*}
for every $0\le \theta<\bar\theta$ and every even $\phi\in D(V_\theta)$.
\end{lemma}
\begin{proof}
\emph{Step 1. A normalized sequence.}
If the conclusion fails, then by homogeneity there are
sequences $\{\theta_j\}_{j\ge1}$ and $\{\phi_j\}_{j\ge1}$ with
$0\le \theta_j<\min\{\theta_*,1/j\}$ and even $\phi_j\in D(V_{\theta_j})$ such that
\begin{equation*}
 \|\phi_j\|_{H^s}=1,\qquad
 \mathcal Q_{V_{\theta_j}}[\phi_j]+j\Bigl|\int_\R \varpi\phi_j\dd x\Bigr|^2<\frac1j.
\end{equation*}
Since $V_{\theta_j}\ge-C_V$ and $\|\phi_j\|_2\le1$, we have
$\mathcal Q_{V_{\theta_j}}[\phi_j]\ge-C_V$, hence
\begin{equation*}
 \Bigl|\int_\R \varpi\phi_j\dd x\Bigr|^2\le\frac{C_V+1}j,\qquad
 \mathcal Q_{V_{\theta_j}}[\phi_j]<\frac1j.
\end{equation*}
Weak compactness in $H^s$ and local compactness \cite[Theorem 7.1]{DNPV}
give a subsequence converging weakly in $H^s$ to an even $\phi$,
strongly in $L^2([-R,R])$ for every $R$, and almost everywhere; and
$\varpi\in L^2$ gives $\int_\R \varpi\phi\dd x=\lim_j\int_\R \varpi\phi_j\dd x=0$.

\emph{Step 2. The limit is zero.}
By \eqref{eq:secant-structure}, the negative parts $(V_{\theta_j})_-$ vanish
outside $[-R_0,R_0]$ and are bounded by $C_V$; since $a\mapsto a_-$ is
$1$-Lipschitz, \eqref{eq:potential-convergence} and the local $L^2$
convergence give
\begin{equation*}
 \begin{split}
 \Bigl|\int_\R(V_{\theta_j})_-\phi_j^2\dd x-\int_\R V_-\phi^2\dd x\Bigr|
 &\le\|V_{\theta_j}-V\|_{L^\infty([-R_0,R_0])}\|\phi_j\|_2^2\\
 &\quad+C_V\|\phi_j-\phi\|_{L^2([-R_0,R_0])}\bigl(\|\phi_j\|_2+\|\phi\|_2\bigr)
 \longrightarrow0\quad\text{as }j\to\infty.
 \end{split}
\end{equation*}
Hence $T_s(\phi_j)+\int_\R(V_{\theta_j})_+\phi_j^2\dd x=\mathcal Q_{V_{\theta_j}}[\phi_j]+\int_\R(V_{\theta_j})_-\phi_j^2\dd x\le1/j+C_V$,
and Fatou's lemma for the positive potential term (using
\eqref{eq:potential-convergence} and almost everywhere convergence)
together with weak lower semicontinuity of $T_s$ gives
\begin{equation*}
 \begin{split}
 T_s(\phi)+\int_\R V_+\phi^2\dd x
 &\le\liminf_{j\to\infty}
      \Bigl(T_s(\phi_j)+\int_\R(V_{\theta_j})_+\phi_j^2\dd x\Bigr)\\
 &=\liminf_{j\to\infty}\Bigl(\mathcal Q_{V_{\theta_j}}[\phi_j]+\int_\R(V_{\theta_j})_-\phi_j^2\dd x\Bigr)
 \le\int_\R V_-\phi^2\dd x,
 \end{split}
\end{equation*}
which is at most $C_V\|\phi\|_2^2<\infty$. Thus $\phi\in D(V)$, and only
now may the finite negative part be subtracted: $\mathcal Q_V[\phi]\le0$. The
positivity assumption and $\int_\R \varpi\phi\dd x=0$ force $\phi=0$.

\emph{Step 3. Contradiction.}
By \eqref{eq:secant-structure} and $\mathcal Q_{V_{\theta_j}}[\phi_j]<1/j$,
\begin{equation*}
 T_s(\phi_j)+\nu_0\int_{|x|>R_0}\phi_j^2\dd x
 \le T_s(\phi_j)+\int_\R(V_{\theta_j})_+\phi_j^2\dd x
 \le\frac1j+C_V\int_{|x|\le R_0}\phi_j^2\dd x\longrightarrow0\quad\text{as }j\to\infty,
\end{equation*}
since $\phi_j\to\phi=0$ in $L^2([-R_0,R_0])$. Hence
\begin{equation*}
 1=\|\phi_j\|_{H^s}^2\le\max\{1,\nu_0^{-1}\}
 \Bigl(T_s(\phi_j)+\nu_0\int_{|x|>R_0}\phi_j^2\dd x\Bigr)
 +\int_{|x|\le R_0}\phi_j^2\dd x\longrightarrow0\quad\text{as }j\to\infty,
\end{equation*}
a contradiction.
\end{proof}

\begin{proof}[Proof of Theorem~\ref{thm:coercivity}]
The potential $V$ is even and locally bounded, and by
\eqref{eq:linearized-potential-bounds} it satisfies
\eqref{eq:secant-structure} (after enlarging $C_V,R_0$ and reducing
$\nu_0$ if necessary, which does not affect \eqref{eq:secant-structure}
for $V_\theta$, $0<\theta<\theta_*$), so the hypotheses of Lemma~\ref{lem:stable-form}
hold for $0\le \theta<\theta_*$; moreover $q_u$ is even, nonzero and in $L^2$.
Lemma~\ref{lem:tangent-positive} gives strict positivity of $\mathcal Q_u$ on the
even $\phi\in D_u$ with $\ell(\phi)=0$, and Lemma~\ref{lem:stable-form} with $\varpi=q_u$
gives \eqref{eq:penalized-coercivity}; the choice $\theta=0$, $\ell(\phi)=0$
gives constrained coercivity. If $|\ell(\phi)|\le\delta\|\phi\|_{H^s}$
with $\delta=\sqrt{c/(2C)}$, then
\begin{equation*}
 \mathcal Q_{V_\theta}[\phi]\ge c\|\phi\|_{H^s}^2-C\delta^2\|\phi\|_{H^s}^2
 =\frac c2\|\phi\|_{H^s}^2.\qedhere
\end{equation*}
\end{proof}

\begin{remark}\label{rem:general-direction}
The constraint term in \eqref{eq:penalized-coercivity} is necessary:
Lemma~\ref{lem:solution-estimates}(ii) gives $u\in D_u$, and
\eqref{eq:weak-all-tests} with $\psi=u$ yields
\begin{equation*}
 \mathcal Q_u[u]=T_s(u)+\int_\R(\lambda-f'(u))u^2\dd x
 =\int_\R\bigl(uf(u)-u^2f'(u)\bigr)\dd x<0
\end{equation*}
by \eqref{eq:f-assumptions}.
\end{remark}

\section{Proof of Theorem~\ref{thm:main}}\label{sec:proof-main}

We show that a symmetric decreasing solution strictly minimizes $S$
among comparison functions with the same value of $\mathcal N$, and
combine this with the strict maximality of $S$ along the rays
$c\mapsto cu$. By Lemma~\ref{lem:symmetry}, a positive solution is
symmetric decreasing after a translation, and then $u'<0$ on
$(0,\infty)$. For $M>0$, the
comparison class in case \eqref{eq:caseA} is
\begin{equation}\label{eq:comparison-class}
 \mathcal C_M=\left\{v\in H^s(\R)\cap C^\infty(\R):
 v>0,\ v\text{ even},\ v'<0\text{ on }(0,\infty),\
 F(v)\in L^1(\R),\ \int_\R v^r\dd x=M\right\};
\end{equation}
every symmetric decreasing positive solution belongs to $\mathcal C_{\mathcal N(u)}$ by
Lemmas~\ref{lem:solution-estimates} and~\ref{lem:symmetry}, and
Lemma~\ref{lem:reconstruction} provides admissible height functions for the
$r$th powers of its elements. In case \eqref{eq:caseB} the comparison
class is $\mathcal A_M$ from Section~\ref{sec:two-level}. We write
$\mathcal K_M=\mathcal C_M$ in case \eqref{eq:caseA} and
$\mathcal K_M=\mathcal A_M$ in case \eqref{eq:caseB}.

The next lemma expresses the difference $S(v)-S(u)$ by a quadratic form
whose potential satisfies the hypotheses of
Theorem~\ref{thm:coercivity}; it gives strict minimality in case
\eqref{eq:caseA} without differentiating the height interpolation at
the endpoint.

\begin{lemma}\label{lem:remainder}
Assume the hypotheses of Theorem~\ref{thm:main}, and let $u$ be a
positive solution. For a bounded positive function $v$, define
\begin{equation*}
 W_v(x)=\lambda-2\int_0^1(1-\tau)
             f'(u(x)+\tau(v(x)-u(x)))\dd\tau.
\end{equation*}
\begin{enumerate}
\item[(i)] If $v\in H^s(\R)$ and $F(v)\in L^1(\R)$, then
$\int_\R|W_v|\,|v-u|^2\dd x<\infty$ and
\begin{equation}\label{eq:exact-remainder}
 2(S(v)-S(u))=T_s(v-u)+\int_\R W_v(v-u)^2\dd x.
\end{equation}
\item[(ii)] Suppose $u$ is symmetric decreasing and $\{v_\theta\}_{0<\theta<\theta_*}$ is a family
of positive even functions in $L^\infty(\R)$ with
$\|v_\theta-u\|_\infty\to0$ as $\theta\to0^+$. Then $W_{v_\theta}$ is even, and there
exist $C_V,R_0,\nu_0>0$ and $\bar\theta\in(0,\theta_*)$ such that, for $0<\theta<\bar\theta$,
\begin{equation}\label{eq:remainder-potential-bounds}
 W_{v_\theta}\ge-C_V\quad\text{a.e. in }\R,\qquad
 W_{v_\theta}(x)\ge\nu_0\quad\text{for a.e. }x\text{ with }|x|\ge R_0,
\end{equation}
and $\lim_{\theta\to0^+}\|W_{v_\theta}-(\lambda-f'(u))\|_{L^\infty([-R,R])}=0$
for every $R>0$.
Thus \eqref{eq:secant-structure} and \eqref{eq:potential-convergence}
hold for the family $\{V_\theta\}_{0<\theta<\bar\theta}$, $V_\theta=W_{v_\theta}$, with
$V_0=\lambda-f'(u)$.
\end{enumerate}
\end{lemma}
\begin{proof}
\emph{Step 1. Integrability of the Taylor remainder.}
Let $g=v-u$ and choose $K,f_K,F_K$ as in \eqref{eq:aux-definition} on a
common interval $[0,m_*]$ containing the values of $u$ and $v$. By
Lemma~\ref{lem:scalar}, $|f'|=|K-f_K'|\le K+f_K'$ and $tf_K(t)\le2F_K(t)$ on
$(0,m_*]$. For $g(x)\ne0$, the substitution $t=u+\tau g$ along the
segment from $u(x)$ to $v(x)$ gives
\begin{equation*}
 g^2\int_0^1f_K'(u+\tau g)\dd\tau=g\int_u^vf_K'(t)\dd t=g\,(f_K(v)-f_K(u))
 =|g|\,|f_K(v)-f_K(u)|,
\end{equation*}
the last equality because $f_K$ is increasing, so that $g$ and
$f_K(v)-f_K(u)$ have the same sign; both sides vanish when $g=0$. Since
$f_K\ge0$ is increasing, $|g|\le\max\{u,v\}$ and
$|f_K(v)-f_K(u)|\le\max\{f_K(u),f_K(v)\}$, we have
$|g|\,|f_K(v)-f_K(u)|\le\max\{u,v\}\max\{f_K(u),f_K(v)\}\le uf_K(u)+vf_K(v)$, and
therefore
\begin{equation}\label{eq:Taylor-integrability}
 \begin{split}
 g^2\int_0^1(1-\tau)|f'(u+\tau g)|\dd\tau
 &\le\frac K2g^2+g^2\int_0^1f_K'(u+\tau g)\dd\tau\\
 &\le\frac K2g^2+uf_K(u)+vf_K(v)
 \le\frac K2g^2+2F_K(u)+2F_K(v)\in L^1(\R),
 \end{split}
\end{equation}
because $F_K(v)=\frac K2v^2-F(v)\in L^1$ by assumption and $F_K(u)\in L^1$ by
Lemma~\ref{lem:solution-estimates}(ii).

\emph{Step 2. An exact formula for $S(v)-S(u)$.}
Taylor's formula with integral remainder on the segment from $u(x)$ to
$v(x)$, which is justified by \eqref{eq:Taylor-integrability}, gives
pointwise
\begin{equation}\label{eq:pointwise-Taylor}
 F(v)-F(u)-f(u)g=g^2\int_0^1(1-\tau)f'(u+\tau g)\dd\tau.
\end{equation}
Expanding the quadratic terms of $S(v)-S(u)$ with $v=u+g$,
\begin{equation*}
 S(v)-S(u)=T_s(u,g)+\frac12T_s(g)+\lambda\int_\R ug\dd x+\frac\lambda2\|g\|_2^2
 -\int_\R\bigl(F(v)-F(u)\bigr)\dd x,
\end{equation*}
and inserting \eqref{eq:pointwise-Taylor}, the terms linear in $g$
combine to $T_s(u,g)+\int_\R(\lambda u-f(u))g\dd x$, which vanishes by
\eqref{eq:weak-all-tests} with $\psi=g\in H^s$. Hence
\begin{equation*}
 S(v)-S(u)=\frac12T_s(g)+\frac\lambda2\|g\|_2^2
 -\int_\R g^2\int_0^1(1-\tau)f'(u+\tau g)\dd\tau\dd x
 =\frac12\Bigl(T_s(g)+\int_\R W_vg^2\dd x\Bigr),
\end{equation*}
which is \eqref{eq:exact-remainder}, and by the definition of $W_v$ and
\eqref{eq:Taylor-integrability},
\begin{equation*}
 \int_\R|W_v|g^2\dd x\le\lambda\|g\|_2^2
 +2\int_\R g^2\int_0^1(1-\tau)|f'(u+\tau g)|\dd\tau\dd x<\infty.
\end{equation*}

\emph{Step 3. Uniform bounds and local convergence of the potential.}
For small $\theta$, $u$ and $v_\theta$ take values in a common interval $(0,m_*]$,
and so do the intermediate values $u+\tau(v_\theta-u)$; hence
\begin{equation*}
 W_{v_\theta}\ge\lambda-2\int_0^1(1-\tau)\sup_{0<t\le m_*}f'(t)\dd\tau
 =\lambda-\sup_{0<t\le m_*}f'(t)=:-C_V>-\infty;
\end{equation*}
here $C_V>0$, since by Lemma~\ref{lem:solution-estimates}(v) there is
$t\in(0,\|u\|_\infty]$ with $f(t)/t\ge\lambda$, and $tf'(t)>f(t)$ gives
$f'(t)>\lambda$.
Choose $\varepsilon,\nu_0>0$ with $f'(t)\le\lambda-\nu_0$ for
$0<t\le\varepsilon$, then $R_0$ with $u\le\varepsilon/2$ on
$\{|x|\ge R_0\}$ ($u\in C_0(\R)$), and then $\bar\theta$ with
$\|v_\theta-u\|_\infty\le\varepsilon/2$ for $\theta<\bar\theta$; then $u,v_\theta\le\varepsilon$
almost everywhere on $\{|x|\ge R_0\}$, every intermediate value
obeys the same bound, and
\begin{equation*}
 W_{v_\theta}(x)\ge\lambda-2\int_0^1(1-\tau)(\lambda-\nu_0)\dd\tau=\nu_0
 \qquad\text{for a.e. }x\text{ with }|x|\ge R_0.
\end{equation*}
This proves \eqref{eq:remainder-potential-bounds}. On a compact
interval $I$, positivity gives $\min_Iu>0$, and for small $\theta$ all
intermediate values lie almost everywhere in $[\frac12\min_Iu,m_*]$,
where $f'$ is uniformly continuous; hence
\begin{equation*}
 \begin{split}
 \|W_{v_\theta}-(\lambda-f'(u))\|_{L^\infty(I)}
 &=\operatorname*{ess\,sup}_{x\in I}\Bigl|2\int_0^1(1-\tau)\bigl(f'(u+\tau(v_\theta-u))-f'(u)\bigr)\dd\tau\Bigr|\\
 &\le\operatorname*{ess\,sup}_{x\in I}\sup_{0\le\tau\le1}
 \bigl|f'(u(x)+\tau(v_\theta(x)-u(x)))-f'(u(x))\bigr|,
 \end{split}
\end{equation*}
which tends to zero as $\theta\to0^+$.
Taking $I=[-R,R]$ gives
$\|W_{v_\theta}-(\lambda-f'(u))\|_{L^\infty([-R,R])}\to0$ as $\theta\to0^+$.
Finally, $W_{v_\theta}(x)$ depends on $x$ only through $u(x)$ and
$v_\theta(x)$, so $W_{v_\theta}$ is even because $u$ and $v_\theta$ are even.
\end{proof}

In case \eqref{eq:caseB} we use instead the following one-sided bound,
which requires only an upper bound for $f'$.

\begin{lemma}\label{lem:one-sided}
Assume the hypotheses of Theorem~\ref{thm:main}, let $u$ be a positive
solution, and let $w\in H^s(\R)\cap L^\infty(\R)$ with $w\ge0$ and
$F(w)\in L^1(\R)$. Let $m_*=\max\{\|u\|_\infty,\|w\|_\infty\}$ and
$C_*=\max\{0,\sup_{0<t\le m_*}f'(t)\}<\infty$. Then
\begin{equation*}
 S(w)-S(u)\ge-\frac{C_*}2\,\|w-u\|_{L^2(\R)}^2.
\end{equation*}
\end{lemma}
\begin{proof}
The constant $C_*$ is finite by \eqref{eq:f-regularity} and
\eqref{eq:f-assumptions}, and $F(u)\in L^1(\R)$ by
Lemma~\ref{lem:solution-estimates}(ii). Expanding the quadratic terms
and applying \eqref{eq:weak-all-tests} with $\psi=w-u$ as in Step~2 of
the proof of Lemma~\ref{lem:remainder},
\begin{equation*}
 S(w)-S(u)=\frac12T_s(w-u)+\frac\lambda2\|w-u\|_2^2
 -\int_\R\bigl[F(w)-F(u)-f(u)(w-u)\bigr]\dd x,
\end{equation*}
where $F(w)-F(u)-f(u)(w-u)\in L^1(\R)$ because $F(w),F(u)\in L^1$,
$f(u)\in L^2$ and $w-u\in L^2$. For $t_1,t_2\in[0,m_*]$, the mean value theorem on $(0,m_*]$ and the
continuity of $f$ at zero give $f(\tau)-f(t_1)\le C_*(\tau-t_1)$ for
$\tau\ge t_1$ and $f(t_1)-f(\tau)\le C_*(t_1-\tau)$ for $\tau\le t_1$;
hence, in both orders of $t_1,t_2$,
\begin{equation*}
 F(t_2)-F(t_1)-f(t_1)(t_2-t_1)=\int_{t_1}^{t_2}\bigl(f(\tau)-f(t_1)\bigr)\dd\tau
 \le\frac{C_*}2(t_2-t_1)^2.
\end{equation*}
Applying this with $t_1=u(x)$, $t_2=w(x)$ and dropping the two
nonnegative terms proves the claim.
\end{proof}

\begin{lemma}\label{lem:own-mass}
Assume the hypotheses of Theorem~\ref{thm:main}. Let $u$ be a
symmetric decreasing positive solution and let $M=\mathcal N(u)$. Then
\begin{equation}\label{eq:own-mass}
 S(v)>S(u)\qquad\text{for every }v\in\mathcal K_M\text{ with }v\ne u.
\end{equation}
\end{lemma}
\begin{proof}
\emph{Case \eqref{eq:caseA}, Step 1. The comparison curve and the
linearized constraint.}
Fix $v\in\mathcal C_M$ with $v\ne u$, and let $u_\theta$ be the height
interpolation from $u$ to $v$, with $H_0,H_1$ the height functions of $u^r,v^r$.
By Lemmas~\ref{lem:height-interpolation} and~\ref{lem:path-finiteness},
$\mathcal N(u_\theta)=M$, $u_\theta\le m_*:=\max\{\|u\|_\infty,\|v\|_\infty\}$,
$u_\theta\in H^s$, $F(u_\theta)\in L^1$, and $g_\theta=u_\theta-u$ satisfies
$\|g_\theta\|_{H^s}+\|g_\theta\|_\infty\to0$ as $\theta\to0^+$. Taylor's
formula for $t\mapsto t^r$ on $[0,m_*]$, integrated over $\R$, and
$\mathcal N(u_\theta)=\mathcal N(u)$ give
\begin{equation*}
 0=\int_\R\bigl((u+g_\theta)^r-u^r\bigr)\dd x
 =r\int_\R u^{r-1}g_\theta\dd x
 +r(r-1)\int_\R g_\theta^2\int_0^1(1-\tau)(u+\tau g_\theta)^{r-2}\dd\tau\dd x,
\end{equation*}
and since $r\ge2$ the inner integral is at most $\frac12m_*^{r-2}$; with
$\ell(g_\theta)=\int_\R u^{r-1}g_\theta\dd x$ this gives
\begin{equation}\label{eq:mass-defect-bound}
 |\ell(g_\theta)|
 \le\frac{r-1}{2}m_*^{r-2}\|g_\theta\|_2^2
 \le\frac{r-1}{2}m_*^{r-2}\|g_\theta\|_{H^s}^2.
\end{equation}
If $g_\theta=0$ for some $\theta\in(0,1]$, then
$\rho_{H_\theta}=u_\theta^r=u^r=\rho_{H_0}$ almost everywhere, so
$H_\theta=H_0$, because $H$ is the height function of $\rho_H$
(Lemma~\ref{lem:height-coordinates}(ii)) and the height function depends only on
the almost everywhere class of the density; since $\theta>0$, this gives
$H_1=H_0$, hence $v^r=\rho_{H_1}=\rho_{H_0}=u^r$ almost everywhere by
Lemma~\ref{lem:height-coordinates}(iii), that is, $v=u$, a contradiction.
So $g_\theta\ne0$ for $0<\theta\le1$, and dividing
\eqref{eq:mass-defect-bound} by $\|g_\theta\|_{H^s}$,
\begin{equation*}
 \frac{|\ell(g_\theta)|}{\|g_\theta\|_{H^s}}
 \le\frac{r-1}2m_*^{r-2}\|g_\theta\|_{H^s}\longrightarrow0\quad\text{as }\theta\to0^+.
\end{equation*}

\emph{Step 2. Strict increase of $S$ near the endpoint.}
By Lemma~\ref{lem:remainder}(ii), applied to the family $\{v_\theta\}_{0<\theta<1}$ with $v_\theta=u_\theta$
(positive, even, bounded, $\|u_\theta-u\|_\infty\to0$), there
is $\bar\theta\in(0,1)$ such that the family $\{V_\theta\}_{0<\theta<\bar\theta}$
with $V_\theta=W_{u_\theta}$ satisfies \eqref{eq:secant-structure} and
\eqref{eq:potential-convergence}. Hence Theorem~\ref{thm:coercivity}
applies to this family, with $\theta_*=\bar\theta$, and yields
$c,C,\delta>0$ such that \eqref{eq:almost-tangent-coercivity} holds for
all sufficiently small $\theta>0$; and
Lemma~\ref{lem:remainder}(i) gives
$\int_\R(W_{u_\theta})_+g_\theta^2\dd x\le\int_\R|W_{u_\theta}|g_\theta^2\dd x<\infty$,
so $g_\theta\in D(V_\theta)$, and $g_\theta$ is even. By Step~1, $|\ell(g_\theta)|\le\delta\|g_\theta\|_{H^s}$
for small $\theta$, and \eqref{eq:almost-tangent-coercivity} together
with \eqref{eq:exact-remainder} gives, for all small $\theta>0$,
\begin{equation*}
 2\bigl(S(u_\theta)-S(u)\bigr)=T_s(g_\theta)+\int_\R W_{u_\theta}g_\theta^2\dd x
 =\mathcal Q_{W_{u_\theta}}[g_\theta]\ge\frac c2\|g_\theta\|_{H^s}^2=:2c_*\|g_\theta\|_{H^s}^2.
\end{equation*}

\emph{Step 3. Strict minimality along the whole curve.}
By \eqref{eq:action-convexity} with $\theta_0=0$, $\theta_1=1$ and
$t=\theta$, $S(u_\theta)\le(1-\theta)S(u)+\theta S(v)$. Combined with
Step~2,
\begin{equation*}
 0<c_*\|g_\theta\|_{H^s}^2\le S(u_\theta)-S(u)\le\theta\bigl(S(v)-S(u)\bigr)
\end{equation*}
for all small $\theta>0$, and dividing by $\theta>0$ proves
\eqref{eq:own-mass} in case \eqref{eq:caseA}.

\emph{Case \eqref{eq:caseB}.}
By Lemmas~\ref{lem:solution-estimates} and~\ref{lem:symmetry}, $u$ is
continuous, bounded, integrable and symmetric decreasing with $\|u\|_\infty>\mu$, so
$M>0$ and $u\in\mathcal A_M$. Let $v\in\mathcal A_M$, let $u_\theta$ be the
two-level interpolation with $u_0=u$ and $u_1=v$, and define
$j(\theta)=S(u_\theta)$. By Lemma~\ref{lem:two-level-convexity}(iii), $j$ is
finite and convex on $[0,1]$. By Lemma~\ref{lem:two-level}(iii),
$\|u_\theta\|_\infty\le\max\{\|u\|_\infty,\|v\|_\infty\}=:m_*$ and
$u_\theta\in L^1$, so $F(u_\theta)\in L^1$ (as $F$ is bounded on $[0,m_*]$
and $|F(t)|\le Ct$ there), and Lemma~\ref{lem:one-sided} applies with
the constant $C_*=\max\{0,\sup_{0<t\le m_*}f'(t)\}$, which does not depend on
$\theta$. Together with Lemma~\ref{lem:two-level-endpoint},
\begin{equation*}
 j(\theta)-j(0)=S(u_\theta)-S(u)\ge-\frac{C_*}2\|u_\theta-u\|_2^2=o(\theta)
 \quad\text{as }\theta\to0^+.
\end{equation*}
Suppose, to the contrary, that $S(v)\le S(u)$ for some
$v\in\mathcal A_M$ with $v\ne u$. Then $j(1)\le j(0)$, and convexity gives
$j(\theta)\le(1-\theta)j(0)+\theta j(1)\le j(0)$ for $0\le\theta\le1$. By
Lemma~\ref{lem:two-level-convexity}(iii), $j$ is not affine, so
$j\not\equiv j(0)$ and there is $\theta_1\in(0,1)$ with $j(\theta_1)<j(0)$.
Convexity on $[0,\theta_1]$ then gives, for $0\le\theta\le\theta_1$,
\begin{equation*}
 j(\theta)-j(0)\le\frac\theta{\theta_1}\bigl(j(\theta_1)-j(0)\bigr),
\end{equation*}
with a negative coefficient of $\theta$, contradicting
$j(\theta)-j(0)\ge o(\theta)$. Therefore $S(v)>S(u)$.
\end{proof}

\begin{lemma}\label{lem:amplitude}
Assume the hypotheses of Theorem~\ref{thm:main}, and let $u$ be a
positive solution. Then $S(cu)$ is finite for every $c>0$, and
$S(cu)<S(u)$ for $c>0$ with $c\ne1$.
\end{lemma}
\begin{proof}
Given $c>0$, fix $0<c_0<\min\{1,c\}$ and $c_1>\max\{1,c\}$, and choose
$K,f_K,F_K$ as in \eqref{eq:aux-definition}
on $[0,\max\{1,c_1\}\|u\|_\infty]$. By \eqref{eq:scalar-scaling},
\begin{equation*}
 |F(cu)|\le\frac K2c^2u^2+F_K(cu)\le\frac K2c^2u^2+c\max\{1,c\}\,uf_K(u)\in L^1(\R),
\end{equation*}
where $uf_K(u)\in L^1$ by Lemma~\ref{lem:solution-estimates}(ii); so
$S(cu)$ is finite. Uniformly for $c\in[c_0,c_1]$,
$|uf(cu)|=|u(Kcu-f_K(cu))|\le Kc_1u^2+\max\{1,c_1\}\,uf_K(u)\in L^1(\R)$,
so $c\mapsto\int_\R F(cu)\dd x$ may be differentiated under the
integral, and with \eqref{eq:weak-all-tests} for $\psi=u$, that is,
$T_s(u)+\lambda\|u\|_2^2=\int_\R uf(u)\dd x$,
\begin{equation}\label{eq:amplitude-derivative}
 \frac d{dc}S(cu)
 =c\bigl(T_s(u)+\lambda\|u\|_2^2\bigr)-\int_\R uf(cu)\dd x
 =c\int_\R u^2\left(\frac{f(u)}u-\frac{f(cu)}{cu}\right)\dd x.
\end{equation}
Since $\frac d{dt}(f(t)/t)=(tf'(t)-f(t))/t^2>0$ by
\eqref{eq:f-assumptions}, the integrand is positive for $0<c<1$ and
negative for $c>1$. Hence
\begin{gather*}
 S(u)-S(cu)=\int_c^1\frac d{d\tau}S(\tau u)\dd\tau>0\qquad(0<c<1),\\
 S(cu)-S(u)=\int_1^c\frac d{d\tau}S(\tau u)\dd\tau<0\qquad(c>1).\qedhere
\end{gather*}
\end{proof}

\begin{proof}[Proof of Theorem~\ref{thm:main}]
By Lemma~\ref{lem:symmetry}, translate two positive solutions $u_0,u_1$
so that both are symmetric decreasing, and let $M_i=\mathcal N(u_i)$. In case
\eqref{eq:caseA}, boundedness, positivity and $r\ge2$ give
$0<M_i=\int_\R u_i^r\dd x\le\|u_i\|_\infty^{r-2}\|u_i\|_2^2<\infty$; in case
\eqref{eq:caseB}, Lemma~\ref{lem:solution-estimates}(iii),(v) gives
$0<M_i\le\|u_i\|_\infty^{r-1}\|u_i\|_1<\infty$. Suppose $M_0<M_1$. We
choose $c>1$ and $c'\in(0,1)$ with
\begin{equation*}
 \mathcal N(cu_0)=M_1,\qquad \mathcal N(c'u_1)=M_0.
\end{equation*}
In case \eqref{eq:caseA}, $\mathcal N(cu_i)=c^rM_i$, so $c=(M_1/M_0)^{1/r}$ and
$c'=c^{-1}$. In case \eqref{eq:caseB}, $c\mapsto \mathcal N(cu_i)$ is continuous on
$(0,\infty)$ by dominated convergence, since
$(cu_i-\mu)_+^r\le c^r\|u_i\|_\infty^{r-1}u_i\in L^1$ locally uniformly in
$c$; it vanishes for $c\le\mu/u_i(0)$, since then $cu_i\le\mu$; it is
nondecreasing in $c$; and it tends to infinity as $c\to\infty$, because
$\mathcal N(cu_i)\ge(cu_i(0)/2-\mu)_+^r\,|\{u_i>u_i(0)/2\}|$. Hence, by the intermediate value theorem on $[1,\infty)$ for $u_0$ and
on $(0,1]$ for $u_1$, there are $c>1$ with $\mathcal N(cu_0)=M_1$ and $c'\in(0,1)$
with $\mathcal N(c'u_1)=M_0$; $c\ne1\ne c'$ because $M_0\ne M_1$. Multiplication by a positive constant preserves positivity,
evenness, smoothness and $u'<0$ on $(0,\infty)$, and
Lemma~\ref{lem:amplitude} gives $F(cu_0),F(c'u_1)\in L^1$; hence
$cu_0\in\mathcal K_{M_1}$ and $c'u_1\in\mathcal K_{M_0}$. Applying
Lemma~\ref{lem:own-mass} to $u_1$ with $v=cu_0$, Lemma~\ref{lem:amplitude}
to $u_0$ with $c>1$, Lemma~\ref{lem:own-mass} to $u_0$ with $v=c'u_1$, and
Lemma~\ref{lem:amplitude} to $u_1$ with $c'<1$, we obtain
\begin{equation*}
 S(u_1)\le S(cu_0)<S(u_0)\le S(c'u_1)<S(u_1),
\end{equation*}
a contradiction. Exchanging $u_0$ and $u_1$ excludes $M_0>M_1$, so
$M_0=M_1$. If the two translated solutions were distinct,
Lemma~\ref{lem:own-mass} applied in both orders would give
$S(u_0)<S(u_1)<S(u_0)$. Hence $u_0=u_1$ after these translations, which proves
Theorem~\ref{thm:main}.
\end{proof}

\section{Nondegeneracy and uniqueness at prescribed \texorpdfstring{$L^2$}{L2} norm}\label{sec:nondegeneracy}

The kernel of $L_u$ splits into even and odd functions. Even elements of
the kernel are excluded by Lemma~\ref{lem:tangent-positive}. For odd
elements we write the quadratic form on the half-line $(0,\infty)$,
where its integral kernel is positive, and use that $-u'>0$ there; this
Perron--Frobenius type argument shows that odd elements of the kernel
are multiples of $u'$. It is the one-dimensional analogue of
\cite {FLS}; see also \cite{FW}. Throughout, $u$ is a symmetric
decreasing positive solution under the hypotheses of
Theorem~\ref{thm:main}, with $V,D_u,\mathcal Q_u,\ell,q_u$ as in
Section~\ref{sec:coercivity}.

\smallskip
\noindent\textit{The operator.}
The space $D_u$ is complete in the norm $\|\cdot\|_{D_u}$: if $\{\phi_j\}_{j\ge1}$ is a Cauchy sequence in $D_u$, then $\phi_j\to\phi$ in
$H^s$ and $\sqrt{V_+}\phi_j\to w$ in $L^2$ as $j\to\infty$, for some
$\phi\in H^s(\R)$ and $w\in L^2(\R)$. Along a subsequence, both
convergences hold almost everywhere, which gives $w=\sqrt{V_+}\phi$.
Hence $\phi\in D_u$ and $\phi_j\to\phi$ in $D_u$ as $j\to\infty$. By
\eqref{eq:linearized-potential-bounds},
\begin{equation*}
 \|\phi\|_{D_u}^2=T_s(\phi)+\|\phi\|_2^2+\int_\R V_+\phi^2\dd x
 \le \mathcal Q_u[\phi]+(C_V+1)\|\phi\|_2^2\le(C_V+1)\|\phi\|_{D_u}^2,
\end{equation*}
so $\mathcal Q_u$ is a densely defined closed form bounded below, and the
representation theorem \cite[Chapter VI, Section 2]{Kato} defines a
self-adjoint operator $L_u$ by
\begin{equation}\label{eq:operator-domain}
 D(L_u)=\left\{\phi\in D_u:\ \text{there is }w\in L^2(\R)
 \text{ with }\mathcal Q_u[\phi,\psi]=\int_\R w\psi\dd x
 \text{ for every }\psi\in D_u\right\},
\end{equation}
with $L_u\phi=w$, which is unique by density. Testing against
$C_c^\infty(\R)$ identifies $L_u\phi=(-\Delta)^s\phi+V\phi$ in
distributions. Conversely, let $\phi\in D_u$ satisfy
$(-\Delta)^s\phi+V\phi=0$ in the sense of distributions, that is,
\begin{equation*}
 \int_\R\phi\,(-\Delta)^s\psi\dd x+\int_\R V\phi\psi\dd x=0
 \qquad\text{for every }\psi\in C_c^\infty(\R).
\end{equation*}
Since $\phi\in H^s(\R)$, Plancherel's identity gives
$\int_\R\phi\,(-\Delta)^s\psi\dd x
=\int_\R|\omega|^{2s}\widehat\phi\,\overline{\widehat\psi}\dd\omega
=T_s(\phi,\psi)$, and the second integral is finite because $V$ is
locally bounded and $\psi$ has compact support. Hence
$\mathcal Q_u[\phi,\psi]=0$ for every $\psi\in C_c^\infty(\R)$. Now let
$\psi\in D_u$. By the density statement following
\eqref{eq:form-continuity-general}, there is a sequence $\{\psi_k\}_{k\ge1}\subset C_c^\infty(\R)$
with $\psi_k\to\psi$ in $D_u$, and by
\eqref{eq:form-continuity-general},
\begin{equation*}
 |\mathcal Q_u[\phi,\psi]|=|\mathcal Q_u[\phi,\psi-\psi_k]|
 \le(1+\|V_-\|_\infty)\|\phi\|_{D_u}\|\psi-\psi_k\|_{D_u}
 \longrightarrow0\quad\text{as }k\to\infty.
\end{equation*}
Thus $\mathcal Q_u[\phi,\psi]=0=\int_\R0\cdot\psi\dd x$ for every
$\psi\in D_u$, and \eqref{eq:operator-domain} with $w=0$ gives
$\phi\in D(L_u)$ and $L_u\phi=0$, that is, $\phi\in\ker L_u$.

\begin{lemma}\label{lem:translation-domain}
Under the hypotheses of Theorem~\ref{thm:nondegeneracy},
$u'\in D(L_u)$ and $L_uu'=0$.
\end{lemma}
\begin{proof}
\emph{Step 1. Difference quotients and a uniform form bound.}
For $0<h<1$ let
\begin{equation*}
 d_h(x)=\frac{u(x+h)-u(x)}h,\qquad
 V_h(x)=\lambda-\int_0^1f'\bigl((1-t)u(x)+tu(x+h)\bigr)\dd t.
\end{equation*}
Subtracting \eqref{eq:weak-all-tests} for $u$ from the same identity
for the translate $u(\cdot+h)$, which is also a solution, and dividing
by $h$, we obtain
\begin{equation}\label{eq:translation-quotient-equation}
 T_s(d_h,\psi)+\int_\R V_hd_h\psi\dd x=0
 \qquad\text{for every }\psi\in H^s(\R),
\end{equation}
where $V_hd_h=\lambda d_h-(f(u(\cdot+h))-f(u))/h$ by the mean value
formula. Here $d_h\in H^s$ and $V_hd_h\in L^2$, so $|V_h|d_h^2\in L^1$ and
the choice $\psi=d_h$ is admissible. As in
\eqref{eq:linearized-potential-bounds}, $V_h\ge-C_1$ with
$C_1=\max\{1,\sup_{0<t\le\|u\|_\infty}f'(t)-\lambda\}$ independent of $h$; and if
$t_1,\nu_1>0$ are chosen with $f'\le\lambda-\nu_1$ on $(0,t_1]$ and $R_1$
with $u\le t_1$ outside $[-R_1,R_1]$, then for $|x|>R_1+1$ we have $u(x)\le t_1$ and
$u(x+h)\le t_1$, so $V_h(x)\ge\nu_1$. Thus $(V_h)_-$ vanishes on $\{|x|>R_1+1\}$ and is
bounded by $C_1$ on $\{|x|\le R_1+1\}$, and
$|d_h(x)|=|\int_0^1u'(x+th)\dd t|\le C_2:=\max_{|y|\le R_1+2}|u'|$ for $|x|\le R_1+1$.
Testing \eqref{eq:translation-quotient-equation} with $\psi=d_h$,
\begin{equation}\label{eq:translation-form-bound}
 T_s(d_h)+\int_\R(V_h)_+d_h^2\dd x
 =\int_\R(V_h)_-d_h^2\dd x\le 2(R_1+1)C_1C_2^2=:C_3,
\end{equation}
and since $(V_h)_+\ge\nu_1$ on $\{|x|>R_1+1\}$, also
$\nu_1\int_{|x|>R_1+1}d_h^2\dd x\le C_3$; with
$\int_{|x|\le R_1+1}d_h^2\dd x\le 2(R_1+1)C_2^2$ this bounds $\|d_h\|_2$, hence
$\|d_h\|_{H^s}$, uniformly in $h$.

\emph{Step 2. Passage to the form domain.}
Since $u$ is smooth and positive on compact intervals, $d_h\to u'$ and
$V_h\to V$ locally uniformly as $h\to0^+$; by the uniform $H^s$ bound,
$d_h\rightharpoonup u'$ weakly in $H^s$ as $h\to0^+$ (every sequence
$\{h_j\}_{j\ge1}$ with $h_j\to0^+$ has a subsequence along which
$d_{h_j}$ converges weakly in $H^s$, and the limit is $u'$), so $u'\in H^s$,
and weak lower semicontinuity of $T_s$ together with Fatou's lemma for
the positive potential term (using the local uniform convergence and
pointwise convergence) give
\begin{equation*}
 T_s(u')+\int_\R V_+(u')^2\dd x
 \le\liminf_{h\to0^+}\Bigl(T_s(d_h)+\int_\R(V_h)_+d_h^2\dd x\Bigr)\le C_3<\infty.
\end{equation*}
Consequently $u'\in D_u$.

\emph{Step 3. The equation on the full form domain.}
Since $u$ and $f(u)$ are smooth, \eqref{eq:main} may be differentiated
in the sense of distributions, giving $(-\Delta)^su'+Vu'=0$. As
$u'\in H^s$ and $V$ is locally bounded, this means
$\mathcal Q_u[u',\psi]=T_s(u',\psi)+\int_\R Vu'\psi\dd x=0$ for every
$\psi\in C_c^\infty(\R)$, and by density in $D_u$ and
\eqref{eq:form-continuity-general} for every $\psi\in D_u$. By
\eqref{eq:operator-domain}, $u'\in D(L_u)$ with $L_uu'=0$.
\end{proof}

\begin{lemma}\label{lem:odd-kernel}
Let $0<s<1$, and let $V$ be even and locally bounded, with $V\ge-C_V$
on $\R$ and $V\ge\nu_0$ on $\R\setminus[-R_0,R_0]$ for some constants
$C_V,R_0,\nu_0>0$. Suppose that an odd
$z\in D(V)$ satisfies $z>0$ almost everywhere on $(0,\infty)$ and
$\mathcal Q_V[z,\psi]=0$ for every $\psi\in D(V)$. Then every odd $\phi\in D(V)$
satisfies
\begin{equation}\label{eq:odd-nonnegative}
 \mathcal Q_V[\phi]\ge0,
 \qquad \mathcal Q_V[\phi]=0\ \text{if and only if }\phi\in\operatorname{span}\{z\}.
\end{equation}
\end{lemma}
\begin{proof}
\emph{Step 1. Comparison of quadratic forms on the positive half-line.}
For $x,y>0$ with $x\ne y$, define the positive kernel
\begin{equation*}
 K_s(x,y)=\frac1{|x-y|^{1+2s}}-\frac1{(x+y)^{1+2s}}>0.
\end{equation*}
For odd $\phi$, we split the double integral defining $T_s(\phi)$ into
the four quadrants. The substitution $(x,y)\mapsto(-x,-y)$ maps the
quadrant $\{x<0,y<0\}$ onto $\{x>0,y>0\}$ and leaves
$(\phi(x)-\phi(y))^2|x-y|^{-1-2s}$ invariant, while $y\mapsto-y$ maps
$\{x>0,y<0\}$ onto $\{x>0,y>0\}$ and turns this integrand into
$(\phi(x)+\phi(y))^2(x+y)^{-1-2s}$, because $\phi(-y)=-\phi(y)$; the same
holds for $\{x<0,y>0\}$ with $x\mapsto-x$. Hence the four quadrants
contribute twice each of the two integrals, and noting the factor
$c_{1,s}/2$ in \eqref{eq:kinetic-double},
\begin{equation*}
 T_s(\phi)=c_{1,s}\int_0^\infty\int_0^\infty
 \left(\frac{(\phi(x)-\phi(y))^2}{|x-y|^{1+2s}}
       +\frac{(\phi(x)+\phi(y))^2}{(x+y)^{1+2s}}\right)\dd y\dd x.
\end{equation*}
We write $(a+b)^2=2(a^2+b^2)-(a-b)^2$ in the second term. The part with
$2(a^2+b^2)$ is, by symmetry in $(x,y)$ and
$\int_0^\infty(x+y)^{-1-2s}\dd y=x^{-2s}/(2s)$,
\begin{equation*}
 \begin{split}
 c_{1,s}\int_0^\infty\!\!\int_0^\infty\frac{2(\phi(x)^2+\phi(y)^2)}{(x+y)^{1+2s}}\dd y\dd x
 &=4c_{1,s}\int_0^\infty\phi(x)^2\int_0^\infty\frac{\dd y}{(x+y)^{1+2s}}\dd x\\
 &=\frac{2c_{1,s}}s\int_0^\infty x^{-2s}\phi(x)^2\dd x,
 \end{split}
\end{equation*}
and the remaining parts combine into the kernel $K_s$:
\begin{equation}\label{eq:odd-kinetic-decomposition}
 T_s(\phi)=c_{1,s}\int_0^\infty\int_0^\infty
 K_s(x,y)(\phi(x)-\phi(y))^2\dd y\dd x
 +\frac{2c_{1,s}}s\int_0^\infty x^{-2s}\phi(x)^2\dd x,
\end{equation}
where the pointwise identity
$\frac{(a-b)^2}{|x-y|^{1+2s}}+\frac{(a+b)^2}{(x+y)^{1+2s}}=K_s(x,y)(a-b)^2+\frac{2(a^2+b^2)}{(x+y)^{1+2s}}$
with both terms on the right nonnegative shows that
\eqref{eq:odd-kinetic-decomposition} holds in $[0,\infty]$, and since
$T_s(\phi)<\infty$ both terms are finite.
Let $\widetilde\phi$ be the odd function equal to $|\phi|$ on $(0,\infty)$,
that is, $\widetilde\phi(x)=\operatorname{sgn}(x)|\phi(x)|$.
Since $||a|-|b||\le|a-b|$ and $(\widetilde\phi)^2=\phi^2$,
\eqref{eq:odd-kinetic-decomposition} gives $T_s(\widetilde\phi)\le T_s(\phi)$,
so $\widetilde\phi\in D(V)$, and subtracting the two finite expressions,
\begin{equation}\label{eq:odd-modulus-difference}
 \begin{split}
 \mathcal Q_V[\phi]-\mathcal Q_V[\widetilde\phi]
 &=c_{1,s}\int_0^\infty\int_0^\infty K_s(x,y)
 \Bigl[(\phi(x)-\phi(y))^2-(|\phi(x)|-|\phi(y)|)^2\Bigr]\dd y\dd x\\
 &=2c_{1,s}\int_0^\infty\int_0^\infty K_s(x,y)
       \bigl(|\phi(x)\phi(y)|-\phi(x)\phi(y)\bigr)\dd y\dd x\ge0.
 \end{split}
\end{equation}
Since $K_s>0$, equality holds if and only if $\phi(x)\phi(y)\ge0$ for
almost every $(x,y)\in(0,\infty)^2$, that is, if and only if $\phi$ has
one sign almost everywhere on $(0,\infty)$, the zero function included.
Indeed, if the sets where $\phi>0$ and where $\phi<0$ in $(0,\infty)$
both have positive measure, then the integrand is strictly positive on
their product, which has positive measure.

\emph{Step 2. Nonnegativity.}
Suppose that
\begin{equation*}
 \Lambda_0=\inf\{\mathcal Q_V[\phi]:\phi\in D(V)\text{ odd},\ \|\phi\|_2=1\}<0,
\end{equation*}
and let $\{\phi_j\}_{j\ge1}$ be a minimizing sequence with $\|\phi_j\|_2=1$. Since
$V_-\le C_V$,
$T_s(\phi_j)+\int_\R V_+\phi_j^2\dd x=\mathcal Q_V[\phi_j]+\int_\R V_-\phi_j^2\dd x\le C$
with $C$ independent of $j$, so along a subsequence $\phi_j$ converges
weakly in $H^s$, strongly in $L^2$ on bounded intervals and almost
everywhere to an odd $\phi$. As $V_-$ is bounded with support in
$[-R_0,R_0]$, $\int_\R V_-\phi_j^2\dd x\to\int_\R V_-\phi^2\dd x$, and
Fatou's lemma for the positive part with weak lower semicontinuity of
$T_s$ give
\begin{equation*}
 T_s(\phi)+\int_\R V_+\phi^2\dd x
 \le\liminf_{j\to\infty}\Bigl(\mathcal Q_V[\phi_j]+\int_\R V_-\phi_j^2\dd x\Bigr)
 =\Lambda_0+\int_\R V_-\phi^2\dd x,
\end{equation*}
so $\phi\in D(V)$, $\mathcal Q_V[\phi]\le\Lambda_0<0$, and $\|\phi\|_2\le1$ by weak
lower semicontinuity; in particular $\phi\ne0$. By the definition of
$\Lambda_0$ applied to $\phi/\|\phi\|_2$,
$\Lambda_0\|\phi\|_2^2\le \mathcal Q_V[\phi]\le\Lambda_0$, and since $\Lambda_0<0$ this forces $\|\phi\|_2=1$ and
$\mathcal Q_V[\phi]=\Lambda_0$, so $\phi$ is a minimizer.  By \eqref{eq:odd-modulus-difference}, $\widetilde\phi$ is also a
minimizer, so we may assume $\phi\ge0$ on $(0,\infty)$. For odd
$\psi\in D(V)$ and $\varepsilon\in\R$, the function $\phi+\varepsilon\psi$
is odd and in $D(V)$, and
$\mathcal Q_V[\phi+\varepsilon\psi]\ge\Lambda_0\|\phi+\varepsilon\psi\|_2^2$;
differentiating at $\varepsilon=0$ gives the Euler--Lagrange equation
\begin{equation*}
 \mathcal Q_V[\phi,\psi]=\Lambda_0\int_\R\phi\psi\dd x
 \qquad\text{for every odd }\psi\in D(V).
\end{equation*}
Taking $\psi=z$ and using $\mathcal Q_V[z,\phi]=0$,
\begin{equation*}
 0=\mathcal Q_V[z,\phi]=\mathcal Q_V[\phi,z]=\Lambda_0\int_\R\phi z\dd x
 =2\Lambda_0\int_0^\infty\phi z\dd x<0,
\end{equation*}
because $\Lambda_0<0$, $\phi\ge0$ and $\phi\not\equiv0$ on $(0,\infty)$, and
$z>0$ there. This contradiction proves $\mathcal Q_V\ge0$ on the odd subspace.

\emph{Step 3. The equality case.}
If an odd $\phi\in D(V)$ satisfies $\mathcal Q_V[\phi]=0$, then Step~2 and
\eqref{eq:odd-modulus-difference} give $0\le \mathcal Q_V[\widetilde\phi]\le \mathcal Q_V[\phi]=0$,
so equality holds in \eqref{eq:odd-modulus-difference} and $\phi$ has
one sign on $(0,\infty)$. Define
\begin{equation*}
 c=\frac{\int_\R\phi z\dd x}{\|z\|_2^2},\qquad\psi=\phi-cz.
\end{equation*}
Then $\psi$ is odd, $\int_\R\psi z\dd x=0$, and the equation for $z$ gives
$\mathcal Q_V[\psi]=\mathcal Q_V[\phi]-2c\mathcal Q_V[\phi,z]+c^2\mathcal Q_V[z]=\mathcal Q_V[\phi]=0$; by the same
argument $\psi$ has one sign on $(0,\infty)$, and orthogonality to the
strictly positive function $z$ on that half-line forces $\psi=0$,
that is, $\phi=cz$. Conversely, every multiple of $z$ has zero quadratic
form.
\end{proof}

\begin{proof}[Proof of Theorem~\ref{thm:nondegeneracy}]
After a translation $u$ is symmetric decreasing and $z=-u'>0$ on $(0,\infty)$
(Lemma~\ref{lem:symmetry}); $z\in D(L_u)$ and $L_uz=0$ by
Lemma~\ref{lem:translation-domain}, that is, $\mathcal Q_u[z,\psi]=0$ for every
$\psi\in D_u$; and \eqref{eq:linearized-potential-bounds} verifies the
hypotheses of Lemma~\ref{lem:odd-kernel}, which gives
\begin{equation}\label{eq:odd-kernel-u}
 \mathcal Q_u[\phi]\ge0\quad\text{for odd }\phi\in D_u,
 \qquad \ker L_u\cap L^2_{\rm odd}(\R)=\operatorname{span}\{u'\}.
\end{equation}

\emph{Step 1. Exclusion of an even kernel element.}
Let $\phi\in\ker L_u$ be even. By Remark~\ref{rem:general-direction},
$u\in D_u$ and $\mathcal Q_u[u]<0$, and $\ell(u)>0$. Let $c=\ell(\phi)/\ell(u)$
and $\psi=\phi-cu$; then $\psi\in D_u$ is even with $\ell(\psi)=0$, and
since $\mathcal Q_u[\phi,\cdot]=0$ on $D_u$,
\begin{equation*}
 \mathcal Q_u[\psi]=\mathcal Q_u[\phi]-2c\mathcal Q_u[\phi,u]+c^2\mathcal Q_u[u]=c^2\mathcal Q_u[u]\le0.
\end{equation*}
By \eqref{eq:strict-tangent}, $\psi=0$, so $\phi=cu$, and then
$0=\mathcal Q_u[\phi]=c^2\mathcal Q_u[u]$ with $\mathcal Q_u[u]<0$ gives $c=0$, that is, $\phi=0$.

\emph{Step 2. The full kernel.}
For $\phi\in D_u$ let $\phi_{\rm e}(x)=\frac12(\phi(x)+\phi(-x))$ and
$\phi_{\rm o}(x)=\frac12(\phi(x)-\phi(-x))$. Evenness of $V$ and reflection
invariance of $T_s$ show that both parts belong to $D_u$ and that
\begin{equation}\label{eq:parity-form-splitting}
 \mathcal Q_u[\phi]=\mathcal Q_u[\phi_{\rm e}]+\mathcal Q_u[\phi_{\rm o}],\qquad
 \mathcal Q_u[\phi_{\rm e},\phi_{\rm o}]=0,
\end{equation}
since $\mathcal Q_u[\phi(-\cdot),\psi(-\cdot)]=\mathcal Q_u[\phi,\psi]$
for all $\phi,\psi\in D_u$, by the reflection invariance of $T_s$ and
the evenness of $V$, and therefore
$\mathcal Q_u[\phi_{\rm e},\phi_{\rm o}]
=\mathcal Q_u[\phi_{\rm e},-\phi_{\rm o}]
=-\mathcal Q_u[\phi_{\rm e},\phi_{\rm o}]$. If $\phi\in\ker L_u$, then for every $\eta\in D_u$
the identity $0=\mathcal Q_u[\phi,\eta]=\mathcal Q_u[\phi_{\rm e},\eta_{\rm e}]+\mathcal Q_u[\phi_{\rm o},\eta_{\rm o}]$
with $\eta$ even, respectively odd, shows that $\phi_{\rm e}$ and
$\phi_{\rm o}$ lie in the kernel. Step~1 gives $\phi_{\rm e}=0$ and
\eqref{eq:odd-kernel-u} gives $\phi_{\rm o}\in\operatorname{span}\{u'\}$,
which together with Lemma~\ref{lem:translation-domain} proves
\eqref{eq:full-kernel}.

\emph{Step 3. Morse index.}
The Morse index is the maximal dimension of a subspace of $D_u$ on
which $\mathcal Q_u$ is negative definite. Let $\phi\in D_u$ with
$\ell(\phi)=0$. Since $q_u$ is even, $\ell(\phi_{\rm o})=0$, and hence
$\ell(\phi_{\rm e})=\ell(\phi)=0$. Then \eqref{eq:strict-tangent} gives
$\mathcal Q_u[\phi_{\rm e}]\ge0$, \eqref{eq:odd-kernel-u} gives
$\mathcal Q_u[\phi_{\rm o}]\ge0$, and \eqref{eq:parity-form-splitting}
gives
\begin{equation*}
 \mathcal Q_u[\phi]=\mathcal Q_u[\phi_{\rm e}]+\mathcal Q_u[\phi_{\rm o}]\ge0.
\end{equation*}
Since $\ell$ is a linear functional on $D_u$ with $\ell(u)>0$, the
hyperplane $\{\ell=0\}$ has codimension one in $D_u$. Hence every
two-dimensional subspace of $D_u$ contains a nonzero $\phi$ with
$\ell(\phi)=0$, for which $\mathcal Q_u[\phi]\ge0$, so $\mathcal Q_u$ is
not negative definite on it. Thus the Morse index is at most one.
Since $\mathcal Q_u[u]<0$ by Remark~\ref{rem:general-direction}, it
equals one.
\end{proof}

\begin{proof}[Proof of Theorem~\ref{thm:normalized}]
\emph{Step 1. Reduction to positive multipliers.}
Let $(\lambda_i,u_i)$, $i=1,2$, be as in the statement. Choose
$k>\max\{0,-\lambda_1,-\lambda_2\}$ and define
\begin{equation*}
 \widetilde f(t)=f(t)+kt,\qquad
 \widetilde F(t)=F(t)+\frac k2t^2,\qquad
 \widetilde\lambda_i=\lambda_i+k>0,
\end{equation*}
so that $\widetilde F$ is the primitive of $\widetilde f$ with
$\widetilde F(0)=0$, and $u_i$ is a positive solution of
\begin{equation}\label{eq:shifted-equation}
 (-\Delta)^sv+\widetilde\lambda_iv=\widetilde f(v)\quad\text{in }\R,
 \qquad v\in H^s(\R),
\end{equation}
for $i=1,2$. We check that, for $i=1,2$, equation
\eqref{eq:shifted-equation} satisfies the hypotheses of
Theorem~\ref{thm:main} in case \eqref{eq:caseA} with $r=2$ and $\mu=0$.

First, \eqref{eq:f-regularity} holds for $\widetilde f$. Moreover,
$t\widetilde f'(t)-\widetilde f(t)=tf'(t)-f(t)>0$ for $t>0$ by
\eqref{eq:normalized-assumptions}, and
\begin{equation*}
 \limsup_{t\to0^+}\widetilde f'(t)=k+\limsup_{t\to0^+}f'(t)
 <k+\lambda_i=\widetilde\lambda_i
\end{equation*}
by \eqref{eq:normalized-multiplier}. Hence \eqref{eq:f-assumptions}
holds for $\widetilde f$ with $\lambda$ replaced by $\widetilde\lambda_i$.
Next, the auxiliary function \eqref{eq:shifted-b} for
\eqref{eq:shifted-equation} with $r=2$ and $\mu=0$ is
\begin{equation*}
 \frac{\widetilde f(\sqrt h)-\widetilde\lambda_i\sqrt h}{\sqrt h}
 =\frac{f(\sqrt h)}{\sqrt h}+k-\lambda_i-k=g(h)-\lambda_i,
 \qquad h>0,
\end{equation*}
whose second derivative is $g''$. Since $1/2\le s<1$, we have
$2\le1/(1-s)$, with equality exactly when $s=1/2$, and
\eqref{eq:normalized-concavity} gives the concavity required in
\eqref{eq:caseA}, including the strictness required when $r=1/(1-s)$.

Consequently, Lemmas~\ref{lem:solution-estimates}
and~\ref{lem:symmetry} apply to $u_i$ as a solution of
\eqref{eq:shifted-equation}. After translating $u_1$ and $u_2$, we may
assume that both are symmetric decreasing, and then each $u_i$ is smooth
and positive, with $u_i'<0$ on $(0,\infty)$. Moreover
$\widetilde F(u_i)\in L^1(\R)$ by Lemma~\ref{lem:solution-estimates}(ii),
and since $u_i\in L^2(\R)$, also
$F(u_i)=\widetilde F(u_i)-\frac k2u_i^2\in L^1(\R)$.

\emph{Step 2. Strict minimality for the energy.}
For functions $v$ with $F(v)\in L^1(\R)$ we define the energy
\begin{equation*}
 \mathcal E(v)=\frac12T_s(v)-\int_\R F(v)\dd x.
\end{equation*}
The action
\eqref{eq:action} associated with \eqref{eq:shifted-equation} is
\begin{equation*}
 \widetilde S_i(v)=\frac12T_s(v)+\frac{\widetilde\lambda_i}2\|v\|_2^2
 -\int_\R\widetilde F(v)\dd x
 =\mathcal E(v)+\frac{\lambda_i}2\|v\|_2^2.
\end{equation*}
For $v\in L^2(\R)$ the conditions $F(v)\in L^1$ and
$\widetilde F(v)\in L^1$ are equivalent, so the class $\mathcal C_M$ of
\eqref{eq:comparison-class} with $r=2$ is the same for $f$ and for
$\widetilde f$. By Step~1, $u_1,u_2\in\mathcal C_M$, since
$\mathcal N(u_i)=\|u_i\|_2^2=M$.

Suppose that $u_1\ne u_2$. Lemma~\ref{lem:own-mass}, applied to
\eqref{eq:shifted-equation} with $i=1$, $u=u_1$ and $v=u_2$, gives
$\widetilde S_1(u_2)>\widetilde S_1(u_1)$. Since
$\|u_1\|_2^2=\|u_2\|_2^2=M$, \begin{equation*}
 \widetilde S_1(u_2)-\widetilde S_1(u_1)
 =\mathcal E(u_2)-\mathcal E(u_1)
 +\frac{\lambda_1}2\bigl(\|u_2\|_2^2-\|u_1\|_2^2\bigr)
 =\mathcal E(u_2)-\mathcal E(u_1),
\end{equation*}
and therefore $\mathcal E(u_2)>\mathcal E(u_1)$.
In the same way, Lemma~\ref{lem:own-mass}, applied to
\eqref{eq:shifted-equation} with $i=2$, $u=u_2$ and $v=u_1$, gives
\begin{equation*}
 \mathcal E(u_1)-\mathcal E(u_2)=\widetilde S_2(u_1)-\widetilde S_2(u_2)>0.
\end{equation*}
These two inequalities contradict each other. Hence $u_1=u_2$, and we
write $u=u_1=u_2$. 

\emph{Step 3. Identification of the multiplier.}
The equations are invariant under translations, so the translated
functions still satisfy $(-\Delta)^su+\lambda_iu=f(u)$ for $i=1,2$.
Subtracting the two equations gives $(\lambda_1-\lambda_2)u=0$ in the
sense of distributions, and since $u>0$, we obtain $\lambda_1=\lambda_2$.
Since $u_1$ and $u_2$ were only translated, the original solutions
differ by a translation, which proves Theorem~\ref{thm:normalized}.
\end{proof}

\section{Extensions of the method}\label{sec:further}

Apart from the regularity of positive solutions, the proof of
Theorem~\ref{thm:main} uses only the following three properties of the
action $S$ and of the quantity $\mathcal N$.

(P1) After a translation, every positive solution belongs to a class
$\mathcal K$ of symmetric decreasing functions which contains $cu$ for
every positive solution $u\in\mathcal K$ and every $c>0$, and
$c\mapsto\mathcal N(cu)$ is continuous and nondecreasing on
$(0,\infty)$, tends to $0$ as $c\to0$ and to $\infty$ as $c\to\infty$.

(P2) Every positive solution $u\in\mathcal K$ satisfies $S(cu)<S(u)$
for $c\ne1$.

(P3) Every positive solution $u\in\mathcal K$ is a strict minimizer of
$S$ on $\{v\in\mathcal K:\mathcal N(v)=\mathcal N(u)\}$.

For any problem with these properties, the argument of
Section~\ref{sec:proof-main} shows that positive solutions are unique
up to translations. Here (P2) is Lemma~\ref{lem:amplitude}, (P3) is
Lemma~\ref{lem:own-mass}, and (P1) follows from
Lemmas~\ref{lem:symmetry} and~\ref{lem:amplitude} and the properties of
$c\mapsto\mathcal N(cu)$ established in the proof of
Theorem~\ref{thm:main}. The main step is (P3). It rests on the convexity
of the action along an interpolation which preserves $\mathcal N$
(Sections~\ref{sec:height} and~\ref{sec:two-level}), and on an analysis
at the endpoint, which uses Theorem~\ref{thm:coercivity} in case
\eqref{eq:caseA} and Lemma~\ref{lem:two-level-endpoint} in case
\eqref{eq:caseB}. The convexity of the kinetic energy holds for
$r\le1/(1-s)$ and does not depend on the nonlinearity, which enters the
convexity of the action only through the concavity of $b_{r,\mu}$. We
therefore expect the method to apply to other problems, some of which
we list below. Except for results proved here or in \cite{JZZ2} and
results quoted from the literature, the statements in (i)--(iv) are
conjectures.

(i) \emph{Higher dimensions.} In a forthcoming paper \cite{JZZ2} we
prove the analogues of Theorems~\ref{thm:main}
and~\ref{thm:nondegeneracy} in $\R^N$, $N\ge2$, for every $0<s<1$,
under the analogue of \eqref{eq:caseB} with $1\le r\le N/(N-s)$.
Positive solutions are then unique up to translations, and every
positive solution $u$ is nondegenerate, the kernel of the linearized
operator being spanned by $\partial_1u,\dots,\partial_Nu$, and has Morse
index one. The interaction of two intervals is replaced by a positive
series for the interaction of two concentric balls. Since
$N/(N-s)<2$, the height interpolation used under \eqref{eq:caseA}
is not available, and pure powers $t^p$ are covered only for some
$p\le(N+s)/(N-s)$. In the spirit of the question raised in \cite{Fr}, we
conjecture that for $N\ge2$, $0<s<1$ and $1<p<(N+2s)/(N-2s)$, every
positive solution of $(-\Delta)^su+u=u^p$ in $\R^N$ is a translate of
the ground state of \cite{FLS}.

(ii) \emph{Choquard equations.} Consider
$(-\Delta)^su+\lambda u=(I_\alpha*u^p)\,u^{p-1}$ in $\R$, where
$I_\alpha(x)=|x|^{-(1-\alpha)}$, $0<\alpha<1$ and $p\ge2$. With the
density $\rho=u^p$, the term $\frac\lambda2\|u\|_2^2$ is convex along the
height interpolation because $p\ge2$, the kinetic energy is convex if
$p\le1/(1-s)$ by Lemma~\ref{lem:kinetic-convexity}, and the action is
strictly maximal along the rays $c\mapsto cu$. The nonlocal term is,
up to sign, the interaction energy of $\rho$ for the radially increasing
kernel $-I_\alpha$, and such interaction energies are convex along the
height interpolation \cite{DYY}. We therefore conjecture that for
$1/2\le s<1$ and $2\le p\le1/(1-s)$, positive solutions of this
equation are unique up to translations and nondegenerate. The new
difficulty concerns (P3), namely the coercivity estimate of
Section~\ref{sec:coercivity}, since the second variation contains a
nonlocal term. For $s=1$ and the
Newtonian kernel in $\R^3$, positive solutions are unique up to
translations for $p=2$ \cite{Li77,MZ}.
Uniqueness of ground states for $p$ close to $2$ was proved in
\cite{Xi}. We refer to \cite{MVS} for a survey.

(iii) \emph{Mixed and multi-order operators.} Let
$\mathscr L=a_0(-\partial_x^2)+\sum_{i=1}^ka_i(-\Delta)^{s_i}$ with
$a_0\ge0$, $a_i>0$ and $1/2\le s_i<1$. By
Lemma~\ref{lem:kinetic-convexity}, each $T_{s_i}$ is convex along the
height interpolation of $\rho=u^r$ for $2\le r\le1/(1-\min_is_i)$, and
so is $\|v'\|_2^2$, by letting $s\to1$. Provided (P1) holds, that is, positive
solutions are symmetric decreasing after a translation, we expect
Theorems~\ref{thm:main} and~\ref{thm:nondegeneracy} to hold for
$\mathscr Lu+\lambda u=f(u)$ under \eqref{eq:caseA} with $1/(1-s)$
replaced by $1/(1-\min_is_i)$. Orders below $1/2$ are not covered,
since the exponent $1-2s_i$ of the two-level interpolation depends on
the order.

(iv) \emph{Prescribed $L^2$ norm for $0<s<1/2$.} By
Remark~\ref{rem:normalized}(i), the conclusion of
Theorem~\ref{thm:normalized} holds for $0<s<1/2$ for pure powers in the
range \eqref{eq:pure-power-range}. We conjecture that it also holds for
the nonlinearities \eqref{eq:mixed-f} with $a=0$ in cases (iii)--(v) of
Corollary~\ref{cor:powers}, for which no scaling argument is available,
see Remark~\ref{rem:normalized}(ii). The obstacle for our method is the
following. In the proof of Theorem~\ref{thm:normalized}, two solutions
with different multipliers are compared through (P3) with
$\mathcal N(v)=\|v\|_2^2$, which is the prescribed quantity and does not
depend on the multiplier. In case \eqref{eq:caseB}, however,
$\mathcal N(v)=\int_\R(v-\mu)_+^r\dd x$ depends on the multiplier through
the level $\mu=t_\lambda$, so that this comparison is not available.

\appendix
\Needspace{14\baselineskip}
\section{Convexity lemmas}\label{app:convexity}

The following three lemmas contain the algebraic calculations used in
Sections~\ref{sec:height}--\ref{sec:coercivity}.

\begin{lemma}\label{lem:monomial}
Let $m\ge1$ be an integer, $q_i\ge0$ for $1\le i\le m$, and
$P\ge1+\sum_{i=1}^m q_i$. Then
\begin{equation*}
 K(X,Y_1,\ldots,Y_m)=X^P\prod_{i=1}^mY_i^{-q_i}
\end{equation*}
is convex on $(0,\infty)^{m+1}$.
\end{lemma}
\begin{proof}
Let $X(t)=X_0+ta$ and $Y_i(t)=Y_{i,0}+tb_i$ be any affine
paths that remain positive on an interval $I$. Define
$k(t)=K(X(t),Y_1(t),\ldots,Y_m(t))$ and
$\sigma=\sum_{i=1}^m q_i$. If $P=1$, the hypothesis forces
$\sigma=0$, so $k=X$ is affine. We may therefore assume $P>1$.
Taking logarithms and differentiating gives
\begin{equation*}
 \frac{d}{dt}\log k(t)=\frac{Pa}{X(t)}
                         -\sum_{i=1}^m\frac{q_i b_i}{Y_i(t)},
\end{equation*}
\begin{equation*}
 \frac{d^2}{dt^2}\log k(t)=-\frac{Pa^2}{X(t)^2}
                         +\sum_{i=1}^m\frac{q_i b_i^2}{Y_i(t)^2}.
\end{equation*}
Since $k=\exp(\log k)$, the product and chain rules yield
\begin{equation*}
 \frac{k''(t)}{k(t)}
 =\left(\frac{d}{dt}\log k(t)\right)^2
                              +\frac{d^2}{dt^2}\log k(t) =\left(\frac{Pa}{X(t)}-\sum_{i=1}^m\frac{q_i b_i}{Y_i(t)}\right)^2
       -\frac{Pa^2}{X(t)^2}+\sum_{i=1}^m\frac{q_i b_i^2}{Y_i(t)^2}.
\end{equation*}
Write $\alpha=a/X(t)$, $\Sigma=\sum_{i=1}^mq_ib_i/Y_i(t)$ and
$T=\sum_{i=1}^mq_ib_i^2/Y_i(t)^2$. Expanding and completing the square in
$\alpha$,
\begin{equation*}
 \begin{split}
 \frac{k''}{k}&=(P\alpha-\Sigma)^2-P\alpha^2+T
 =P(P-1)\alpha^2-2P\alpha\Sigma+\Sigma^2+T\\
 &=P(P-1)\Bigl(\alpha-\frac{\Sigma}{P-1}\Bigr)^2-\frac{P}{P-1}\Sigma^2+\Sigma^2+T,
 \end{split}
\end{equation*}
and $-\frac P{P-1}\Sigma^2+\Sigma^2=-\frac1{P-1}\Sigma^2$. This gives the
exact identity
\begin{equation*}
 \frac{k''(t)}{k(t)}
 ={}P(P-1)\left(\frac a{X(t)}
             -\frac1{P-1}\sum_{i=1}^m\frac{q_i b_i}{Y_i(t)}\right)^2 +\sum_{i=1}^m\frac{q_i b_i^2}{Y_i(t)^2}
       -\frac1{P-1}\left(\sum_{i=1}^m\frac{q_i b_i}{Y_i(t)}\right)^2.
\end{equation*}
Weighted Cauchy--Schwarz gives, also when some or all $q_i$ vanish,
\begin{equation*}
 \left(\sum_{i=1}^m\frac{q_i b_i}{Y_i(t)}\right)^2
 \le\sigma\sum_{i=1}^m\frac{q_i b_i^2}{Y_i(t)^2}.
\end{equation*}
Dropping the first square and using $P-1\ge\sigma$ therefore gives
\begin{equation}\label{eq:monomial-quantitative}
 \frac{k''(t)}{k(t)}
 \ge\frac{P-1-\sigma}{P-1}
       \sum_{i=1}^m\frac{q_i b_i^2}{Y_i(t)^2}\ge0.
\end{equation}
Thus $k''\ge0$ on every affine segment in the positive orthant,
which proves convexity.
\end{proof}

\begin{lemma}\label{lem:strict-monomial}
Let $s>\beta>0$, and let
\begin{equation*}
 X(t)=X_0+ta,\qquad Y_1(t)=Y_{1,0}+tb,\qquad Y_2(t)=Y_{2,0}+tc
\end{equation*}
be positive on an open interval $I$, where $a,b,c\in\R$. Define
\begin{equation*}
 k_0(t)=X(t)^{1+2s}Y_1(t)^{-\beta}Y_2(t)^{-\beta}.
\end{equation*}
Then, for every $t\in I$,
\begin{equation}\label{eq:strict-monomial}
 \frac1{k_0(t)}\frac{d^2k_0}{dt^2}(t)
 \ge\frac{\beta(s-\beta)}s
       \left(\frac{b^2}{Y_1(t)^2}+\frac{c^2}{Y_2(t)^2}\right).
\end{equation}
\end{lemma}
\begin{proof}
Apply \eqref{eq:monomial-quantitative} with $m=2$, $P=1+2s$,
and $q_1=q_2=\beta$. Here
$\sigma=2\beta<2s=P-1$, so
\begin{equation*}
 \frac{k_0''(t)}{k_0(t)}
 \ge\frac{2s-2\beta}{2s}
       \left(\frac{\beta b^2}{Y_1(t)^2}+\frac{\beta c^2}{Y_2(t)^2}\right)
 =\frac{\beta(s-\beta)}s
       \left(\frac{b^2}{Y_1(t)^2}+\frac{c^2}{Y_2(t)^2}\right).
\end{equation*}
This is \eqref{eq:strict-monomial}.
\end{proof}

\begin{lemma}\label{lem:nonpositive-monomial}
Let $m\ge1$ be an integer and $q_0,\ldots,q_m\ge0$. Then
\begin{equation*}
 K(Y_0,\ldots,Y_m)=\prod_{i=0}^mY_i^{-q_i}
\end{equation*}
is convex on $(0,\infty)^{m+1}$. More precisely, if
$Y_i(t)=Y_{i,0}+tb_i$ are positive on an interval $I$ and
$k(t)=K(Y_0(t),\ldots,Y_m(t))$, then
\begin{equation}\label{eq:nonpositive-monomial}
 \frac{k''(t)}{k(t)}\ge\sum_{i=0}^m\frac{q_ib_i^2}{Y_i(t)^2}
 \qquad(t\in I).
\end{equation}
In particular, if $q_0,q_1>0$ and $(b_0,b_1)\ne(0,0)$, then $k''>0$ on $I$.
\end{lemma}
\begin{proof}
As in the proof of Lemma~\ref{lem:monomial},
\begin{equation*}
 \frac{k''(t)}{k(t)}=\left(\frac d{dt}\log k(t)\right)^2
 +\frac{d^2}{dt^2}\log k(t)
 =\left(\sum_{i=0}^m\frac{q_ib_i}{Y_i(t)}\right)^2
 +\sum_{i=0}^m\frac{q_ib_i^2}{Y_i(t)^2},
\end{equation*}
which gives \eqref{eq:nonpositive-monomial}.
\end{proof}

\section{Verification of the examples}\label{app:examples}

In this appendix we prove Corollaries~\ref{cor:powers}
and~\ref{cor:concave} and the assertions of Remarks~\ref{rem:pure-powers}
and~\ref{rem:saturable}. 

\begin{lemma}\label{lem:criterion}
Let $\mu\ge0$, $r\ge1$, and
$\Phi(\tau)=\lambda(\mu+\tau)-f(\mu+\tau)$ for $\tau>0$. Then
$b_{r,\mu}''\le0$ on $(0,\infty)$ if and only if
\begin{equation}\label{eq:Leps}
 \mathcal L_r\Phi(\tau):=\tau^2\Phi''(\tau)
 -3(r-1)\tau\Phi'(\tau)+(r-1)(2r-1)\Phi(\tau)\ge0
 \qquad(\tau>0),
\end{equation}
and $b_{r,\mu}''<0$ on $(0,\infty)$ if and only if
$\mathcal L_r\Phi>0$ on $(0,\infty)$. The operator
$\mathcal L_r=(\tau\partial_\tau-r+1)(\tau\partial_\tau-2r+1)$
is linear, commutes with the dilations $\tau\mapsto c\tau$, and
satisfies
$\mathcal L_r\tau^\alpha=(\alpha-r+1)(\alpha-2r+1)\tau^\alpha$.
\end{lemma}
\begin{proof}
With $h=\tau^r$ we have $h^{1/r}=\tau$ and $h^\beta=\tau^{r-1}$,
so $b_{r,\mu}(h)=-\tau^{1-r}\Phi(\tau)$, and $\dd h=r\tau^{r-1}\dd\tau$
gives $\frac d{dh}=(r\tau^{r-1})^{-1}\frac d{d\tau}$. Applying this twice,
\begin{equation*}
 b_{r,\mu}''(h)=-\frac{\tau^{1-r}}{r^2}\frac d{d\tau}
 \Bigl(\tau^{1-r}\frac d{d\tau}\bigl(\tau^{1-r}\Phi\bigr)\Bigr).
\end{equation*}
Here
\begin{align*}
 \tau^{1-r}\frac{d}{d\tau}\bigl(\tau^{1-r}\Phi\bigr)
 &=\tau^{1-2r}\bigl(\tau\Phi'-(r-1)\Phi\bigr),\\
 \frac{d}{d\tau}\Bigl[\tau^{1-2r}\bigl(\tau\Phi'-(r-1)\Phi\bigr)\Bigr]
 &=-(2r-1)\tau^{-2r}\bigl(\tau\Phi'-(r-1)\Phi\bigr)
   +\tau^{1-2r}\bigl(\tau\Phi''+(2-r)\Phi'\bigr)\\
 &=\tau^{-2r}\bigl[\tau^2\Phi''-3(r-1)\tau\Phi'
   +(r-1)(2r-1)\Phi\bigr]
  =\tau^{-2r}\,\mathcal L_r\Phi(\tau),
\end{align*}
since $-(2r-1)+(2-r)=-3(r-1)$. Hence
\begin{equation*}
 b_{r,\mu}''(h)=-\frac{\tau^{1-3r}}{r^2}\,\mathcal L_r\Phi(\tau),
\end{equation*}
and the sign claims follow because $\tau^{1-3r}/r^2>0$. The
factorization of $\mathcal L_r$ follows from
$(\tau\partial_\tau)^2=\tau^2\partial_\tau^2+\tau\partial_\tau$, which gives
$(\tau\partial_\tau-r+1)(\tau\partial_\tau-2r+1)
=\tau^2\partial_\tau^2+\tau\partial_\tau-(3r-2)\tau\partial_\tau+(r-1)(2r-1)$,
and $\mathcal L_r\tau^\alpha=(\alpha-r+1)(\alpha-2r+1)\tau^\alpha$ because
$\tau\partial_\tau\tau^\alpha=\alpha\tau^\alpha$.
\end{proof}

For $p>1$ define $\Phi_p(\tau)=(1+\tau)-(1+\tau)^p$ for $\tau\ge0$, and
let $\Psi(\tau)=(1+\tau)\log(1+\tau)$. The function $\Phi_p$ is the
function $\Phi$ of Lemma~\ref{lem:criterion} for $f(t)=t^p$ with
$\mu=t_\lambda=1$, and $\Psi$ is, up to a constant factor, the
contribution of the logarithmic term in \eqref{eq:mixed-f}; see
\eqref{eq:L-split} below. We write $y=1+\tau$. Since $\Phi_p'=1-py^{p-1}$
and $\Phi_p''=-p(p-1)y^{p-2}$,
\begin{equation*}
 \mathcal L_r\Phi_p
 =-p(p-1)\tau^2y^{p-2}-3(r-1)\tau+3(r-1)p\,\tau y^{p-1}
  +(r-1)(2r-1)(y-y^p),
\end{equation*}
and the two terms without a factor $y^{p-2}$ combine to
$(r-1)(2r-1)(1+\tau)-3(r-1)\tau
=(r-1)(2r-1)-2(r-1)(2-r)\tau$; hence
\begin{equation}\label{eq:Lphip}
 \mathcal L_r\Phi_p(\tau)=(r-1)(2r-1)
 -2(r-1)(2-r)\tau
 -y^{p-2}\bigl[p(p-1)\tau^2-3(r-1)p\,\tau y
 +(r-1)(2r-1)y^2\bigr].
\end{equation}
Likewise $\Psi'=1+\log y$ and $\Psi''=1/y$, so
$\mathcal L_r\Psi=\frac{\tau^2}y-3(r-1)\tau(1+\log y)
+(r-1)(2r-1)y\log y$, and collecting the terms with $\log y$,
\begin{equation}\label{eq:Lpsi}
 \mathcal L_r\Psi(\tau)=\frac{\tau^2}y-3(r-1)\tau
 +(r-1)\bigl[(2r-1)-2(2-r)\tau\bigr]\log y,
\end{equation}
where $(2r-1)(1+\tau)-3\tau=(2r-1)-2(2-r)\tau$.

\begin{lemma}\label{lem:monotone-p}
Let $1<r<2$, $\gamma_1=2(r-1)(2-r)$ and
$\gamma_2=(r-1)(2r-1)$.
\begin{enumerate}
\item[(i)] For $p>1$ and $\tau>0$,
\begin{equation}\label{eq:Jp}
\begin{gathered}
 \mathcal L_r\Phi_p(\tau)=p(p-1)\tau^2(1+\tau)^{p-2}j_p(\tau),\\[2pt]
 j_p(\tau)=\frac{\gamma_1(1+\tau)^{2-p}}{p\tau}
 +\int_0^1(\gamma_1+\gamma_2x)\Bigl(\frac{1+\tau}{1+\tau x}\Bigr)^{2-p}\dd x-1,
\end{gathered}
\end{equation}
and for every $\tau>0$ the function $p\mapsto j_p(\tau)$ is
nonincreasing on $(1,\infty)$. Consequently, if
$\mathcal L_r\Phi_{\bar p}\ge0$ on $(0,\infty)$ for some $\bar p>1$,
then $\mathcal L_r\Phi_p\ge0$ on $(0,\infty)$ for every
$p\in(1,\bar p]$.
\item[(ii)] If $\mathcal L_r\Phi_{\bar p}\ge0$ on $(0,\infty)$
for some $\bar p>1$, then $\mathcal L_r\Psi\le0$ on $(0,\infty)$.
\end{enumerate}
\end{lemma}
\begin{proof}
(i) Taylor's formula gives $(1+\tau)^p=1+p\tau+p(p-1)Z(\tau)$ with
$Z(\tau)=\int_0^\tau(\tau-\sigma)(1+\sigma)^{p-2}\dd\sigma$, so that
$\Phi_p=-(p-1)(\tau+pZ)$. By Lemma~\ref{lem:criterion},
$\mathcal L_r\tau=(2-r)(2-2r)\tau=-\gamma_1\tau$.
Since $Z'(\tau)=\int_0^\tau(1+\sigma)^{p-2}\dd\sigma$ and
$Z''(\tau)=(1+\tau)^{p-2}$,
\begin{equation*}
 \begin{split}
 \mathcal L_r Z(\tau)
 &=\tau^2(1+\tau)^{p-2}-3(r-1)\tau\int_0^\tau(1+\sigma)^{p-2}\dd\sigma
   +\gamma_2\int_0^\tau(\tau-\sigma)(1+\sigma)^{p-2}\dd\sigma\\
 &=\tau^2(1+\tau)^{p-2}
  -\int_0^\tau\bigl((3(r-1)-\gamma_2)\tau+\gamma_2\sigma\bigr)(1+\sigma)^{p-2}\dd\sigma,
 \end{split}
\end{equation*}
and $3(r-1)-\gamma_2=(r-1)(3-2r+1)=\gamma_1$. Therefore
$\mathcal L_r\Phi_p=-(p-1)(\mathcal L_r\tau+p\,\mathcal L_r Z)$
equals
\begin{equation*}
 \mathcal L_r\Phi_p(\tau)=(p-1)\Bigl[\gamma_1\tau
 +p\int_0^\tau(\gamma_1\tau+\gamma_2\sigma)(1+\sigma)^{p-2}\dd\sigma
 -p\tau^2(1+\tau)^{p-2}\Bigr].
\end{equation*}
The substitution $\sigma=\tau x$ turns the integral into
$\tau^2\int_0^1(\gamma_1+\gamma_2x)(1+\tau x)^{p-2}\dd x$, and factoring
$p(p-1)\tau^2(1+\tau)^{p-2}$ out of the bracket gives \eqref{eq:Jp}: the
first term becomes $\gamma_1(1+\tau)^{2-p}/(p\tau)$, the integral becomes
$\int_0^1(\gamma_1+\gamma_2x)((1+\tau)/(1+\tau x))^{2-p}\dd x$, and the last
term becomes $-1$. Since
$1+\tau\ge1+\tau x\ge1$ for $x\in[0,1]$, the functions
$p\mapsto(1+\tau)^{2-p}$, $p\mapsto((1+\tau)/(1+\tau x))^{2-p}$ and
$p\mapsto1/p$ are nonincreasing, and $\gamma_1,\gamma_2>0$; this proves
the monotonicity, and the last assertion follows because the factor
$p(p-1)\tau^2(1+\tau)^{p-2}$ in \eqref{eq:Jp} is positive.

(ii) By (i), $\mathcal L_r\Phi_p\ge0$ for every $p\in(1,\bar p]$.
As $p\to1$, for every $\tau>0$,
$\Phi_p(\tau)/(p-1)\to-\Psi(\tau)$,
$\Phi_p'(\tau)/(p-1)=(1-p(1+\tau)^{p-1})/(p-1)\to-\Psi'(\tau)$ and
$\Phi_p''(\tau)/(p-1)=-p(1+\tau)^{p-2}\to-\Psi''(\tau)$, because
$\Psi'=1+\log(1+\tau)$ and $\Psi''=1/(1+\tau)$. Hence
$\mathcal L_r\Phi_p(\tau)/(p-1)\to-\mathcal L_r\Psi(\tau)$,
and the limit of nonnegative numbers is nonnegative.
\end{proof}

\begin{lemma}\label{lem:explicit}
The following inequalities hold on $(0,\infty)$:
(i) $\mathcal L_r\Phi_{2r-1}>0$ if
$3/2\le r<2$;
(ii) $\mathcal L_{4/3}\Phi_{3/2}>0$;
(iii) $\mathcal L_{5/4}\Phi_{5/4}>0$;
(iv) $\mathcal L_{6/5}\Psi<0$.
\end{lemma}
\begin{proof}
We write $y=1+\tau$.

(i) Let $p=2r-1\in[2,3)$, so that $r-1=(p-1)/2$. Then
$(r-1)(2r-1)=(r-1)p=\frac12p(p-1)$,
$3(r-1)p=\frac32p(p-1)$ and
$2(r-1)(2-r)=(r-1)(3-p)$, so the bracket
in \eqref{eq:Lphip} equals
$\frac12p(p-1)(2\tau^2-3\tau y+y^2)=(r-1)p(y-\tau)(y-2\tau)=(r-1)p(1-\tau)$,
because $y-\tau=1$. Hence \eqref{eq:Lphip} reads
$\mathcal L_r\Phi_p=(r-1)p-(r-1)(3-p)\tau-(r-1)p(1-\tau)y^{p-2}
=(r-1)[p-(3-p)\tau+p(\tau-1)y^{p-2}]$.
If $\tau\ge1$, then $y^{p-2}\ge1$ and the bracket is at least
$(2p-3)\tau>0$. If $0<\tau<1$, concavity of $y\mapsto y^{p-2}$ gives
$y^{p-2}\le1+(p-2)\tau$, so the bracket is at least
$\tau[(2p-3)-p(p-2)(1-\tau)]\ge\tau(p-1)(3-p)>0$.

(ii) Here $r=\frac43$, $p=\frac32$,
$(r-1)(2r-1)=\frac59$, $2(r-1)(2-r)=\frac49$,
$p(p-1)=\frac34$ and $3(r-1)p=\frac32$. With $\tau=y-1$ the bracket in
\eqref{eq:Lphip} is
$\frac34(y-1)^2-\frac32(y-1)y+\frac59y^2=-\frac7{36}y^2+\frac34$, so
\begin{equation*}
 \mathcal L_{4/3}\Phi_{3/2}
 =\frac59-\frac49(y-1)-y^{-1/2}\Bigl(-\frac7{36}y^2+\frac34\Bigr)
 =1-\frac49y+\frac7{36}y^{3/2}-\frac34y^{-1/2}.
\end{equation*}
With $z=y^{1/2}>1$ and $w=z-1>0$,
$36z\,\mathcal L_{4/3}\Phi_{3/2}=7z^4-16z^3+36z-27=w\,(7w^3+12w^2-6w+16)>0$,
because $12w^2-6w+16>0$; the factorization is checked by expanding
$7(w+1)^4-16(w+1)^3+36(w+1)-27=7w^4+12w^3-6w^2+16w$.

(iii) Here $r=\frac54$, $p=\frac54$,
$(r-1)(2r-1)=2(r-1)(2-r)=\frac38$, $p(p-1)=\frac5{16}$ and
$3(r-1)p=\frac{15}{16}$. The bracket in \eqref{eq:Lphip} is
$\frac5{16}(y-1)^2-\frac{15}{16}(y-1)y+\frac38y^2=-\frac14y^2+\frac5{16}y+\frac5{16}$,
so
\begin{equation*}
 16\,\mathcal L_{5/4}\Phi_{5/4}
 =16\Bigl[\frac38-\frac38(y-1)-y^{-3/4}\Bigl(-\frac14y^2+\frac5{16}y+\frac5{16}\Bigr)\Bigr]
 =12-6y+4y^{5/4}-5y^{1/4}-5y^{-3/4}.
\end{equation*}
With $z=y^{1/4}>1$ and $w=z-1>0$,
$16z^3\mathcal L_{5/4}\Phi_{5/4}=4z^8-6z^7-5z^4+12z^3-5
=w\,(4w^7+26w^6+70w^5+98w^4+65w^3+6w^2-8w+6)>0$,
because $6w^2-8w+6>0$.  The factorization follows by expanding
$4(w+1)^8-6(w+1)^7-5(w+1)^4+12(w+1)^3-5$ with the binomial theorem, which
gives $4w^8+26w^7+70w^6+98w^5+65w^4+6w^3-8w^2+6w$.

(iv) Here $r=\frac65$, and \eqref{eq:Lpsi} gives
$25\,\mathcal L_{6/5}\Psi=25\tau^2/y-15\tau+(7-8\tau)\log y$. For
$\tau\ge0$ we have $\frac{2\tau}{2+\tau}\le\log y\le\frac{\tau(2+\tau)}{2(1+\tau)}$:
the three expressions vanish at $\tau=0$ and their derivatives are
ordered, since
$\frac1{1+\tau}-\frac4{(2+\tau)^2}=\frac{\tau^2}{(1+\tau)(2+\tau)^2}$ and
$\frac{2+2\tau+\tau^2}{2(1+\tau)^2}-\frac1{1+\tau}=\frac{\tau^2}{2(1+\tau)^2}$.
If $\tau\ge7/8$, then $7-8\tau\le0$, so the lower bound for $\log y$
gives $(7-8\tau)\log y\le(7-8\tau)\frac{2\tau}{2+\tau}$ and
\begin{equation*}
 25\,\mathcal L_{6/5}\Psi
 \le\frac{25\tau^2(2+\tau)-15\tau(1+\tau)(2+\tau)+2\tau(7-8\tau)(1+\tau)}{(1+\tau)(2+\tau)}
 =\frac{\tau(-6\tau^2+3\tau-16)}{(1+\tau)(2+\tau)}<0.
\end{equation*}
If $0<\tau<7/8$, then $7-8\tau>0$, so the upper bound gives
$(7-8\tau)\log y\le(7-8\tau)\frac{\tau(2+\tau)}{2(1+\tau)}$ and
\begin{equation*}
 25\,\mathcal L_{6/5}\Psi
 \le\frac{50\tau^2-30\tau(1+\tau)+\tau(7-8\tau)(2+\tau)}{2(1+\tau)}
 =\frac{\tau(-8\tau^2+11\tau-16)}{2(1+\tau)}<0.
\end{equation*}
In both cases the quadratic polynomial in the numerator is negative,
since both quadratics have negative discriminant and negative leading
coefficient.
\end{proof}

\begin{proof}[Proof of Corollary~\ref{cor:powers}]
\emph{Step 1. The structural conditions.}
Let $f$ be given by \eqref{eq:mixed-f}. Then
$f\in C([0,\infty))\cap C^\infty((0,\infty))$, $f(0)=0$, and
$f'(t)=a(\log t^2+2)+\sum_jc_jp_jt^{p_j-1}$, so that
$\limsup_{t\to0^+}f'(t)$ equals $-\infty$ if $a>0$ and $0$ if $a=0$;
moreover $tf'(t)-f(t)=2at+\sum_jc_j(p_j-1)t^{p_j}>0$ for $t>0$. Thus
\eqref{eq:f-regularity} and \eqref{eq:f-assumptions} hold for every
$\lambda>0$. Since $f(t)/t\to\infty$ as $t\to\infty$, the level
$t_\lambda$ exists, and
\begin{equation}\label{eq:lambda-split}
 \lambda=\frac{f(t_\lambda)}{t_\lambda}
 =a\log t_\lambda^2+\sum_{j=1}^mc_jt_\lambda^{\,p_j-1}.
\end{equation}

\emph{Step 2. Cases (i) and (ii).}
Let $1/2\le s<1$ and $\mu=0$. For $r\ge2$,
\begin{equation*}
 b_r(h)=\frac{2a}r\,h^{2/r-1}\log h+\sum_{j=1}^mc_jh^{\alpha_j}-\lambda h^{\gamma},
 \qquad \alpha_j=\frac{p_j-r+1}r,\quad \gamma=\frac{2-r}r.
\end{equation*}
In case (i), $r=2$, $\gamma=0$ and $\alpha_j=(p_j-1)/2\in(0,1]$, so
\begin{equation}\label{eq:b2-mixed}
 b_2(h)=a\log h+\sum_{j=1}^mc_jh^{\alpha_j}-\lambda,\qquad
 b_2''(h)=-\frac a{h^2}+\sum_{j=1}^mc_j\alpha_j(\alpha_j-1)h^{\alpha_j-2},
\end{equation}
which is nonpositive because $0<\alpha_j\le1$, and negative because
$a>0$. Hence \eqref{eq:caseA} holds with $r=2\le1/(1-s)$. In case (ii), $a=0$,
$0\le\alpha_j\le1$ by $r-1\le p_j\le2r-1$, and $\gamma\le0$, so
$b_r''(h)=\sum_jc_j\alpha_j(\alpha_j-1)h^{\alpha_j-2}-\lambda\gamma(\gamma-1)h^{\gamma-2}\le0$.
If $r>2$, the last term is negative and $b_r''<0$; if $r=2$, then
$b_r''<0$ as soon as some $p_j<3$. Thus \eqref{eq:caseA} holds, the
strictness at $r=1/(1-s)$ being automatic for $r>2$ and required for
$r=2$ only when $s=1/2$, where some $p_j<3$ is assumed.

\emph{Step 3. Cases (iii)--(vi): reduction to the model functions.}
Let $0<s<1/2$, $\mu=t_\lambda$ and $r\ge1$, where $r$ is chosen in
Step~4, and let $\Phi(\tau)=\lambda(\mu+\tau)-f(\mu+\tau)$. By \eqref{eq:lambda-split}
and $\lambda\mu=f(\mu)$,
\begin{equation*}
 \begin{split}
 \Phi(\tau)&=a(\mu+\tau)\bigl(\log\mu^2-\log(\mu+\tau)^2\bigr)
 +\sum_{j=1}^mc_j\bigl(\mu^{p_j-1}(\mu+\tau)-(\mu+\tau)^{p_j}\bigr)\\
 &=-2a\mu\,\Psi(\tau/\mu)+\sum_{j=1}^mc_j\mu^{p_j}\Phi_{p_j}(\tau/\mu),
 \end{split}
\end{equation*}
and since $\mathcal L_r$ is linear and commutes with dilations,
\begin{equation}\label{eq:L-split}
 \mathcal L_r\Phi(\tau)
 =-2a\mu\,(\mathcal L_r\Psi)(\tau/\mu)
 +\sum_{j=1}^mc_j\mu^{p_j}(\mathcal L_r\Phi_{p_j})(\tau/\mu).
\end{equation}
By Lemma~\ref{lem:criterion}, and since $\mu=t_\lambda>0$ and
$f(\mu)=\lambda\mu$ are already in place, \eqref{eq:caseB} holds as soon as
$r\le1/(1-s)$, $\mathcal L_r\Phi_{p_j}\ge0$ on $(0,\infty)$ for
every $j$, and $\mathcal L_r\Psi\le0$ on $(0,\infty)$ if $a>0$;
by Lemma~\ref{lem:monotone-p}, all of these follow from
$\mathcal L_r\Phi_{\bar p}\ge0$ for a single $\bar p\ge p_*$, $\bar p>1$,
and if $m=0$ it suffices that $\mathcal L_r\Psi\le0$.

\emph{Step 4. Cases (iii)--(vi): the choice of $r$.}
In case (iii) we take $r=\max\{3/2,(p_*+1)/2\}$. Then $r\in[3/2,2)$ because
$p_*\le(1+s)/(1-s)<3$; $r\le1/(1-s)$ reads $s\ge1/3$ if $p_*\le2$ and
$p_*\le(1+s)/(1-s)$ if $p_*>2$; and since $p_*\le2r-1$,
Lemma~\ref{lem:explicit}(i) applies with $\bar p=2r-1$. In case
(iv) we take $r=4/3$. Then $r\le1/(1-s)$ because $s\ge1/4$, and
Lemma~\ref{lem:explicit}(ii) applies with $\bar p=3/2\ge p_*$. In case (v) we take
$r=5/4$. Then $r\le1/(1-s)$ because $s\ge1/5$, and
Lemma~\ref{lem:explicit}(iii) applies with $\bar p=5/4\ge p_*$. In case
(vi) we have $m=0$, and we take $r=6/5$. Then $r\le1/(1-s)$ because
$s\ge1/6$, and
Lemma~\ref{lem:explicit}(iv) gives $\mathcal L_{6/5}\Psi<0$. In all four
cases Step~3 gives \eqref{eq:caseB} with $\mu=t_\lambda$, and
Theorems~\ref{thm:main} and~\ref{thm:nondegeneracy} give uniqueness up
to translations, nondegeneracy and Morse index one.

\emph{Step 5. Pure powers (Remark~\ref{rem:pure-powers}).}
For $f(t)=t^p$ take $a=0$, $m=1$, $c_1=1$. For $1/2\le s<1$ and
$1<p\le(1+s)/(1-s)$, let $r=\max\{2,(p+1)/2\}$; then $r\le1/(1-s)$
(as $r=2$ if $p\le3$ and $r=(p+1)/2$ if $p>3$), $r-1\le p\le2r-1$, and
$p<3$ if $s=1/2$ by assumption, so case (ii) applies. For
$1/3\le s<1/2$, $1/4\le s<1/3$ and $1/5\le s<1/4$, cases (iii), (iv) and
(v) apply with $t_\lambda=\lambda^{1/(p-1)}$. Ground states exist by
\cite{FL}: for $s\ge1/2$ every $p>1$ is $H^s$-subcritical, and for
$s<1/2$ the $H^s$-subcritical range $1<p<(1+2s)/(1-2s)$ contains the stated
exponents, since $(1+s)/(1-s)<(1+2s)/(1-2s)$, $3/2<3\le(1+2s)/(1-2s)$
for $s\ge1/4$, and $5/4<7/3\le(1+2s)/(1-2s)$ for $s\ge1/5$. By
Theorem~\ref{thm:main}, every positive solution is a translate of the
ground state.

\emph{Step 6. Prescribed $L^2$ norm.}
In case (i), and in case (ii) with $p_j\le3$ for all $j$, the function
$g(h)=f(\sqrt h)/\sqrt h=a\log h+\sum_jc_jh^{\alpha_j}$,
$\alpha_j=(p_j-1)/2$, satisfies $g''\le0$, and $g''<0$ if $a>0$ or some
$p_j<3$, by \eqref{eq:b2-mixed}; and \eqref{eq:normalized-assumptions}
holds by Step~1, with $\limsup_{t\to0^+}f'(t)=-\infty$ if $a>0$ and $=0$
if $a=0$. Hence Theorem~\ref{thm:normalized} applies for every real
multiplier if $a>0$, and for every positive multiplier if $a=0$.
\end{proof}

\begin{proof}[Proof of Corollary~\ref{cor:concave}]
For $r=1$ and $\mu=t_\lambda$ we have
$b_{1,t_\lambda}(h)=f(t_\lambda+h)-\lambda(t_\lambda+h)$, which is
concave on $(0,\infty)$ if and only if $f$ is concave on
$[t_\lambda,\infty)$. Since $1\le1/(1-s)$ and $f(t_\lambda)=\lambda t_\lambda$,
condition \eqref{eq:caseB} holds for every $0<s<1/2$, and
Theorems~\ref{thm:main} and~\ref{thm:nondegeneracy} apply.

We now verify the assertions of Remark~\ref{rem:saturable}. Let
$f(t)=ct^3/(1+t^2)$ with $c>0$. Then $f\ge0$, $f'(0)=0<\lambda$, and
$tf'(t)-f(t)=2ct^3/(1+t^2)^2>0$ for $t>0$, so \eqref{eq:f-regularity}
and \eqref{eq:f-assumptions} hold for every $\lambda>0$. If
$\lambda\ge c$, then $f(t)<\lambda t$ for all $t>0$, and there is no
positive solution by Lemma~\ref{lem:solution-estimates}(v), whose
hypothesis \eqref{eq:small-growth} holds with $p_0=1$ if $s<1/2$, while
for $s\ge1/2$ condition \eqref{eq:caseA} holds as shown below. If
$\lambda<c$, then $t_\lambda^2=\lambda/(c-\lambda)$, and
$f''(t)=-2ct(t^2-3)/(1+t^2)^3$ is negative for $t>\sqrt3$ and positive
for $0<t<\sqrt3$, so $f$ is concave on $[t_\lambda,\infty)$ if and only
if $t_\lambda\ge\sqrt3$, that is, $\lambda\ge3c/4$. For $1/2\le s<1$ and $r=2$, $b_2(h)=ch/(1+h)-\lambda$
and $b_2''(h)=-2c/(1+h)^3<0$, so \eqref{eq:caseA} holds for every
$\lambda>0$.
\end{proof}

\section{Proofs of the preliminary lemmas}\label{app:prelim}

This appendix contains the proofs of
Lemmas~\ref{lem:growth}--\ref{lem:symmetry}. 

\begin{proof}[Proof of Lemma~\ref{lem:growth}]
Since $(f(t)/t)'=(tf'(t)-f(t))/t^2>0$, we have $f(t)\ge f(1)t$ for
$t\ge1$. Concavity gives $b_r(h)\le b_r(1)+b_r'(1)(h-1)$ and
$b_r'(h)\le b_r'(1)$ for $h\ge1$, so that
\begin{gather*}
 f(1)t\le f(t)=\lambda t+t^{r-1}b_r(t^r)\le Ct^{2r-1},\\
 f(1)<f'(t)=\lambda+(r-1)t^{r-2}b_r(t^r)+rt^{2r-2}b_r'(t^r)\le Ct^{2r-2}.
\end{gather*}
In case \eqref{eq:caseB}, $f(t)\le\lambda t$ for $0<t\le\mu$, and by
concavity there are $A_1,A_2\in\R$ with $b_{r,\mu}(h)\le A_1+A_2h$ for
$h>0$. For $t>\mu$ let $h=(t-\mu)^r\le t^r$; then
$f(t)=\lambda t+h^\beta b_{r,\mu}(h)\le\lambda t+|A_1|t^{r-1}+|A_2|t^{2r-1}$,
and $t^{r-1}\le\mu^{-1}\max\{1,\mu\}\,t$ for $t\ge\mu$ since
$0\le r-1<1$. The bounds for $2r-1$ follow from $1\le r\le1/(1-s)$ and
$(1+s)(1-2s)\le1-s$.
\end{proof}

\begin{proof}[Proof of Lemma~\ref{lem:scalar}]
We have $f_K(0)=0$, $f_K'>0$ and $f_K(t)-tf_K'(t)=tf'(t)-f(t)>0$, so $f_K$ is
increasing and $f_K(t)/t$ is decreasing. Hence
$\frac\tau tf_K(t)\le f_K(\tau)\le f_K(t)$ for $0<\tau<t$, and integration in
$\tau$ proves \eqref{eq:scalar-bounds}. For $0<c\le1$, monotonicity
gives $f_K(ct)\le f_K(t)$; for $c\ge1$, the decrease of $f_K(t)/t$ gives
$f_K(ct)\le cf_K(t)$. Finally $F_K(ct)\le ct\,f_K(ct)\le c\max\{1,c\}\,tf_K(t)$.
\end{proof}

\begin{proof}[Proof of Lemma~\ref{lem:kernel}]
(i) Let $p_\tau$ be the density of the symmetric $2s$-stable law,
$\sqrt{2\pi}\,\widehat{p_\tau}(\omega)=e^{-\tau|\omega|^{2s}}$. Since
$y\mapsto e^{-\tau y^s}$ is completely monotone, Bernstein's theorem
shows that $p_\tau$ is a mixture of centered Gaussian densities, hence
positive and symmetric decreasing; moreover
$p_\tau(x)=\tau^{-1/(2s)}p_1(\tau^{-1/(2s)}x)$ and
$p_1(x)\le C(1+|x|)^{-1-2s}$ \cite{BG}, so that
$p_\tau(x)\le C\min\{\tau^{-1/(2s)},\tau|x|^{-1-2s}\}$. Define
$G_\lambda=\int_0^\infty e^{-\lambda\tau}p_\tau\dd\tau$. It has the
stated Fourier transform, is positive and symmetric decreasing, and
$\int_\R G_\lambda\dd x=\int_0^\infty e^{-\lambda\tau}\dd\tau=1/\lambda$ by
Tonelli's theorem. For $0<|x|\le1$ we split the $\tau$-integral at
$|x|^{2s}$ and use the two bounds on $p_\tau$:
\begin{equation*}
 G_\lambda(x)\le C\int_0^{|x|^{2s}}\tau|x|^{-1-2s}\dd\tau
 +C\int_{|x|^{2s}}^\infty\tau^{-1/(2s)}\dd\tau
 =\frac C2|x|^{2s-1}+\frac{C}{1/(2s)-1}|x|^{2s-1},
\end{equation*}
where $1/(2s)>1$ was used to integrate $\tau^{-1/(2s)}$ at infinity;
for $|x|\ge1$ we use the second bound on $p_\tau$ for all $\tau$:
$G_\lambda(x)\le C|x|^{-1-2s}\int_0^\infty\tau e^{-\lambda\tau}\dd\tau
=C\lambda^{-2}|x|^{-1-2s}$.
For \eqref{eq:kernel-modulus} with $0<h\le1$, the function
$y\mapsto G_\lambda(y+h)-G_\lambda(y)$ is nonnegative for $y<-h/2$,
nonpositive for $y>-h/2$, and has integral zero, so its $L^1$ norm is
$2\int_{-h/2}^{h/2}G_\lambda\dd y\le Ch^{2s}$; the case $h<0$
follows by the substitution $y\mapsto y-h$.

(ii) By scaling, $(-\Delta)^s\Psi_\sigma(x)=\sigma^{-1-4s}((-\Delta)^s\Psi_1)(x/\sigma)$,
so it suffices to prove $|(-\Delta)^s\Psi_1|\le C_\Psi\Psi_1$. Since
$s<1/2$ and $\Psi_1$ is smooth with bounded derivative,
$(-\Delta)^s\Psi_1(y)=c_{1,s}\int_\R(\Psi_1(y)-\Psi_1(z))|y-z|^{-1-2s}\dd z$
converges absolutely and is bounded, which handles $|y|\le2$, where
$\Psi_1$ is bounded below. Let $|y|\ge2$. If $|z-y|<|y|/2$, every point $w$
between $y$ and $z$ satisfies $|w|\ge|y|/2\ge1$, so
$|\Psi_1'(w)|=(1+2s)|w|(1+w^2)^{-(3+2s)/2}\le(1+2s)|w|^{-2-2s}\le C|y|^{-2-2s}$,
and the mean value theorem gives
$|\Psi_1(y)-\Psi_1(z)|\le C|y-z|\,|y|^{-2-2s}$; since
$\int_{|z-y|<|y|/2}|y-z|^{-2s}\dd z=\frac{2(|y|/2)^{1-2s}}{1-2s}$, this part of
the integral is at most $C|y|^{-1-4s}$. On $\{|z-y|\ge|y|/2\}$ we estimate
the two terms of $\Psi_1(y)-\Psi_1(z)$ separately:
\begin{gather*}
 \int_{|z-y|\ge|y|/2}\frac{\Psi_1(y)}{|y-z|^{1+2s}}\dd z=\Psi_1(y)\frac{(|y|/2)^{-2s}}s,\\
 \int_{|z-y|\ge|y|/2}\frac{\Psi_1(z)}{|y-z|^{1+2s}}\dd z\le\Bigl(\frac{|y|}2\Bigr)^{-1-2s}\int_\R\Psi_1\dd z.
\end{gather*}
Since $\Psi_1(y)\ge(2y^2)^{-(1+2s)/2}=2^{-(1+2s)/2}|y|^{-1-2s}$ for $|y|\ge1$
and $|y|^{-2s}\le1$, each of the three contributions is bounded by
$C\Psi_1(y)$, which proves the claim.
\end{proof}

\begin{proof}[Proof of Lemma~\ref{lem:solution-estimates}]
\emph{Step 1. Boundedness.}
In case \eqref{eq:caseA} with $s>1/2$, $\|u\|_\infty\le C\|u\|_{H^s}$ by
the Sobolev embedding. If $s=1/2$, then $r=2$, Lemma~\ref{lem:growth}
gives $|f(t)|\le C(1+t^3)$, and $u\in L^6(\R)$ by the critical
embedding; writing the equation as
$(-\Delta)^{1/2}u+u=f(u)+(1-\lambda)u$ with right side in $L^2(B_2(x_0))$
uniformly in $x_0$, the local estimate \cite[Lemma B.1]{FLS} (with
$N=1$, $p=2$, $q=\infty$) gives
$\|u\|_{L^\infty(B_1(x_0))}\le C(1+\|u\|_6^3+\|u\|_2)$ uniformly in
$x_0$.

Under \eqref{eq:small-growth}, $u\in L^2\cap L^{2/(1-2s)}$ by the
embedding, and $g=f(u)_+\le C_1(u+u^{p_0})$ belongs to $L^2$ because
$2p_0\le2/(1-2s)$; also $f(u)_-\le\max_{[0,1]}|f|+|f(1)|u$, since
$f(t)\ge f(1)t$ for $t\ge1$, so $f(u)\in L^2_{\rm loc}$ and the
distributional equation is meaningful. Testing with nonnegative $\psi\in H^s(\R)$
(approximated by a sequence $\{\psi_n\}_{n\ge1}$ of nonnegative functions in $C_c^\infty(\R)$) gives
$T_s(u,\psi)+\lambda\int_\R u\psi\dd x\le\int_\R g\psi\dd x$. The
function $z=G_\lambda*g\in H^{2s}(\R)$ satisfies the corresponding
equality by Plancherel's identity, and testing the difference with
$\psi=(u-z)_+$, together with $T_s(u-z,(u-z)_+)\ge T_s((u-z)_+)$, gives
$(u-z)_+=0$:
\begin{equation*}
 u\le G_\lambda*g\le C_1G_\lambda*(u+u^{p_0})\qquad\text{a.e. in }\R.
\end{equation*}
We now improve the integrability of $u$ step by step. Suppose that
$u\in L^2\cap L^q$ for some $q\ge2/(1-2s)$. Then $2\le q/p_0\le q$, as
$p_0\le1/(1-2s)$, so $u\in L^{q/p_0}$ by interpolation and
$u^{p_0}\in L^{q/p_0}$ because $u\in L^q$; thus $g_1=u+u^{p_0}\in L^2\cap L^{q/p_0}$.
By Lemma~\ref{lem:kernel}, $G_\lambda\in L^{q_1}$ whenever $0\le1-1/q_1<2s$.
If $p_0/q\ge2s$, we take $1-1/q_1=2s-\eta$ with a small $\eta>0$ and apply
Young's inequality
$\|G_\lambda*g_1\|_{\tilde q}\le\|G_\lambda\|_{q_1}\|g_1\|_{q/p_0}$, where
$1+1/\tilde q=1/q_1+p_0/q$, that is,
\begin{equation*}
 \frac1{\tilde q}=\frac{p_0}q-\Bigl(1-\frac1{q_1}\Bigr)=\frac{p_0}q-(2s-\eta)>0;
\end{equation*}
hence $u\le C_1G_\lambda*g_1\in L^{\tilde q}$, and $u\in L^2$ since $u\in H^s(\R)$. Using
$1/q\le(1-2s)/2$,
\begin{equation*}
 \frac1q-\frac1{\tilde q}=2s-\eta-\frac{p_0-1}q
 \ge2s-\eta-\frac{(p_0-1)(1-2s)}2=:\kappa_0,
\end{equation*}
and $\kappa_0>0$ for small $\eta$, because $(p_0-1)(1-2s)<4s$ is equivalent
to $p_0<(1+2s)/(1-2s)$, which holds since $p_0\le1/(1-2s)$. Thus $1/q$
decreases by at least the fixed amount $\kappa_0$ at each step, while
$\tilde q>q\ge2/(1-2s)$ keeps the hypothesis; after finitely many steps
$p_0/q<2s$. Then Hölder's inequality with $1-1/q_1=p_0/q$, that is,
$1/q_1+p_0/q=1$, gives $G_\lambda*g_1\in L^\infty$ with
$\|G_\lambda*g_1\|_\infty\le\|G_\lambda\|_{q_1}\|g_1\|_{q/p_0}$. Hence
$u\in L^\infty(\R)$.

\emph{Step 2. Truncated test functions.}
Fix $m_*\ge\|u\|_\infty$ and $K,f_K,F_K$ as in \eqref{eq:aux-definition}, and
for $n\ge1$ let $f_{K,n}(t)=\int_0^t\min\{f_K'(\tau),n\}\dd\tau$ on $[0,m_*]$,
extended by zero on $(-\infty,0)$ and by $f_{K,n}(m_*)$ on $(m_*,\infty)$. Then
$f_{K,n}$ is nondecreasing, Lipschitz with constant $n$, and
$0\le f_{K,n}\le\min\{f_K,nt\}$ on $[0,m_*]$, so $f_{K,n}(u)\in H^s(\R)$ with
$\|f_{K,n}(u)\|_2\le n\|u\|_2$ and $T_s(f_{K,n}(u))\le n^2T_s(u)$; moreover
$f_{K,n}(t)\le f_{K,n+1}(t)$ and $f_{K,n}(t)\to f_K(t)$ as $n\to\infty$, by monotone
convergence, since $\int_0^tf_K'(\tau)\dd\tau=f_K(t)$. Let $\chi\in C_c^\infty(\R)$ be symmetric decreasing, with
$0\le\chi\le1$, $\chi=1$ on $[-1,1]$ and $\chi=0$ outside $(-2,2)$, and define
$\chi_R(x)=\chi(x/R)$ for $R\ge1$. We record the cutoff approximation
used here and in Section~\ref{sec:coercivity}: the bound
$|\chi_R(x)-\chi_R(y)|\le C\min\{1,|x-y|/R\}$ gives
\begin{equation}\label{eq:cutoff-approximation}
 T_s(\chi_Rv)\le2T_s(v)+CR^{-2s}\|v\|_2^2
 \quad(v\in H^s(\R),\ R\ge1),\qquad
 \lim_{R\to\infty}\|\chi_Rv-v\|_{H^s}=0,
\end{equation}
the limit following from the uniform bound, the density of
$C_c^\infty(\R)$ in $H^s(\R)$ and $\chi_Rv_j=v_j$ for large $R$, where
$\{v_j\}_{j\ge1}\subset C_c^\infty(\R)$ converges to $v$ in $H^s(\R)$.

\emph{Step 3. Integrability of the nonlinear terms.}
Since $u$ and $f_K(u)$ are locally bounded, the function
$\chi_Rf_{K,n}(u)$ belongs to $H^s(\R)$ and is an admissible test function
in the equation
$(-\Delta)^su+f_K(u)=(K-\lambda)u$ (approximate it in $H^s$ by
$C_c^\infty$ functions supported in a fixed interval, on which
$\lambda u-f(u)$ is bounded):
\begin{equation*}
 T_s(u,\chi_Rf_{K,n}(u))+\int_\R\chi_Rf_K(u)f_{K,n}(u)\dd x
 =(K-\lambda)\int_\R\chi_Ru\,f_{K,n}(u)\dd x.
\end{equation*}
As $R\to\infty$, \eqref{eq:cutoff-approximation} gives
$T_s(u,\chi_Rf_{K,n}(u))\to T_s(u,f_{K,n}(u))$, dominated convergence applies
to the right side since $0\le\chi_Ru\,f_{K,n}(u)\le nu^2$, and monotone
convergence gives $f_K(u)f_{K,n}(u)\in L^1$ and
\begin{equation*}
 T_s(u,f_{K,n}(u))+\int_\R f_K(u)f_{K,n}(u)\dd x
 =(K-\lambda)\int_\R u\,f_{K,n}(u)\dd x.
\end{equation*}
Since $f_{K,n}$ is nondecreasing, $T_s(u,f_{K,n}(u))\ge0$, and with $f_{K,n}\le f_K$
and Cauchy--Schwarz,
$\|f_{K,n}(u)\|_2^2\le\int_\R f_K(u)f_{K,n}(u)\dd x\le(K-\lambda)\|u\|_2\|f_{K,n}(u)\|_2$.
Letting $n\to\infty$ by monotone convergence proves \eqref{eq:a-L2}.
By \eqref{eq:scalar-bounds}, $|f(t)|\le Kt+f_K(t)$,
$t|f'(t)|\le Kt+f_K(t)$ and $|F(t)|\le\frac K2t^2+tf_K(t)$ on $(0,m_*]$, so
$f(u),uf'(u)\in L^2$ and $F(u),uf(u),u^2f'(u)\in L^1$. Since
$f(u)\in L^2$, every term of the weak equation is continuous in the
$H^s$ norm of the test function, and density gives
\eqref{eq:weak-all-tests}.

\emph{Step 4. Global regularity and decay.}
By (ii), $f(u)-\lambda u\in L^2$, and the Fourier transform of the
equation, $|\omega|^{2s}\widehat u=\widehat{f(u)-\lambda u}$, holds
almost everywhere; by Plancherel, $u\in H^{2s}(\R)$. In case
\eqref{eq:caseA} we have $2s\ge1$, so $u\in H^1(\R)$, $u'\in L^2(\R)$ and
$|u(x)-u(y)|\le\|u'\|_2|x-y|^{1/2}$; moreover $\widehat u\in L^1(\R)$ by
Cauchy--Schwarz, and $u\in C_0(\R)$ by Fourier inversion and the
Riemann--Lebesgue lemma.

Under \eqref{eq:small-growth}, $f(u)\in L^\infty$ since $f$ is
continuous on $[0,\|u\|_\infty]$, and
$(|\omega|^{2s}+\lambda)\widehat u=\widehat{f(u)}$ gives
$u=G_\lambda*f(u)$. By \eqref{eq:kernel-modulus},
$|u(x+h)-u(x)|\le C\|f(u)\|_\infty|h|^{2s}$ for $|h|\le1$, so $u$
is uniformly continuous, hence in $C_0(\R)$ as $u\in L^2$. For the
decay rate, choose $\varepsilon_0,\delta>0$ with
$f(t)\le(\lambda-\delta)t$ for $0\le t\le\varepsilon_0$, and $R>0$ with
$u\le\varepsilon_0$ on $\{|x|>R\}$. Then \eqref{eq:weak-all-tests} gives
$T_s(u,\psi)+\delta\int_\R u\psi\dd x\le0$ for every nonnegative
$\psi\in H^s(\R)$ vanishing on $\{|x|\le R\}$. By
Lemma~\ref{lem:kernel}(ii) choose $\sigma>0$ with
$(-\Delta)^s\Psi_\sigma+\delta\Psi_\sigma\ge0$, and $A_0>0$ with
$A_0\Psi_\sigma\ge\|u\|_\infty$ on $\{|x|\le R\}$. For $\phi\in C_c^\infty(\R)$ the double integral obtained by inserting
the singular integral for $(-\Delta)^s\Psi_\sigma$ into
$\int_\R\phi\,(-\Delta)^s\Psi_\sigma\dd x$ converges absolutely, since
$\int_\R|\Psi_\sigma(y)-\Psi_\sigma(z)|\,|y-z|^{-1-2s}\dd z$ is bounded in
$y$ (proof of Lemma~\ref{lem:kernel}(ii)); Fubini's theorem and
symmetrization give
$\int_\R\phi\,(-\Delta)^s\Psi_\sigma\dd x=T_s(\Psi_\sigma,\phi)$. As
$\Psi_\sigma\in H^2(\R)$, the right side equals
$\int_\R\phi\,\mathcal F^{-1}(|\omega|^{2s}\widehat{\Psi_\sigma})\dd x$ with
$\mathcal F^{-1}(|\omega|^{2s}\widehat{\Psi_\sigma})\in L^2(\R)$, so the
pointwise $(-\Delta)^s\Psi_\sigma$ is this $L^2$ function, and by density
\begin{equation*}
 T_s(\Psi_\sigma,\psi)+\delta\int_\R\Psi_\sigma\psi\dd x
 =\int_\R\bigl((-\Delta)^s\Psi_\sigma+\delta\Psi_\sigma\bigr)\psi\dd x\ge0
 \qquad\text{for every nonnegative }\psi\in H^s(\R).
\end{equation*} Testing both inequalities with
$\psi=w_+$, $w=u-A_0\Psi_\sigma$, which vanishes on $\{|x|\le R\}$, and
using $T_s(w,w_+)\ge T_s(w_+)$ gives $w_+=0$. Thus
$u\le A_0\Psi_\sigma\le C(1+|x|)^{-1-2s}$.

\emph{Step 5. Positivity and smoothness.}
Since $u$ is continuous and positive almost everywhere, $u\ge0$
everywhere. By (i)--(iii), $u$ is a bounded continuous weak solution
with continuous right side $f(u)-\lambda u$, hence a viscosity solution
on every bounded interval by \cite[Theorem 1]{SV}. If $u(x_0)=0$, the
function equal to $0$ on $(x_0-\rho,x_0+\rho)$ and to $u$ elsewhere
touches $u$ from below at $x_0$, and the viscosity supersolution
inequality gives
$0=f(0)-\lambda u(x_0)\le-c_{1,s}\int_{|y-x_0|\ge\rho}u(y)|x_0-y|^{-1-2s}\dd y<0$,
a contradiction; so $u>0$ everywhere. On bounded open intervals
$I,J$ with $\overline I\subset J$ we then have $\min_{\overline J}u>0$, all derivatives of
$f$ are bounded on $[\min_{\overline J}u,\|u\|_\infty]$, and the interior
estimate in \cite{RS}
\begin{equation*}
 \|u\|_{C^{\alpha+2s}(I)}
 \le C\Bigl(\int_\R\frac{|u(x)|}{(1+|x|)^{1+2s}}\dd x
       +\|u\|_{C^\alpha(J)}+\|f(u)-\lambda u\|_{C^\alpha(J)}\Bigr)
 \qquad(\alpha,\alpha+2s\notin\mathbb N),
\end{equation*}
which extends to weak solutions by convolution as in the proof of
\cite[Lemma 2.9]{RS}, can be iterated starting from the Hölder
continuity obtained in Step~4 ($u\in C^{\alpha_0}_{\rm loc}$ for
$\alpha_0<1/2$ in case \eqref{eq:caseA}, for $\alpha_0<2s$ under
\eqref{eq:small-growth}), gaining more than $s/2$ at each step while
avoiding integer exponents. Thus $u\in C^\infty(\R)$.

\emph{Step 6. The maximum of $u$.}
Testing \eqref{eq:weak-all-tests} with $\psi=u$ gives
$T_s(u)+\lambda\|u\|_2^2=\int_\R uf(u)\dd x$. If $f(t)\le\lambda t$ for
$0<t\le\|u\|_\infty$, then $T_s(u)\le0$, so $u$ is constant, hence zero,
a contradiction. If $\mu>0$ and $f(\mu)\le\lambda\mu$, then
$f(t)\le\lambda t$ on $(0,\mu]$ because $f(t)/t$ is increasing, and the
first assertion gives $\|u\|_\infty>\mu$.
\end{proof}

The symmetry proof follows the method of the antisymmetric maximum
principles of Chen--Li--Li \cite{CLL} and uses the
Hopf lemma of Li--Chen \cite{LC}. These results are
linear statements for antisymmetric functions, and the reflected
equation \eqref{eq:reflected-equation} below is of this form, with a
coefficient bounded from below since $f'$ is bounded from above on
bounded subsets of $(0,\infty)$. Since the proof
needs only a negative-minimum estimate and a zero-point argument, we
give them directly, which also keeps the constants uniform in the
reflection point.

\begin{proof}[Proof of Lemma~\ref{lem:symmetry}]
\emph{Step 1. Negative minima of the reflected difference.}
By Lemma~\ref{lem:solution-estimates}, $u\in C^\infty(\R)\cap C_0(\R)$
and $u>0$. For $\ell\in\R$ let $x^\ell=2\ell-x$ and
$w_\ell(x)=u(x^\ell)-u(x)$. Subtracting the equations,
\begin{equation}\label{eq:reflected-equation}
 (-\Delta)^sw_\ell+c_\ell(x)w_\ell=0,\qquad
 c_\ell(x)=\lambda-\int_0^1f'((1-t)u(x)+tu(x^\ell))\dd t.
\end{equation}
The function $w_\ell$ is smooth, antisymmetric about $\ell$, and
vanishes at infinity, and
$\int_\R|w_\ell(x)|\,(1+|x|)^{-1-2s}\dd x\le C\|u\|_2$ by
Cauchy--Schwarz. By
\eqref{eq:f-assumptions}, $c_\ell\ge-C_f$ with
$C_f=\sup_{0<t\le\|u\|_\infty}f'(t)-\lambda$, which does not depend on
$\ell$ and is positive by Lemma~\ref{lem:solution-estimates}(v) and
\eqref{eq:f-assumptions}. Choose $t_0,\nu,R>0$ with $f'(t)\le\lambda-\nu$ for $0<t\le t_0$
and $u\le t_0$ on $\{|x|>R\}$; then $c_\ell(x)\ge\nu$ at every point with
$w_\ell(x)<0$ and $|x|>R$, because there $0<u(x^\ell)<u(x)\le t_0$.

Suppose $w_\ell<0$ somewhere on $(-\infty,\ell)$. Then $w_\ell$ attains
a negative minimum at some $x_\ell<\ell$. We split the singular integral
$(-\Delta)^sw_\ell(x_\ell)=c_{1,s}\int_\R(w_\ell(x_\ell)-w_\ell(y))|x_\ell-y|^{-1-2s}\dd y$
into the two half-lines. In the integral over $(\ell,\infty)$ we
substitute $y=z^\ell=2\ell-z$ with $z<\ell$ and use
$w_\ell(z^\ell)=-w_\ell(z)$:
\begin{equation*}
 \int_\ell^\infty\frac{w_\ell(x_\ell)-w_\ell(y)}{|x_\ell-y|^{1+2s}}\dd y
 =\int_{-\infty}^\ell\frac{w_\ell(x_\ell)+w_\ell(z)}{|x_\ell-z^\ell|^{1+2s}}\dd z
 =\int_{-\infty}^\ell\frac{-(w_\ell(x_\ell)-w_\ell(z))+2w_\ell(x_\ell)}{|x_\ell-z^\ell|^{1+2s}}\dd z.
\end{equation*}
Adding the integral over $(-\infty,\ell)$ and renaming $z$ as $y$,
\begin{equation*}
 \begin{split}
 (-\Delta)^sw_\ell(x_\ell)
 ={}&c_{1,s}\int_{-\infty}^\ell
  (w_\ell(x_\ell)-w_\ell(y))
  \Bigl(\frac1{|x_\ell-y|^{1+2s}}-\frac1{|x_\ell-y^\ell|^{1+2s}}\Bigr)\dd y\\
 &+2c_{1,s}w_\ell(x_\ell)\int_{-\infty}^\ell\frac{\dd y}{(2\ell-x_\ell-y)^{1+2s}},
 \end{split}
\end{equation*}
both integrals being convergent since $w_\ell'(x_\ell)=0$. The first
integrand is nonpositive, because $w_\ell(x_\ell)\le w_\ell(y)$ and, for
$y<\ell$, $|x_\ell-y^\ell|=(\ell-x_\ell)+(\ell-y)>|(\ell-x_\ell)-(\ell-y)|=|x_\ell-y|$;
and since $2\ell-x_\ell-y>\ell-x_\ell>0$ for $y<\ell$, the substitution
$z=2\ell-x_\ell-y$ gives
$\int_{-\infty}^\ell(2\ell-x_\ell-y)^{-1-2s}\dd y=\int_{\ell-x_\ell}^\infty z^{-1-2s}\dd z
=(2s(\ell-x_\ell)^{2s})^{-1}$. Hence, by \eqref{eq:reflected-equation},
\begin{equation}\label{eq:negative-minimum-estimate}
 c_\ell(x_\ell)=\frac{(-\Delta)^sw_\ell(x_\ell)}{-w_\ell(x_\ell)}
 \le-\frac{c_{1,s}}{s(\ell-x_\ell)^{2s}}<0.
\end{equation}
Thus $|x_\ell|\le R$, because $c_\ell(x)\ge\nu$ at every point with
$w_\ell(x)<0$ and $|x|>R$, and
$c_{1,s}/(s(\ell-x_\ell)^{2s})\le C_f$. With
$d_0=(c_{1,s}/(sC_f))^{1/(2s)}$,
\begin{equation}\label{eq:negative-minimum-location}
 x_\ell\in[-R,R],\qquad \ell-x_\ell\ge d_0.
\end{equation}

\emph{Step 2. Symmetry and strict monotonicity.}
By \eqref{eq:negative-minimum-location}, $w_\ell\ge0$ on
$(-\infty,\ell)$ for every $\ell\le-R$. Let $\ell_*$ be the supremum of
the $\ell$ such that $w_b\ge0$ on $(-\infty,b)$ for every $b\le\ell$;
it is finite, since otherwise $u$ would be nondecreasing, contradicting
$u>0$ and $u\in C_0(\R)$. By continuity $w_{\ell_*}\ge0$ on
$(-\infty,\ell_*)$. At an interior zero $x$ of a nonnegative
antisymmetric solution $w_b$ of \eqref{eq:reflected-equation}, pairing
the half-lines gives
$0=-c_{1,s}\int_{y<b}w_b(y)(|x-y|^{-1-2s}-|x-y^b|^{-1-2s})\dd y$ with a
strictly positive kernel, so $w_b\equiv0$; hence $w_b>0$ or $w_b\equiv0$
on $(-\infty,b)$. If $w_{\ell_*}>0$, the definition of $\ell_*$ gives a sequence
$\{\ell_j\}_{j\ge1}$ with $\ell_j>\ell_*$ and $\ell_j\to\ell_*$ such that $w_{\ell_j}$ has a
negative minimum at some $x_j$; by \eqref{eq:negative-minimum-location}, a
subsequence of $\{x_j\}_{j\ge1}$ converges to some $x_*\in[-R,R]$ with $x_*\le\ell_*-d_0$, and
continuity gives $0<w_{\ell_*}(x_*)=\lim w_{\ell_j}(x_j)\le0$, a
contradiction. Hence $w_{\ell_*}\equiv0$, that is, $u$ is symmetric about
$\ell_*$. For $b<\ell_*$ we cannot have $w_b\equiv0$, since then
$u(x+2(\ell_*-b))=u(x)$ and $u\in C_0(\R)$ would force $u\equiv0$; so
$w_b>0$ on $(-\infty,b)$, and taking $b=(x+y)/2$ for $x<y\le\ell_*$
gives $u(y)-u(x)=w_b(x)>0$. Translate $\ell_*$ to zero.

\emph{Step 3. Strict negativity of the derivative.}
Let $z=-u'$; it is smooth, odd, nonnegative and not identically zero
on $(0,\infty)$, and $\|z\|_1=2u(0)$. Differentiating the equation,
\begin{equation}\label{eq:derivative-equation}
 (-\Delta)^sz+(\lambda-f'(u))z=0,
\end{equation}
which holds pointwise since $z$ is smooth with
$\int_\R|z|(1+|x|)^{-1-2s}\dd x<\infty$. If $z(x_0)=0$ for some $x_0>0$,
then $z'(x_0)=0$, and pairing the half-lines gives the absolutely
convergent identity
$0=(-\Delta)^sz(x_0)=-c_{1,s}\int_0^\infty z(y)(|x_0-y|^{-1-2s}-(x_0+y)^{-1-2s})\dd y<0$,
a contradiction. Thus $u'<0$ on $(0,\infty)$.

\emph{Step 4. Nondegeneracy of the maximum.}
Now $z>0$ on $(0,\infty)$, and since $u(0)>0$ the coefficient of
\eqref{eq:derivative-equation} is bounded near zero, bounded from below
on $(0,\infty)$ by $-C_f$ (as in Step~1, since $0<u\le\|u\|_\infty$), and
satisfies $x^2(\lambda-f'(u(x)))\to0$ as $x\to0^+$. Together with the smoothness,
oddness and tail bound of Step~3, this verifies the hypotheses of the
antisymmetric Hopf lemma \cite{LC} on the half-line
$(0,\infty)$, which gives $-u''(0)=z'(0)>0$.
\end{proof}

% ---- Acknowledgments and statements (to be completed before submission) ----
%\section*{Acknowledgments}
%The authors thank ... . This work was supported by ... (grant numbers).
%
%\subsection*{Data availability statement}
%Data sharing is not applicable to this article as no data sets were
%generated or analysed during the current study.
%
%\subsection*{Conflict of interest}
%The authors declare that they have no conflict of interest.
% ----------------------------------------------------------------------------

\raggedbottom


\begin{thebibliography}{99}
\raggedright

\bibitem{AT}
C.~J. Amick and J.~F. Toland,
\emph{Uniqueness and related analytic properties for the
Benjamin--Ono equation--a nonlinear Neumann problem in the plane},
Acta Math. \textbf{167} (1991), 107--126.

\bibitem{APZ}
X. An, S. Peng, and F. Zhong,
\emph{Quantitative analysis of ground states for the fractional
logarithmic Schr\"odinger equation},
preprint (2026).
\href{https://arxiv.org/abs/2609.10254}{arXiv:2609.10254}.

\bibitem{Ar}
A.~H. Ardila,
\emph{Existence and stability of standing waves for nonlinear fractional
Schr\"odinger equation with logarithmic nonlinearity},
Nonlinear Anal. \textbf{155} (2017), 52--64.

\bibitem{Ben}
T.~B. Benjamin,
\emph{Internal waves of permanent form in fluids of great depth},
J. Fluid Mech. \textbf{29} (1967), 559--592.

\bibitem{BDG}
M. Bhakta, P. Das, and D. Ganguly,
\emph{Fractional Schr\"odinger equations with mixed nonlinearities:
asymptotic profiles, uniqueness and nondegeneracy of ground states},
J. Geom. Anal. \textbf{35} (2025), no.~8, Paper No.~245.

\bibitem{BBM}
I. Bia\l ynicki-Birula and J. Mycielski,
\emph{Nonlinear wave mechanics},
Ann. Physics \textbf{100} (1976), 62--93.

\bibitem{BG}
R.~M. Blumenthal and R.~K. Getoor,
\emph{Some theorems on stable processes},
Trans. Amer. Math. Soc. \textbf{95} (1960), 263--273.

\bibitem{BF}
L. Brasco and G. Franzina,
\emph{Convexity properties of Dirichlet integrals and Picone-type
inequalities},
Kodai Math. J. \textbf{37} (2014), 769--799.

\bibitem{BO}
H. Brezis and L. Oswald,
\emph{Remarks on sublinear elliptic equations},
Nonlinear Anal. \textbf{10} (1986), 55--64.

\bibitem{CS}
L. Caffarelli and L. Silvestre,
\emph{An extension problem related to the fractional Laplacian},
Comm. Partial Differential Equations \textbf{32} (2007), 1245--1260.

\bibitem{CL}
T. Cazenave and P.-L. Lions,
\emph{Orbital stability of standing waves for some nonlinear
Schr\"odinger equations},
Comm. Math. Phys. \textbf{85} (1982), no.~4, 549--561.

\bibitem{CLL}
W. Chen, C. Li, and Y. Li,
\emph{A direct method of moving planes for the fractional Laplacian},
Adv. Math. \textbf{308} (2017), 404--437.

\bibitem{Cof}
C.~V. Coffman,
\emph{Uniqueness of the ground state solution for $\Delta u-u+u^3=0$ and
a variational characterization of other solutions},
Arch. Rational Mech. Anal. \textbf{46} (1972), 81--95.

\bibitem{DSZ}
P. d'Avenia, M. Squassina, and M. Zenari,
\emph{Fractional logarithmic Schr\"odinger equations},
Math. Methods Appl. Sci. \textbf{38} (2015), no.~18, 5207--5216.

\bibitem{DdPW}
J. D\'avila, M. del Pino, and J. Wei,
\emph{Concentrating standing waves for the fractional nonlinear
Schr\"odinger equation},
J. Differential Equations \textbf{256} (2014), no.~2, 858--892.

\bibitem{DV}
M.~G. Delgadino and M. Vaughan,
\emph{Continuous symmetrizations and uniqueness of solutions to nonlocal equations},
Anal. PDE \textbf{18} (2025), no.~10, 2325--2366.

\bibitem{DYY}
M.~G. Delgadino, X. Yan, and Y. Yao,
\emph{Uniqueness and nonuniqueness of steady states of
aggregation-diffusion equations},
Comm. Pure Appl. Math. \textbf{75} (2022), no.~1, 3--59.

\bibitem{DSa}
J.~I. D\'iaz and J.~E. Sa\'a,
\emph{Existence et unicit\'e de solutions positives pour certaines
\'equations elliptiques quasilin\'eaires},
C. R. Acad. Sci. Paris S\'er. I Math. \textbf{305} (1987), 521--524.

\bibitem{DNPV}
E. Di Nezza, G. Palatucci, and E. Valdinoci,
\emph{Hitchhiker's guide to the fractional Sobolev spaces},
Bull. Sci. Math. \textbf{136} (2012), no.~5, 521--573.

\bibitem{FV}
M.~M. Fall and E. Valdinoci,
\emph{Uniqueness and nondegeneracy of positive solutions of
$(-\Delta)^su+u=u^p$ in $\mathbb R^N$ when $s$ is close to $1$},
Comm. Math. Phys. \textbf{329} (2014), no.~1, 383--404.

\bibitem{FW}
M.~M. Fall and T. Weth,
\emph{Nondegeneracy properties and uniqueness of positive solutions to a
class of fractional semilinear equations},
Ann. Inst. H. Poincar\'e C Anal. Non Lin\'eaire \textbf{43} (2026),
no.~6, 1259--1298.

\bibitem{FQT}
P. Felmer, A. Quaas, and J. Tan,
\emph{Positive solutions of the nonlinear Schr\"odinger equation with the
fractional Laplacian},
Proc. Roy. Soc. Edinburgh Sect. A \textbf{142} (2012), 1237--1262.

\bibitem{Fr}
R.~L. Frank,
\emph{On the uniqueness of ground states of non-local equations},
Journ\'ees \'Equations aux D\'eriv\'ees Partielles (2011),
Exp. No.~5, 10~pp.

\bibitem{FL}
R.~L. Frank and E. Lenzmann,
\emph{Uniqueness of non-linear ground states for fractional Laplacians in $\mathbb R$},
Acta Math. \textbf{210} (2013), no.~2, 261--318.

\bibitem{FLS}
R.~L. Frank, E. Lenzmann, and L. Silvestre,
\emph{Uniqueness of radial solutions for the fractional Laplacian},
Comm. Pure Appl. Math. \textbf{69} (2016), no.~9, 1671--1726.

\bibitem{GNN}
B. Gidas, W.-M. Ni, and L. Nirenberg,
\emph{Symmetry of positive solutions of nonlinear elliptic equations in
$\mathbb R^n$},
in: Mathematical Analysis and Applications, Part A,
Adv. in Math. Suppl. Stud. \textbf{7a}, Academic Press, New York, 1981,
369--402.

\bibitem{GSS}
M. Grillakis, J. Shatah, and W. Strauss,
\emph{Stability theory of solitary waves in the presence of symmetry.~I},
J. Funct. Anal. \textbf{74} (1987), no.~1, 160--197.

\bibitem{Je}
L. Jeanjean,
\emph{Existence of solutions with prescribed norm for semilinear elliptic
equations},
Nonlinear Anal. \textbf{28} (1997), 1633--1659.

\bibitem{JZZ2}
Y. Jiang, X. Zeng, and H. Zhou,
\emph{Uniqueness and nondegeneracy of positive solutions for a class of fractional
Schr\"odinger equations in higher
dimensions}, in preparation.

\bibitem{Kato}
T. Kato,
\emph{Perturbation Theory for Linear Operators},
second edition, Grundlehren der mathematischen Wissenschaften,
Vol.~132, Springer-Verlag, Berlin--New York, 1976.

\bibitem{KMR}
C.~E. Kenig, Y. Martel, and L. Robbiano,
\emph{Local well-posedness and blow-up in the energy space for a class
of $L^2$ critical dispersion generalized Benjamin--Ono equations},
Ann. Inst. H. Poincar\'e C Anal. Non Lin\'eaire \textbf{28} (2011),
no.~6, 853--887.

\bibitem{KLT}
P. Kfoury, S. Le Coz, and T.-P. Tsai,
\emph{Analysis of stability and instability for standing waves of the
double power one dimensional nonlinear Schr\"odinger equation},
C. R. Math. Acad. Sci. Paris \textbf{360} (2022), 867--892.

\bibitem{Kw}
M.~K. Kwong,
\emph{Uniqueness of positive solutions of $\Delta u-u+u^p=0$ in
$\mathbb R^n$},
Arch. Rational Mech. Anal. \textbf{105} (1989), 243--266.

\bibitem{Las00}
N. Laskin,
\emph{Fractional quantum mechanics and L\'evy path integrals},
Phys. Lett. A \textbf{268} (2000), 298--305.

\bibitem{Las02}
N. Laskin,
\emph{Fractional Schr\"odinger equation},
Phys. Rev. E \textbf{66} (2002), 056108.

\bibitem{LC}
C. Li and W. Chen,
\emph{A Hopf type lemma for fractional equations},
Proc. Amer. Math. Soc. \textbf{147} (2019), no.~4, 1565--1575.

\bibitem{Li77}
E.~H. Lieb,
\emph{Existence and uniqueness of the minimizing solution of Choquard's
nonlinear equation},
Stud. Appl. Math. \textbf{57} (1977), 93--105.

\bibitem{LZ}
H. Luo and Z. Zhang,
\emph{Normalized solutions to the fractional Schr\"odinger equations
with combined nonlinearities},
Calc. Var. Partial Differential Equations \textbf{59} (2020),
Paper No.~143.

\bibitem{MZ}
L. Ma and L. Zhao,
\emph{Classification of positive solitary solutions of the nonlinear
Choquard equation},
Arch. Ration. Mech. Anal. \textbf{195} (2010), 455--467.

\bibitem{McL}
K. McLeod,
\emph{Uniqueness of positive radial solutions of $\Delta u+f(u)=0$ in
$\mathbb R^n$.~II},
Trans. Amer. Math. Soc. \textbf{339} (1993), no.~2, 495--505.

\bibitem{MVS}
V. Moroz and J. Van Schaftingen,
\emph{A guide to the Choquard equation},
J. Fixed Point Theory Appl. \textbf{19} (2017), no.~1, 773--813.

\bibitem{Ono}
H. Ono,
\emph{Algebraic solitary waves in stratified fluids},
J. Phys. Soc. Japan \textbf{39} (1975), 1082--1091.

\bibitem{RS}
X. Ros-Oton and J. Serra,
\emph{The Dirichlet problem for the fractional Laplacian:
regularity up to the boundary},
J. Math. Pures Appl. \textbf{101} (2014), no.~3, 275--302.

\bibitem{ST}
J. Serrin and M. Tang,
\emph{Uniqueness of ground states for quasilinear elliptic equations},
Indiana Univ. Math. J. \textbf{49} (2000), 897--923.

\bibitem{SV}
R. Servadei and E. Valdinoci,
\emph{Weak and viscosity solutions of the fractional Laplace equation},
Publ. Mat. \textbf{58} (2014), no.~1, 133--154.

\bibitem{So}
N. Soave,
\emph{Normalized ground states for the NLS equation with combined
nonlinearities},
J. Differential Equations \textbf{269} (2020), no.~9, 6941--6987.

\bibitem{SZ}
L. Song and W. Zou,
\emph{Uniqueness of positive solutions of double-power nonlinear
stationary Schr\"odinger equations: the fixed frequency case and the
fixed mass case},
preprint (2026).
\href{https://arxiv.org/abs/2609.15882}{arXiv:2609.15882}.


\bibitem{Wei87}
M.~I. Weinstein,
\emph{Existence and dynamic stability of solitary wave solutions of
equations arising in long wave propagation},
Comm. Partial Differential Equations \textbf{12} (1987), 1133--1173.

\bibitem{Xi}
C.-L. Xiang,
\emph{Uniqueness and nondegeneracy of ground states for Choquard
equations in three dimensions},
Calc. Var. Partial Differential Equations \textbf{55} (2016), no.~6,
Paper No.~134.

\end{thebibliography}
\end{document}